\documentclass[preprint,3p,12pt]{elsarticle}

\usepackage[utf8]{inputenc}
\usepackage[T1]{fontenc}
\usepackage{calc}
\usepackage{subfigure}
\usepackage{titlesec}
\usepackage{graphicx}
\usepackage{mathtools}
\usepackage{xcolor}
\usepackage{multirow}
\usepackage{comment}
\usepackage{ulem}

\usepackage[linesnumbered,ruled,vlined]{algorithm2e}
\usepackage{amssymb}
\newcommand{\X}{\textbf{X}}
\newcommand{\ud}{\textbf{u}}
\newcommand{\eps}{\boldsymbol{\varepsilon}}
\newcommand{\epse}{\boldsymbol{{\varepsilon}}_{e}}
\newcommand{\epsp}{\boldsymbol{{\varepsilon}}_{p}}

\newcommand{\sig}{\boldsymbol{\sigma}}

\newcommand{\Fm}{\textbf{F}}

\journal{Elsevier}
\SetKwInput{KwInput}{Input}                % Set the Input
\SetKwInput{KwOutput}{Output}              % set the Output

\graphicspath{{img/}}

\begin{document}

\begin{frontmatter}


  \title{Extended tensor decomposition model reduction methods: training, prediction, and design under uncertainty}


  \author[label2]{Ye Lu\corref{cor1}}
  \cortext[cor1]{Corresponding author}
%  \ead{yeluu90@gmail.com}
  \ead{yelu@umbc.edu}
  \author[label3]{Satyajit Mojumder}
  \author[label3]{Jiachen Guo}
  \author[label1]{Yangfan Li}
  \author[label1]{Wing Kam Liu\corref{cor1}}
  \ead{w-liu@northwestern.edu}

 \address[label2]{Department of Mechanical Engineering, University of Maryland Baltimore County, Baltimore, USA}
 \address[label3]{Theoretical and Applied Mechanics Program, Northwestern
University, Evanston, USA}
  \address[label1]{Department of Mechanical Engineering, Northwestern
University, Evanston, USA}

%\address{California, United States}

  \begin{abstract}
  This paper introduces an extended tensor decomposition (XTD) method for model reduction. The proposed method is based on a sparse non-separated enrichment to the conventional tensor decomposition, which is expected to improve the approximation accuracy and the reducibility (compressibility) in highly nonlinear and singular cases. The proposed XTD method can be a powerful tool for solving nonlinear space-time parametric problems. The method has been successfully applied to parametric elastic-plastic problems and real time additive manufacturing residual stress predictions with uncertainty quantification. Furthermore, a combined XTD-SCA (self-consistent clustering analysis) strategy {is} presented for multi-scale material modeling, which enables real time multi-scale multi-parametric simulations. The efficiency of the method is demonstrated with comparison to finite element analysis. The proposed method enables a novel framework for fast manufacturing and material design with uncertainties.  
  
  \end{abstract}


\end{frontmatter}

%\linenumbers

\newcommand{\revision}[2]{\sout{#1} \textcolor{red}{(#2)}}

\section{Introduction}
Simulation-based engineering and science, including design and uncertainty quantification, usually necessitate fast responses of numerical models. Despite the significant improvement of the computer hardware over the last decades, real-time simulations of large-scale systems are still intractable with conventional finite element analysis (FEA). The repetitive evaluation of the numerical model using FEA makes the design process computationally prohibitive. Furthermore,  this kind of problem can involve solutions in a large-scale space-time-parameter domain, which challenges the metamodeling or surrogate modeling approaches \cite{wang2007review} that succeeded in relatively small size problems. For example, the design of additive manufacturing processes can involve a large number of parameters related to the process and materials, such as the input heat source power, laser beam size, scan velocity,  material absorptivity, powder size and properties. Optimizing the design parameters demands space-time simulations of the process for different parameter combinations and thus the space-time-parametric solutions. Due to the computational expense of the space-time process simulation, the  metamodeling methods for this design problem require the development of powerful sampling methods for the space-time-parameter domain, which seems still an open question.

Exploring efficient space-time-parametric solution strategies has been one major task in the active area of model order reduction (or model reduction  for short). The model reduction was originated for the purpose of reducing the computational cost of numerical simulations in physical space-time domains, by reducing the degrees of freedom of underlying numerical systems. Self-consistent clustering analysis (SCA) \cite{liu2016self,yu2021multiresolution}, proper orthogonal decomposition (POD) \cite{willcox2002balanced,goury2016automatised,kerfriden2011bridging,lu2018space}, hyper reduction \cite{ryckelynck2009hyper,carlberg2013gnat,zhang2017efficient,lu2020adaptive}, and space-time LATIN proper generalized decomposition (PGD) \cite{ladeveze2010latin,boucinha2014ideal,giacoma2016efficient} are some of the most successful approaches in terms of space-time computation acceleration. For accelerating the computation of parametric solutions, the PGD method has been extended to consider the parameters as extra-coordinates of solutions \cite{ammar2006new,chinesta2011short,chinesta2013pgd} and decompose them like space and time in the approximation. This leads to a new definition of PGD and provides an attractive way to compute space-time-parametric solutions. Specifically, the basic idea is to decompose the multidimensional space-time-parametric solution into several separated 1D functions  and thus reduce the complexity of solving the multi-dimensional problem. However, it seems difficult to apply this method to nonlinear  problems, like plasticity and additive manufacturing. Differently, the so-called LATIN-PGD \cite{relun2013model,heyberger2012multiparametric,ladeveze2018extended} has also been extended for parametric studies, which relies on the LATIN solver \cite{ladeveze2010latin} with a point-wise evaluation of parameters. The method has shown the efficiency for some parametric elastic-plastic problems, thanks to the non-incremental nature of the space-time LATIN solver. In the context of additive manufacturing or other welding-like processes, the nonlinear parametric model reduction becomes even more challenging, especially in terms of problem reducibility. This has motivated the development of the non-intrusive data-driven PGD methods, i.e. high order PGD (HOPGD) \cite{lu2018multi,modesto2015proper} and its sparse counterpart \cite{lu2018adaptive,lu2019datadriven}. The HOPGD can be used for data learning (regression) and compression and has a direct connection to machine learning techniques. 

Machine learning has shown its superior performance in regression, especially for  complex discrete data. The universal approximation theorem \cite{hornik1989multilayer} and the natural scalability of deep learning neural networks have motivated the application to the function approximation for large systems. Many attempts were made to construct  machine learning based surrogate models for fast  space-time-parametric predictions (see e.g. \cite{lu2019data}). However, purely data-driven machine learning needs a huge amount of data that sometimes are not easy to obtain. To alleviate the need of data, the physics-informed type neural networks \cite{raissi2019physics} {was} proposed by incorporating the physically based governing equations into  the machine learning loss functions. \cite{liu2019deep} proposed to train a deep material network for nonlinear response predictions by using only elastic response data.  More recently, the Hierarchical Deep-learning Neural Network (HiDeNN) method \cite{zhang2021hierarchical,saha2021hierarchical}
 {was} developed by constraining the weights and biases of DNN to mesh coordinates to build the shape functions for FEA. The method allowed the automatic mesh adaptivity and showed good potential to prevent large numerical systems. Nevertheless, the machine learning-based methods {still suffer from significant computational cost}. A solution to reduce the cost is combining the machine learning with model reduction methods. As demonstrated by \cite{zhang2021hidenn}, the coupled HiDeNN-PGD method shows a significant improvement in terms of both efficiency and accuracy by taking advantages of the HiDeNN and PGD methods, compared to conventional FEA. From a general point of view, the synergy of machine learning and model reduction can give rise to a new type of reduced-order machine learning method and opens very interesting opportunities in computational science and engineering. The development of efficient nonlinear parametric model reduction methods is crucial to enable a wider application of such methods.

This work proposes a model reduction method, the so-called extended tensor decomposition (XTD), that provides a new way to solve nonlinear space-time-parametric problems. The XTD stands for a tensor decomposition with a sparse non-separated enrichment. The motivation is that the conventional tensor decomposition may require a strong assumption on the smoothness or the separability of the approximated functions, which can be violated in highly nonlinear or singular cases, thus the enrichment is expected to improve the approximation accuracy and the reducibility. The sparsity of the enrichment is to ensure the localization of the enrichment and minimize the associated computational cost.   Like other model reduction methods, the XTD can be computed from either numerical models or data. If the data is provided, the XTD can be used as a regression (or compression) method. Otherwise, the XTD can perform a reduced-order modeling with the governing equations. In any cases, the overall XTD framework follows a general offline-online strategy, in which the training is performed offline for a large space-time-parameter domain so that the online prediction is real time. The proposed method has been applied to solve nonlinear elastic-plastic problems by performing reduced-order modeling. Different challenging applications, including additive manufacturing type thermal residual stress and multi-scale composite modeling enabled by a coupled XTD-SCA approach, will be presented. It is found that the XTD method shows a good performance in terms of efficiency and accuracy, especially for large scale problems. 

This paper is organized as follows. Section 2 presents the overall offline-online XTD framework. Section 3 presents the formulation of the problem and the XTD solution strategy. Section 4 presents the numerical examples concerning additive manufacturing and XTD-SCA based multi-scale material modeling. Finally, the paper closes with some concluding remarks.

\section{General offline-online XTD learning framework}
\subsection{Extended tensor decomposition}
Let us consider a general  multidimensional space-time-parametric function. 
\begin{equation}
\displaystyle
\label{eq:XTD0} \ud(\X, t,{\mu_1},\dots,\mu_n)
\end{equation}
where $\X$, $t$ are the usual space and time variables. $\mu_i|_{i=1,\dots,n}$ denote general parameters for a given problem. The function $\ud$ can represent any quantity of interest of engineering problems, which is varying in a space-time-parameter domain. Now, we can consider a discrete representation of the above function using a high order array,
\begin{equation}
\displaystyle
\label{eq:XTDdisc} \boldsymbol{\mathcal{U}}\in \mathbb{R}^{{N_\X}\times{N_t}\times{N_{\mu_1}}\times\cdots\times{N_{\mu_n}}}
\end{equation}
with $N_\X$ denotes the total number degrees of freedom in the discretized physical domain, $N_{t}$ is the number of time steps, 
$N_{\mu_i}$ is the number of discrete points for the parameter $\mu_i$. This array can also be considered as a tensor of order $2+n$ (see e.g. \cite{kolda2009tensor}). Then, we define the following decomposition as the eXtended Tensor Decomposition (XTD).
\begin{equation}
\displaystyle
\label{eq:XTDdef} \boldsymbol{\mathcal{U}}\approx\sum_{m=1}^{M}
\boldsymbol{\mathcal{A}}^{(m)}\otimes \boldsymbol{T}^{(m)}\otimes\boldsymbol{g}_1^{(m)}\otimes\cdots\otimes\boldsymbol{g}_n^{(m)}+\sum_{k=1}^{K}\Tilde{\boldsymbol{\mathcal{U}}}^{(k)}
\end{equation}
where $\boldsymbol{\mathcal{A}}^{(m)}\in\mathbb{R}^{N_\X}$, $\boldsymbol{T}^{(m)}\in\mathbb{R}^{N_t}$, $\boldsymbol{g}_i^{(m)}\in\mathbb{R}^{N_{\mu_i}}$ are vectors, $\otimes$ denotes the vector outer product. $\Tilde{\boldsymbol{\mathcal{U}}}^{(k)}$ is a sparse array which has the same order as the original $\boldsymbol{\mathcal{U}}$. $M$ and $K$ are the numbers of modes, which should be relatively small compared to the total size of the problem. It is known that the first part of equation \eqref{eq:XTDdef} is called canonical tensor decomposition (Generally, it can be any other tensor decomposition format, such as  Tucker  decomposition {\cite{tucker1966some,kiers2000towards}}). The $\Tilde{\boldsymbol{\mathcal{U}}}^{(k)}$ appear as an enrichment, which is expected to remove the potential singularity or the high non-linearity from the original function (i.e. multidimensional array). {More specifically, the singularity is defined as potential  discontinuity or non-smoothness of the original function $\boldsymbol{\mathcal{U}}$, which might appear in both physical and parameter space, from a general point of view. In physical (space-time)  space, the singularity can be due to fracture or  bi-material interfaces, which causes the spatial displacement discontinuity. In  parameter space, the singularity can be the sharp changes of solutions with respect to parameter changes, which can be caused by nonlinear material phase changes. For example, let us consider a parametric thermal problem with the material melting temperature as a varying parameter, whose value determines the phase (solid or liquid) of the material system. The parametric solutions of this problem  may have a sharp change when the value of melting temperature varies from one to another, since one material system with a higher melting temperature may be pure solid while the other one with a lower melting point is pure liquid. Similar situations can happen for elastoplastic problems with a varying yield stress, since the solution with a higher yield stress can be purely elastic while the other one with a lower elastic limit may include nonlinear plastic responses. These examples of singularity in physical or parameter space are expected to appear locally in the global space-time-parameter domain.}  Therefore, the enrichment array should be sparse. We expect that this enrichment can help improve the reducibility of problem and the overall approximation accuracy of tensor decomposition, especially when singularity appears.

Come back to the continuous representation, the XTD can be written as \begin{equation}
\displaystyle
\label{eq:XTDdefcontinu}
\ud\approx \underbrace{\sum_{m=1}^{M}\boldsymbol{a}^{(m)}(\X)T^{(m)}(t)g_1^{(m)}(\mu_1)\cdots g_n^{(m)}(\mu_n) + \sum_{k=1}^{K}\Tilde{\boldsymbol{u}}^{(k)}(\X,t, {\mu_1},\dots,\mu_n)}_{\ud^{\text{XTD}}}
\end{equation}
where $\boldsymbol{a}^{(m)} (\X),T^{(m)}(t), g_i^{(m)}(\mu_i)$ are separated functions. The enrichment  $\Tilde{\boldsymbol{u}}^{(k)}(\X, t,{\mu_1},\dots,\mu_n)$ equals zero except for  selected local regions in the space-time-parameter domain.

It should be noticed that equations \eqref{eq:XTDdef} and \eqref{eq:XTDdefcontinu} only define the format of the XTD, in which all the terms are a priori unknown. They should be obtained from a training stage, which is described in the following.

\subsection{Offline training}
The XTD model can be trained either from data or directly from the physics-based governing equation of a given problem. The general training equation reads then
\begin{equation}
\displaystyle
\label{eq:XTDtraining}
\mathcal{L}(\ud^{\text{XTD}}(\X,t,{\mu_1},\dots,\mu_n))=0
\end{equation}
If the data of $\ud$ is provided, the above equation becomes
\begin{equation}
\displaystyle
\label{eq:XTDdata}
\int_{\Omega}(\delta\ud^{\text{XTD}})^T(\ud-\ud^{\text{XTD}})\ d\X\ dt\ d\mu_1\cdots d\mu_{n}=0
\end{equation}
where $\Omega$ stands for the general space-time-parameter domain.

Otherwise, we can consider $\mathcal{L}$ as a general differential operator, which represents the governing equation of a given problem. 

The solution strategy for solving equation \eqref{eq:XTDdata} is presented in \ref{XTDdata}. It confirms the performance of XTD for dealing with highly nonlinear (or singular) problems. In real applications, a large variation of parameters can involve non-smooth transitions of solutions, thus the XTD is a good choice for solving large space-time-parametric problems. 
We will focus on using the governing equation to train the XTD model in this work. The detailed algorithm will be explained in section \ref{section:XTDformulation}.

We remark that the training should cover a very large space-time-parameter domain so that it can be done once-for-all in the offline stage. The trained XTD solutions will be stored as a database that can be used repetitively in the online phase without expensive computations. 

\subsection{Online prediction}
Once the training stage is finished, we obtain an XTD model $\ud^{\text{XTD}}(\X,t,{\mu_1},\dots,\mu_n)$ containing the general parametric solutions for the underlying problem. Hence, for any unseen parameter value, we can use the XTD model to predict the solution, like neural networks. This online prediction is very fast, and is as cheap as a simple interpolation. This can enable many applications in the online stage, such as real-time system monitoring and control, uncertainty quantification of parameters, inverse identification and design, which would require prohibitive computational costs with conventional approaches.

As an example, we can formulate a general design problem as follows
\begin{equation}
\displaystyle
\label{eq:design}
\begin{cases}
\begin{aligned}
&\text{Find}\ \mu_1,\dots,\mu_n\\
&\text{Minimize}\ J(\ud^{\text{XTD}}(\X,t,{\mu_1},\dots,\mu_n))\\
&\text{Subject to:}\\
&c_1(\X,t,{\mu_1},\dots,\mu_n)=0\\
&c_2(\X,t,{\mu_1},\dots,\mu_n)\leq0
\end{aligned}
\end{cases}
\end{equation}
where $J$ is the objective function that needs to be evaluated repetitively during the optimization, $c_1,c_2$ are general constraints. Since the XTD model is very cheap to evaluate once the training is finished, we can expect a very fast design in the online stage. This can lead to orders of magnitude speed-ups compared to a conventional design strategy.

In what follows, we will present the XTD method to elastic-plastic problems as applications.
\section{The XTD method for nonlinear parametric elastic-plastic problems}
\label{section:XTDformulation}
\subsection{Problem formulation}
Considering a physical spatial domain $\Omega_\X$ with a prescribed external force $\overline{\textbf{d}}$ and displacement $\overline{\ud}$ respectively on the boundaries $\partial_{\textbf{F}}{\Omega_\X}$ and $\partial_{\textbf{u}}{\Omega_\X}$, and a design parameter vector $\boldsymbol{\mu}\in\Omega_{\boldsymbol{\mu}}$ that is related to materials, loading, or geometry. The quasi-static mechanical analysis consists in seeking the admissible stress, displacement, and plastic strain fields $\sig$, $\ud$ and $\epsp$ satisfying the following balance equations and boundary conditions
\begin{equation}
\displaystyle
\label{eq:MSF}
\begin{cases}
\begin{aligned}
&\nabla\cdot\sig (\X,t;\boldsymbol{\mu})=0\quad \text{in}\ \Omega_\X\\
&\sig(\X,t;\boldsymbol{\mu})\cdot\textbf{n}(\X,t;\boldsymbol{\mu})=\overline{\textbf{d}}(\X,t;\boldsymbol{\mu})\quad   \text{on}\ \partial_{\textbf{F}}{\Omega_\X} \\
&\ud(\X,t;\boldsymbol{\mu})=\overline{\ud}(\X,t;\boldsymbol{\mu})\quad \text{on}\ \partial_{\textbf{u}}{\Omega_\X}
\end{aligned}
\end{cases}
\end{equation}
where {the body forces are neglected}. Under the infinitesimal strain assumption, the stress field $\sig$ reads
\begin{equation}
\displaystyle
\label{eq:stress}
\sig =\textbf{D}:\epse =\textbf{D}:(\eps -\epsp)
\end{equation}
where the total strain $ \eps=\nabla_s \ud(\X,t)$,  $\epse$ and $\epsp$ denote respectively the elastic and the plastic strains,  and $\textbf{D}$ denotes the fourth order elasticity tensor.  The complementary plastic behavior law with an isotropic hardening R reads
\begin{equation}
\displaystyle
\label{eq:MBL}
\begin{cases}
\begin{aligned}
&f(\sig ,\text{p})=\parallel {{\sig }_{d}}\parallel -{\sigma }_{y}-\text{R}(\text{p})\\
&{\dot{\eps }_p}=\mathcal{H}(f)\frac{<{{\sig }_{d}}:{{{\dot{\sig }}}_{d}}>}{\text{R}g}\frac{{{\sig}_{d}}}{\parallel {{\sig }_{d}}\parallel } 
\end{aligned}
\end{cases}
\end{equation}
where ${\sigma }_{y}$ denotes the initial yield stress, $\mathcal{H}(\bullet)$ the Heaviside function, ${\sig }_{d}$  the deviatoric stress tensor, $<A>$ the positive part of A, and p  the equivalent plastic strains and $g=\frac{d\text{R}}{d\text{p}}$.

The weak form of \eqref{eq:MSF} can be written as
\begin{equation}
\displaystyle
\label{eq:weakform}
-\int_{\Omega_\X}{\sig(\ud) :\eps ({{\ud}^{*}})d\X}+\int_{\partial_{\textbf{F}} {{\Omega_\X}}}{\bar{\textbf{d}}\cdot{{\ud}^{*}}ds=0}
\end{equation}
with ${\ud}^{*}$ denoting the test function.

The problem can be discretized by the finite element method and solved by a nonlinear solver (e.g. Newton-Raphson). This consists in solving, at each time step $k$ and for a given parameter vector $\boldsymbol{\mu}$, the following discretized system
\begin{equation}
\displaystyle
\label{eq:FEdiscret}
\begin{aligned}
 \textbf{R}(t_k;\boldsymbol{\mu})={\Fm}_{\text{ext}}(t_k;\boldsymbol{\mu})-{\Fm}_{\text{int}}(t_k;\boldsymbol{\mu})=0
\end{aligned}
\end{equation}
where 
\begin{equation}
\displaystyle
\label{eq:FEforce}
\begin{aligned}
&{\Fm}_{\text{ext}}(t_k;\boldsymbol{\mu}) =\int_{\partial_{\textbf{F}} {{\Omega_\X}}}{\textbf{N}^{T}\bar{\textbf{d}}}\ ds\\
&{\Fm}_{\text{int}}(t_k;\boldsymbol{\mu}) =\int_{\Omega_\X}{\textbf{B}^{T}\sig\ d\X}\\
\end{aligned}
\end{equation}
with the usual shape function and its gradient denoted by $\textbf{N}$ and $\textbf{B}$ respectively. ${\Fm}_{\text{ext}}$ and ${\Fm}_{\text{int}}$ stand for the external and internal nodal forces. Considering the equation \eqref{eq:stress}, the discretized form \eqref{eq:FEdiscret} becomes
\begin{equation}
\displaystyle
\label{eq:FEdiscret2}
\begin{aligned}
 \textbf{R}(t_k;\boldsymbol{\mu})={\Fm}_{\text{ext}}(t_k;\boldsymbol{\mu})-\textbf{K}(\boldsymbol{\mu})\textbf{U}(t_k;\boldsymbol{\mu})+{\Fm}_{\text{pl}}(t_k;\boldsymbol{\mu})=0
\end{aligned}
\end{equation}
where 
\begin{equation}
\displaystyle
\label{eq:Fpl}
\begin{aligned}
&\textbf{K} =\int_{\Omega_\X}{\textbf{B}^{T}\textbf{D}\textbf{B}\ d\X}\\
&{\Fm}_{\text{pl}} =\int_{\Omega_\X}{\textbf{B}^{T}\textbf{D}:\epsp\ d\X}\\
\end{aligned}
\end{equation}
with $\textbf{U}$ as the usual nodal displacement vector. As the solution is computed for a specific parameter $\boldsymbol{\mu}$ in this case, Eq. \eqref{eq:FEdiscret2} has to be solved many times for different $\boldsymbol{\mu}$ to obtain the parametric solutions of the problem. This is the conventional FEA-based approach for computing parametric solutions.
%\subsection{Challenges with PGD based parametric model reduction}
%The usual way to solve the parametric system \eqref{eq:FEdiscret} can be seen as a point-wise evaluation in the parameter space $\Omega_{\boldsymbol{\mu}}$, which may result in a prohibitive computational cost if an exhaustive exploration of the space is required.
%\subsection{General parametric problems for nonlinear solids}
%In a more general sense, the materials can involve other nonlinear effects, like the viscosity. Denoting by $\eps_{nl}$ the part of strain induced by nonlinear effects, the stress can be then written as
%\begin{equation}
%\displaystyle
%\label{eq:stressgeneral}
%\sig =\textbf{D}:\epse =\textbf{D}:(\eps -\eps_{nl})
%\end{equation}
%and the complementary constitutive law can be written implicitly as
%\begin{equation}
%\displaystyle
%\dot{\eps}_{nl} = f_{ev} (\sig,\eps_{nl}) 
%\end{equation}

\subsection{{Extended} tensor decomposition for elastic-plastic problems}

Different from FEA, the XTD method can compute all the parametric solutions simultaneously. For solving the aforementioned parametric space-time problem, 
we chose the increment of the displacement field as the separated variable, which results in
\begin{equation}
\displaystyle
\label{eq:XTDdu}
\Delta\ud(\X,{\mu_1},\dots,\mu_n)= \sum_{m=1}^{M}\boldsymbol{a}^{(m)}(\X)g_1^{(m)}(\mu_1)\cdots g_n^{(m)}(\mu_n) + \sum_{k=1}^{K}\Tilde{\boldsymbol{u}}^{(k)}(\X,{\mu_1},\dots,\mu_n)
\end{equation}
where the $\Delta\ud$ stands for the displacement increment at each time step. The $\mu_i|_{i=1,\dots,n}$ are the components of the design parameter vector $\boldsymbol{\mu}$, and therefore the underlying variable parameters that relate to materials, loading, or geometry. Each parameter $\mu_i$ belongs to a  given range $\Omega_{\mu_i}$, and $\Omega_{\boldsymbol{\mu}}=\Omega_{\mu_1}\times \Omega_{\mu_2}\times\cdots\times\Omega_{\mu_n}$. With this definition, the global computational domain of XTD becomes a space-parameter domain, i.e., $\Omega=\Omega_\X\times\Omega_{\boldsymbol{\mu}}$.

In general, we can also choose the full displacement $\ud$ for the decomposition. The choice in our work seems appropriate for having a good reducibility of the problem and fits well the actual FE solution scheme for plasticity. In addition, since the sparsity of $\Tilde{\boldsymbol{u}}^{(k)}$ is essential to ensure the efficiency of the XTD method, we assume that $\Tilde{\boldsymbol{u}}^{(k)}=0$ holds true for the most part of the space-parameter domain $\Omega=\Omega_\X\times\Omega_{\boldsymbol{\mu}}$ for a displacement increment. We expect that the non-zero terms of $\Tilde{\boldsymbol{u}}^{(k)}$ are related to the contribution of plastic strains. Hence, the plastic region (at least the increment of the plastic region) can be relatively small compared to the global space-parameter domain $\Omega=\Omega_\X\times\Omega_{\boldsymbol{\mu}}$. 

The overall idea for finding the XTD model solution is to  compute alternatively the separated terms and the non-separated terms in equation \eqref{eq:XTDdu}. In particular, when computing the separated terms, the problem is considered as a linear problem with plasticity-induced internal forces as external loads. Computing the non-separated terms requires solving the nonlinear plastic problem but only for a small portion of the global domain. The overall solution procedure follows an incremental strategy, which starts from $M=1$ and $K=1$. For a general purpose, let us assume that the $M-1$ separated modes and $K-1$ non-separated extended modes have been computed, the (to be computed) unknown modes  are  denoted by $\boldsymbol{a},g_1,\cdots, g_n$ and $\Tilde{\boldsymbol{u}}$. Then
\begin{equation}
\displaystyle
\Delta\ud(\X,{\mu_1},\dots,\mu_n)= \sum_{m=1}^{M-1}\boldsymbol{a}^{(m)}g_1^{(m)}\cdots g_n^{(m)} + \boldsymbol{a}g_1\cdots g_n + \sum_{k=1}^{K-1}\Tilde{\boldsymbol{u}}^{(k)}+\Tilde{\boldsymbol{u}}
\end{equation}
The {weak} form \eqref{eq:weakform} of the problem for training the XTD model reads then
\begin{equation}
\displaystyle
-\int_{\Omega}{\sig(\ud_0+\Delta\ud) :\eps ({\Delta{\ud}^{*}})\ d\X d\mu_1\cdots d\mu_n}+\int_{\partial_{\textbf{F}} {{\Omega}}}{\bar{\textbf{d}}\cdot{\Delta{\ud}^{*}}\ dsd\mu_1\cdots d\mu_n=0}
\end{equation}
where $\partial_{\textbf{F}} {{\Omega}}=\partial_{\textbf{F}} {{\Omega_\X}}\times\Omega_{\boldsymbol{\mu}}$, $\Delta\ud$ is the displacement increment at the current time step $t_k$, $\ud_0$ is the displacement at the previous step $t_{k-1}$. The test function
$\Delta\ud^*=\boldsymbol{a}^*g_1\cdots g_n+\boldsymbol{a}g_1^*\cdots g_n+\cdots+\boldsymbol{a}g_1\cdots g_n^*+\Tilde{\boldsymbol{u}}^*$.

\subsubsection{Formulation to compute the separated modes \label{formulationsepmode}}
The unknown modes are computed alternatively by considering the others are fixed. For example, we compute $\boldsymbol{a}$ by considering that $g_1,\cdots, g_n$ and  $\Tilde{\boldsymbol{u}}$ are given. Then with the FE discretization, we have
\begin{equation}
\displaystyle
\begin{aligned}
&\int_{\Omega}{g_n^{T}\cdots g_1^{T}\boldsymbol{\mathcal{A}}^{*T}\sum_{m=1}^{M-1}\textbf{B}^T \textbf{D}\textbf{B}\boldsymbol{\mathcal{A}}^{(m)}g_1^{(m)}\cdots g_n^{(m)}\ d\X d\mu_1\cdots d\mu_n}\\
&+\int_{\Omega}{g_n^{T}\cdots g_1^{T}\boldsymbol{\mathcal{A}}^{*T}\textbf{B}^T\textbf{D}\textbf{B}\boldsymbol{\mathcal{A}}g_1\cdots g_n\ d\X d\mu_1\cdots d\mu_n}\\
&+\int_{\Omega}{g_n^{T}\cdots g_1^{T}\boldsymbol{\mathcal{A}}^{*T}\sum_{k=1}^{K-1}\textbf{B}^T \textbf{D}\textbf{B}\Tilde{\boldsymbol{\mathcal{U}}}^{(k)}\ d\X d\mu_1\cdots d\mu_n}\\
&-\int_{\Omega_{\boldsymbol{\mu}}}{g_n^{T}\cdots g_1^{T}\boldsymbol{\mathcal{A}}^{*T}\Delta\textbf{F}_{\text{pl}}\ d\mu_1\cdots d\mu_n}\\
&=\int_{\partial_{\textbf{F}} {{\Omega}}}{g_n^{T}\cdots g_1^{T}\boldsymbol{\mathcal{A}}^{*T}\textbf{N}^T\Delta\bar{\textbf{d}}\ ds d\mu_1\cdots d\mu_n}
\end{aligned}
\end{equation}
where $\boldsymbol{\mathcal{A}}^{(m)}, \boldsymbol{\mathcal{A}}$ are the discrete representation of  $\boldsymbol{a}^{(m)},\boldsymbol{a}$ as introduced in \eqref{eq:XTDdef}. The detailed derivation can be found in \ref{XTDmodel}. Here, we assume that an estimate of $\Delta\textbf{F}_{\text{pl}}$, {with $\Delta\textbf{F}_{\text{pl}}=\int_{\Omega_\X}{\textbf{B}^{T}\textbf{D}:\Delta\epsp\ d\X}$,} is also given. Therefore, the only unknown in the above equation is $\boldsymbol{\mathcal{A}}$.

For computing $g_1$, we have 
\begin{equation}
\displaystyle
\begin{aligned}
&\int_{\Omega}{g_n^{T}\cdots g_1^{*T}\boldsymbol{\mathcal{A}}^{T}\sum_{m=1}^{M-1}\textbf{B}^T \textbf{D}\textbf{B}\boldsymbol{\mathcal{A}}^{(m)}g_1^{(m)}\cdots g_n^{(m)}\ d\X d\mu_1\cdots d\mu_n}\\
&+\int_{\Omega}{g_n^{T}\cdots g_1^{*T}\boldsymbol{\mathcal{A}}^{T}\textbf{B}^T\textbf{D}\textbf{B}\boldsymbol{\mathcal{A}}g_1\cdots g_n\ d\X d\mu_1\cdots d\mu_n}\\
&+\int_{\Omega}{g_n^{T}\cdots g_1^{*T}\boldsymbol{\mathcal{A}}^{T}\sum_{k=1}^{K-1}\textbf{B}^T \textbf{D}\textbf{B}\Tilde{\boldsymbol{\mathcal{U}}}^{(k)}\ d\X d\mu_1\cdots d\mu_n}\\
&-\int_{\Omega_{\boldsymbol{\mu}}}{g_n^{T}\cdots g_1^{*T}\boldsymbol{\mathcal{A}}^{T}\Delta\textbf{F}_{\text{pl}}\ d\mu_1\cdots d\mu_n}\\
&=\int_{\partial_{\textbf{F}} {{\Omega}}}{g_n^{T}\cdots g_1^{*T}\boldsymbol{\mathcal{A}}^{T}\textbf{N}^T\Delta\bar{\textbf{d}}\ ds d\mu_1\cdots d\mu_n}
\end{aligned}
\end{equation}
with that the $g_1$ is the only unknown variable.

Analogically, for computing $g_n$, we have 
\begin{equation}
\displaystyle
\begin{aligned}
&\int_{\Omega}{g_n^{*T}\cdots g_1^{T}\boldsymbol{\mathcal{A}}^{T}\sum_{m=1}^{M-1}\textbf{B}^T \textbf{D}\textbf{B}\boldsymbol{\mathcal{A}}^{(m)}g_1^{(m)}\cdots g_n^{(m)}\ d\X d\mu_1\cdots d\mu_n}\\
&+\int_{\Omega}{g_n^{*T}\cdots g_1^{T}\boldsymbol{\mathcal{A}}^{T}\textbf{B}^T\textbf{D}\textbf{B}\boldsymbol{\mathcal{A}}g_1\cdots g_n\ d\X d\mu_1\cdots d\mu_n}\\
&+\int_{\Omega}{g_n^{*T}\cdots g_1^{T}\boldsymbol{\mathcal{A}}^{T}\sum_{k=1}^{K-1}\textbf{B}^T \textbf{D}\textbf{B}\Tilde{\boldsymbol{\mathcal{U}}}^{(k)}\ d\X d\mu_1\cdots d\mu_n}\\
&-\int_{\Omega_{\boldsymbol{\mu}}}{g_n^{*T}\cdots g_1^{T}\boldsymbol{\mathcal{A}}^{T}\Delta\textbf{F}_{\text{pl}}\ d\mu_1\cdots d\mu_n}\\
&=\int_{\partial_{\textbf{F}} {{\Omega}}}{g_n^{*T}\cdots g_1^{T}\boldsymbol{\mathcal{A}}^{T}\textbf{N}^T\Delta\bar{\textbf{d}}\ ds d\mu_1\cdots d\mu_n}
\end{aligned}
\end{equation}
The above equations can be solved repetitively until $|\boldsymbol{a}g_1\cdots g_n|$ converges. {From our experience, the convergence is usually achieved within 4 or 5 iterations with a relative tolerance $10^{-4}$. The computational cost can be very cheap with simple matrix-vector multiplications.} The final converged solution is the solution for one mode. We can incrementally increase the number $M$ to enrich the displacement solution. The converged number can be the one for which the enrichment becomes very small.

\subsubsection{Formulation to compute the non-separated extended modes \label{formulationextmode}}
Once the separated modes are computed, we can enrich the extended modes. This step also consists in correcting the estimated plastic strain and internal forces with the updated displacement increment. 

Taking the current displacement as an initial guess  $\Delta\ud^{(0)}$, we can {determine the local plastic region by computing the increment of plastic strain due to the displacement increment or simply checking the yield criterion. This can be done by a point-wise like evaluation in the space-parameter domain in a parallel setting. The plastic region is then the local enrichment region in this case, which reads}
 \begin{equation}
    \displaystyle
        \Omega_\text{pl}=\{(\X,\boldsymbol{\mu})\in \Omega_\X\times\Omega_{\boldsymbol{\mu}}\ | \ \Delta\epsp(\X,\boldsymbol{\mu})\neq 0\}
    \end{equation}
Then, Newton's method can be applied to solve an updated $\Delta\ud$ satisfying the balance equation 
 \begin{equation}
    \displaystyle
    {{\Fm}_{\text{ext}}-{\Fm}_{\text{int}}=0 \quad \text{in}\quad \Omega_{\text{pl}} \quad \text{with } \quad \Delta\ud|_{\partial \Omega_{\text{pl}} }=\Delta\ud|_{\partial \Omega_{\text{el}} }}
 \end{equation}
 where $\Omega_{\text{el}}=(\Omega_\X\times\Omega_{\boldsymbol{\mu}})\backslash \Omega_{\text{pl}}$. {Note that the above equation can be solved in a standard way, like FEA. The internal variables, such as plastic strain $\epsp$, can be computed using a classic return-mapping algorithm. The history dependant behavior of materials is naturally taken into account.} The converged solution gives a new extended mode as
  \begin{equation}
    \displaystyle
     \Tilde{\boldsymbol{u}}=\Delta\ud-\Delta\ud^{(0)}
 \end{equation}
 with the updated $\epsp$, $\Delta\epsp$ and $\Delta\textbf{F}_{\text{pl}}$. They are sparse arrays as the computations are performed locally.  A notable advantage of this formulation is that it enables an  easily implemented parallelization and the possibility to incorporate with an existing FE code.

\subsubsection{Overall solution procedure}
The XTD model is trained by solving alternatively the separated and extended modes using the previous formulation. The overall solution procedure is summarized below.
\begin{enumerate}
    \item  For a given time step $t=t_k$, starting from $M=1, K=1$ and assuming
    \begin{equation}
    \displaystyle
        \Tilde{\boldsymbol{u}}=0, \quad \Delta\epsp=0, \quad \forall (\X,\boldsymbol{\mu})\in \Omega_\X\times\Omega_{\boldsymbol{\mu}}
    \end{equation}
    \item Increase $M$ and compute incrementally the separated functions until a predefined convergence as explained in \ref{formulationsepmode} 
    \begin{equation}
    \displaystyle
        {\boldsymbol{a}^{(m)}, g_1^{(m)},\cdots, g_n^{(m)} }
    \end{equation}
    and update the nodal displacement increment $\Delta\ud$.
    \item Update the local plastic region $\Omega_\text{pl}$ with $\Delta\ud$ 
    \item Update the extended mode and the plastic strain increment as explained in \ref{formulationextmode}
    \begin{equation}
    \displaystyle
        \Tilde{\boldsymbol{u}}, \quad \Delta\epsp, \quad \forall (\X,\boldsymbol{\mu})\in \Omega_{\text{pl}}
    \end{equation}
    \item Check the convergence on $\Tilde{\boldsymbol{u}}$
    \begin{equation}
\displaystyle
\frac{\|\Tilde{\boldsymbol{u}}\|_\infty}{\|\Delta\ud\|_\infty}\leq\epsilon_{\text{XTD}}
\end{equation}
    \item If the convergence is achieved, stop computations. Otherwise, $K=K+1$ and go to step 2.
\end{enumerate}

We remark that in the XTD model, the internal variables like $\epsp(\X,\boldsymbol{\mu})$, $\sig(\X,\boldsymbol{\mu})$ are not stored in the separated form. Hence, the storage memory for these variables can be high when dealing with high dimensional space-time-parametric problems. To reduce the storage memory, we can perform a post-treatment using HOPGD \cite{lu2018multi} or the data-driven XTD (see \ref{XTDdata}) for data compression.

\section{Numerical examples}
\subsection{Elastic-plastic test with a variable loading}
\begin{figure}[htbp]
\centering 
\def\svgscale{0.4}
{%% Creator: Inkscape inkscape 0.92.4, www.inkscape.org
%% Accompanies image file 'Ex1domLoad.pdf' (pdf, eps, ps)
\begingroup%
  \makeatletter%
  \providecommand\color[2][]{%
    \errmessage{(Inkscape) Color is used for the text in Inkscape, but the package 'color.sty' is not loaded}%
    \renewcommand\color[2][]{}%
  }%
  \providecommand\transparent[1]{%
    \errmessage{(Inkscape) Transparency is used (non-zero) for the text in Inkscape, but the package 'transparent.sty' is not loaded}%
    \renewcommand\transparent[1]{}%
  }%
  \providecommand\rotatebox[2]{#2}%
  \newcommand*\fsize{\dimexpr\f@size pt\relax}%
  \newcommand*\lineheight[1]{\fontsize{\fsize}{#1\fsize}\selectfont}%
  \ifx\svgwidth\undefined%
    \setlength{\unitlength}{751.63508858bp}%
    \ifx\svgscale\undefined%
      \relax%
    \else%
      \setlength{\unitlength}{\unitlength * \real{\svgscale}}%
    \fi%
  \else%
    \setlength{\unitlength}{\svgwidth}%
  \fi%
  \global\let\svgwidth\undefined%
  \global\let\svgscale\undefined%
  \makeatother%
  \begin{picture}(1,0.43372107)%
    \lineheight{1}%
    \setlength\tabcolsep{0pt}%
    \put(0,0){\includegraphics[width=\unitlength,page=1]{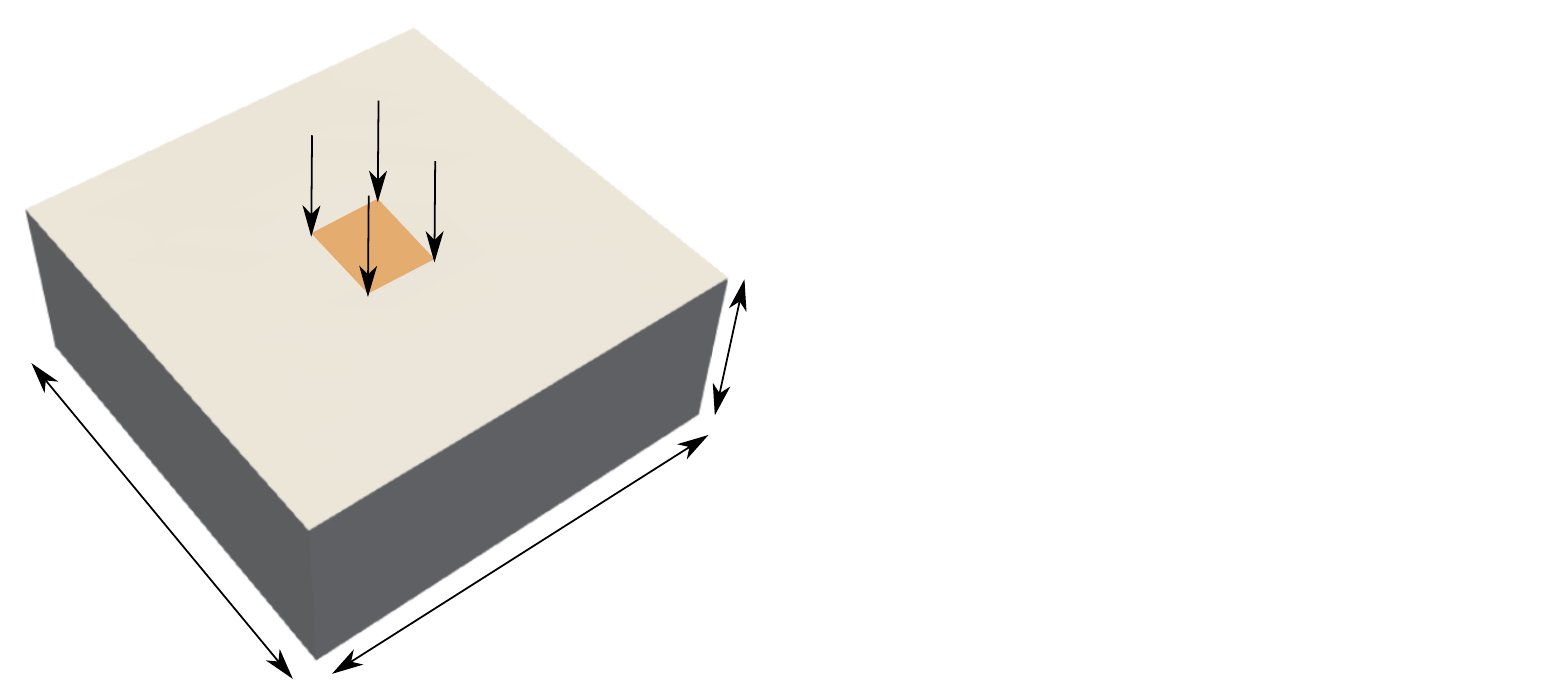}}%
    \put(0.04347986,0.07194613){\color[rgb]{0,0,0}\makebox(0,0)[lt]{\lineheight{1.25}\smash{\begin{tabular}[t]{l}$l_1$\end{tabular}}}}%
    \put(0.3467331,0.04477185){\color[rgb]{0,0,0}\makebox(0,0)[lt]{\lineheight{1.25}\smash{\begin{tabular}[t]{l}$l_2$\end{tabular}}}}%
    \put(0.48697658,0.19476561){\color[rgb]{0,0,0}\makebox(0,0)[lt]{\lineheight{1.25}\smash{\begin{tabular}[t]{l}$l_3$\end{tabular}}}}%
    \put(0.26582466,0.35793107){\color[rgb]{0,0,0}\makebox(0,0)[lt]{\lineheight{1.25}\smash{\begin{tabular}[t]{l}$\bar{\ud}$\end{tabular}}}}%
    \put(0.6852628,0.26880074){\color[rgb]{0,0,0}\makebox(0,0)[lt]{\lineheight{1.25}\smash{\begin{tabular}[t]{l}$\bar{\ud}$\end{tabular}}}}%
    \put(0,0){\includegraphics[width=\unitlength,page=2]{Ex1domLoad.pdf}}%
  \end{picture}%
\endgroup%
}
\caption{Geometry of the material domain and the loading history (with respect to time).}
\label{fig:Ex1doman}
\end{figure}

The first example consists in testing the proposed XTD method on elastic-plastic materials with cyclic loading. The geometry of the material domain and the loading history are depicted in \figurename~\ref{fig:Ex1doman}. A concentrated displacement load is applied on the top-surface of the domain with the fixed boundary on the bottom surface. The following loading increment is used as a cycle: $\Delta \bar{u}_z=0.03$ m for the first 2 time steps, $\Delta \bar{u}_z=-0.03$ m for the step 3 to step 5, $\Delta \bar{u}_z=0.03$ m for the step 6. The Young's modulus and Poisson's ratio are respectively 210 GPa and 0.3.  Two design parameters are considered: the initial yield stress $\sigma_y\in [\sigma_y^{\min},\ \sigma_y^{\max}]$ and the linear hardening coefficient $H\in [H^{\min},\ H^{\max}]$. Then the XTD model aims at computing the following parametric solutions for each time step

\begin{equation}
\displaystyle
\label{eq:XTDex1}
\Delta\ud(\X,\sigma_y,H)= \sum_{m=1}^{M}\boldsymbol{a}^{(m)}(\X)g_1^{(m)}(\sigma_y)g_2^{(m)}(H) + \sum_{k=1}^{K}\Tilde{\boldsymbol{u}}^{(k)}(\X,{\sigma_y},H)
\end{equation}

For evaluating the performance of the solution strategy, we used 7 test cases as shown in \tablename~\ref{table:ex1test}. It can be noticed that the parameter space is assumed to be extremely large, although it may not be realistic. Nevertheless, we consider that the offline computations should cover as large as possible the design parameter space so that no more online computations will be required. This is also consistent with the once-for-all concept in the offline-online strategy. From this consideration, we tested the XTD method with large parameter variations.

\begin{table}[htbp]
\caption{Test cases with different meshes and discretizations in the parameter space}
\centering
\begin{tabular}{|c|c|c|c|c|c|c|c|}
\hline
Case & $l_1\times l_2 \times l_3$ $(m^3)$ &  Mesh &$ [\sigma_y^{\min},\ \sigma_y^{\max}]$ & Grid $\sigma_y$ &$ [H^{\min},\ H^{\max}]$& Grid  $H$ \\ \hline
1 & 6 $\times$ 6 $\times$ 2.5  &  13 $\times$ 13 $\times$ 6 & $[300,\ 900]$ MPa &11  &  $[21,\ 63]$ GPa & 5  \\ \hline
2 & 6 $\times$ 6 $\times$ 2.5  &  13 $\times$ 13 $\times$ 6 & $[300,\ 1800]$ MPa &26  &  $[21,\ 105]$ GPa & 9  \\ \hline
3 & 6 $\times$ 6 $\times$ 2.5  &  13 $\times$ 13 $\times$ 6 & $[300,\ 1800]$ MPa &51  &  $[21,\ 105]$ GPa & 17  \\ \hline
4 & 12$\times$ 12 $\times$ 2.5  &  25 $\times$ 25 $\times$ 6 & $[300,\ 900]$ MPa &11  &  $[21,\ 63]$ GPa & 5  \\ \hline
5 & 12$\times$ 12 $\times$ 2.5  &  25 $\times$ 25 $\times$ 6 & $[300,\ 1800]$ MPa &26  &  $[21,\ 105]$ GPa & 9   \\ \hline
6 & 12$\times$ 12 $\times$ 2.5  &  25 $\times$ 25 $\times$ 6& $[300,\ 1800]$ MPa &51  &  $[21,\ 105]$ GPa & 9  \\ \hline
7 & 12$\times$ 12 $\times$ 10  &  25 $\times$ 25 $\times$ 21& $[300,\ 1800]$ MPa &51  &  $[21,\ 105]$ GPa & 9  \\ \hline
\end{tabular}
\label{table:ex1test}\\
 Note: each design parameter is discretized using a uniform grid
\end{table}

\begin{figure}[htbp]
\centering
{\includegraphics[scale=0.45]{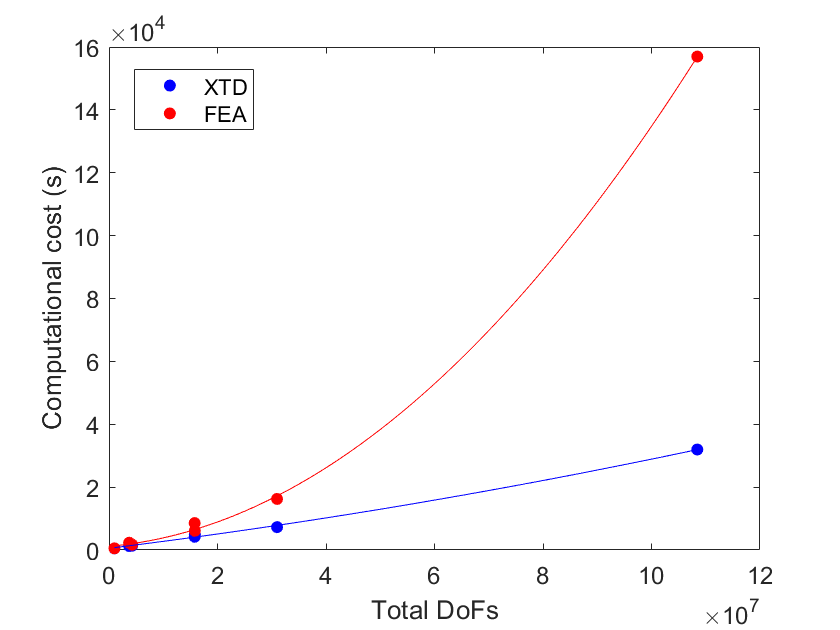}}
\caption{Computational cost against DoFs}
\label{fig:ex1cost}
\end{figure}

\tablename~\ref{table:ex1cost} reports the offline computational cost for different cases. For better illustrating the computational complexity, we compute the total degrees of freedom (DoFs) using the following equation
\begin{equation}
\displaystyle
\label{eq:dof}
\text{Total DoFs}=3\times\text{Mesh}\times \text{Time steps} \times \text{Grid }\ \sigma_y\times\text{Grid}\ H
\end{equation}
where 3 is the DoFs of the displacement per node, the number of time steps is 6. For comparison purposes, the finite element analysis (FEA) is run for the 7 test cases in a sequential setting. \figurename~\ref{fig:ex1cost} depicts the computational cost comparison between FEA and XTD in the 7 test cases. It is shown that the XTD has a significant speed-up when dealing with large DoFs. For small size problems (in terms of spatial mesh and discretization in parameter space), the computational cost is almost equivalent between the two methods. This is expected, as the XTD method is designed for computations in a large space-time-parameter domain. This cost analysis is only for the offline training. Once it is done, the online prediction is very fast. For the prediction on arbitrary parameter values, the online cost is less than 0.01 s. This kind of real time prediction is usually intractable with standard FEA.

\begin{table}[htbp]
\caption{Computational cost and prediction error for the test cases}
\centering
\begin{tabular}{|c|c|c|c|c|c|}
\hline
Case & Total DoFs &  Offline & Online& Err stress& Err plastic strain\\ \hline
1 & 1003860 & 502 s & $<0.01$ s & $<2\%$& $<2\%$\\ \hline
2 & 4270968 & 1355 s & $<0.01$ s & $1.8\%$& $1.5\%$\\ \hline
3 & 15824484 & 4965 s & $<0.01$ s & $1.8\%$& $1.8\%$\\ \hline
4 & 3712500 & 1264 s & $<0.01$ s & $<2\%$& $<2\%$\\ \hline
5 & 15795000 & 4207 s & $<0.01$ s & $1.9\%$& $1.4\%$\\ \hline
6 & 30982500 & 7276 s & $<0.01$ s & $1.9\%$& $1.5\%$\\ \hline
7 & 108438750 & 31959 s & $<0.01$ s & $<2\%$& $<2\%$\\ \hline
\end{tabular}
\label{table:ex1cost}\\
Note: Err stands for the mean error of the equivalent stress and plastic strain to finite element solutions on the test data set. 
\end{table}

\begin{figure}[htbp]
\centering
\subfigure[FEA, $t_3$]{\includegraphics[scale=0.2]{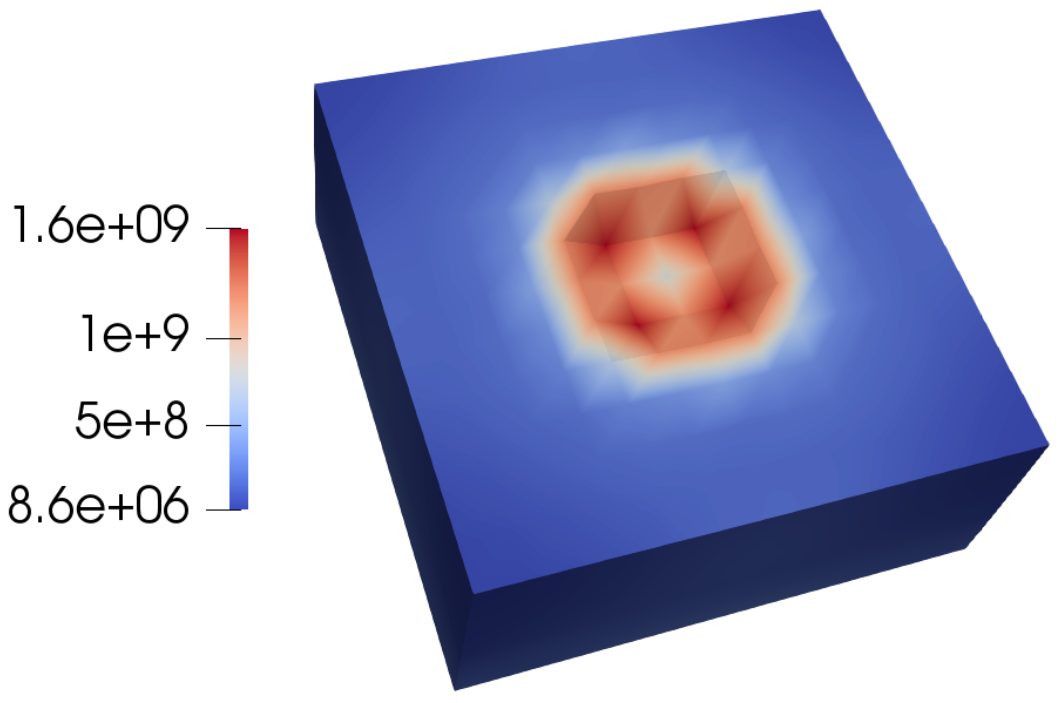}}
\subfigure[XTD, $t_3$]{\includegraphics[scale=0.2]{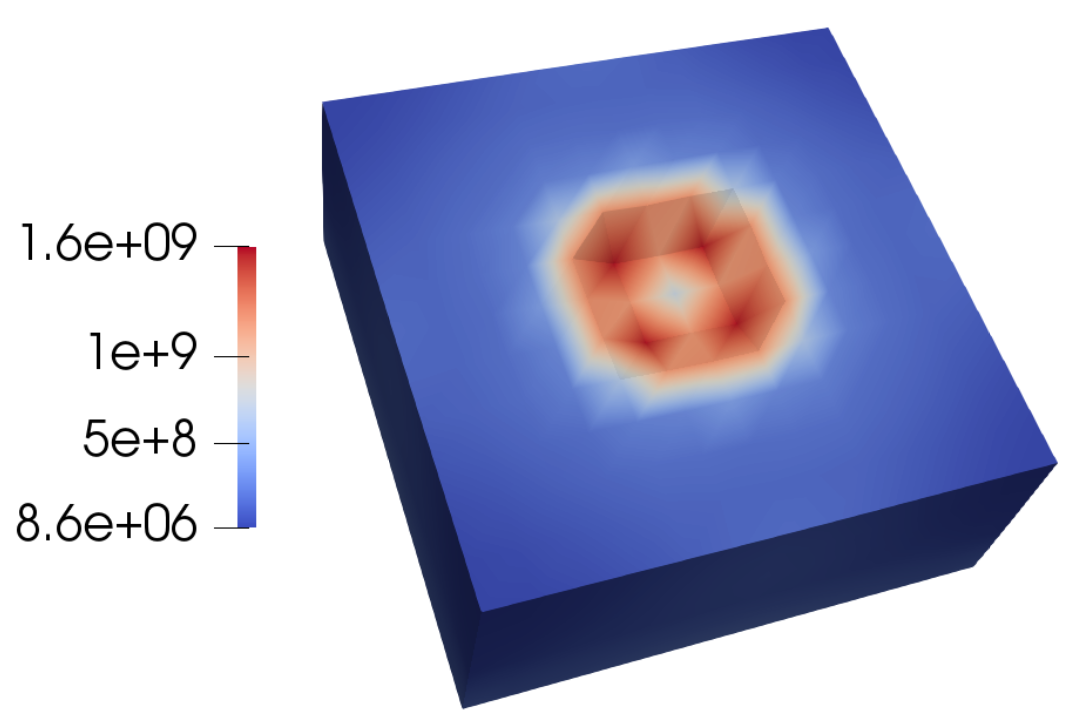}}
\subfigure[FEA,  $t_5$]{\includegraphics[scale=0.2]{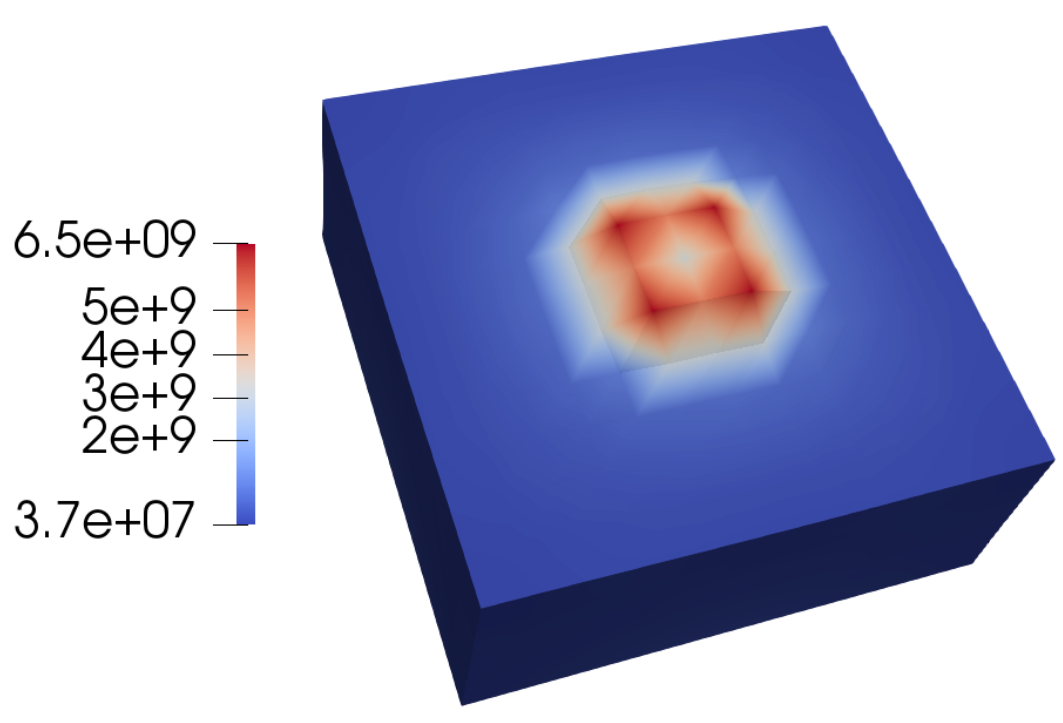}}
\subfigure[XTD, $t_5$ ]{\includegraphics[scale=0.2]{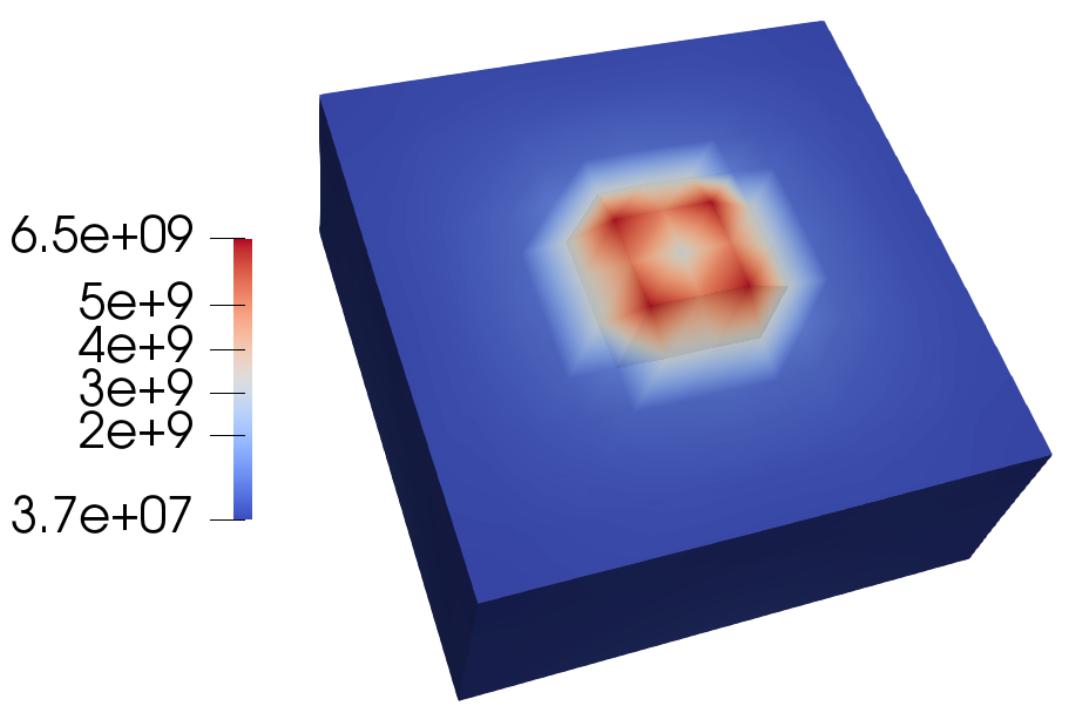}}\\
\subfigure[FEA, $t_3$]{\includegraphics[scale=0.2]{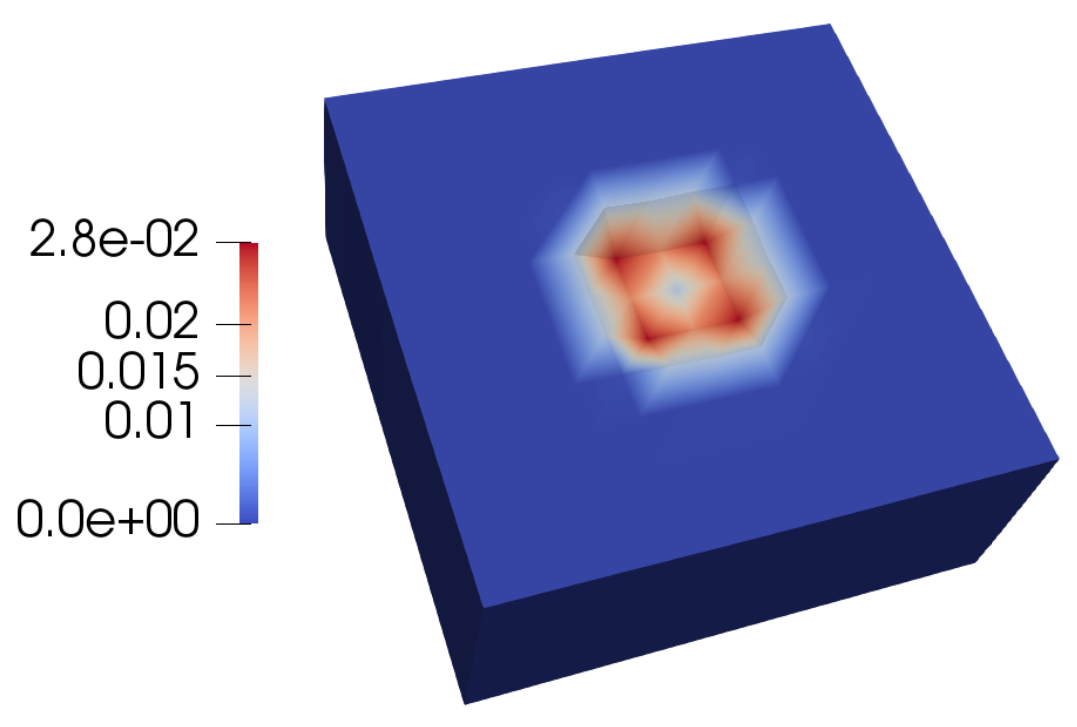}}
\subfigure[XTD, $t_3$]{\includegraphics[scale=0.2]{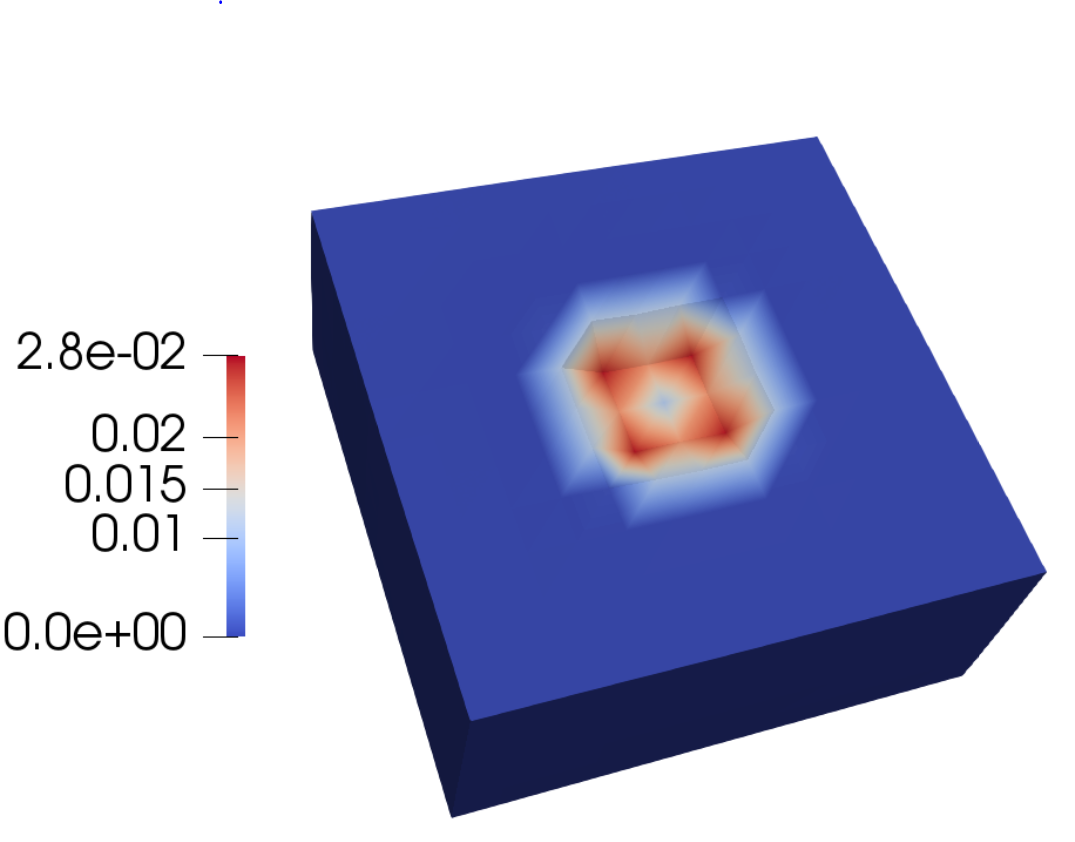}}
\subfigure[FEA,  $t_5$]{\includegraphics[scale=0.2]{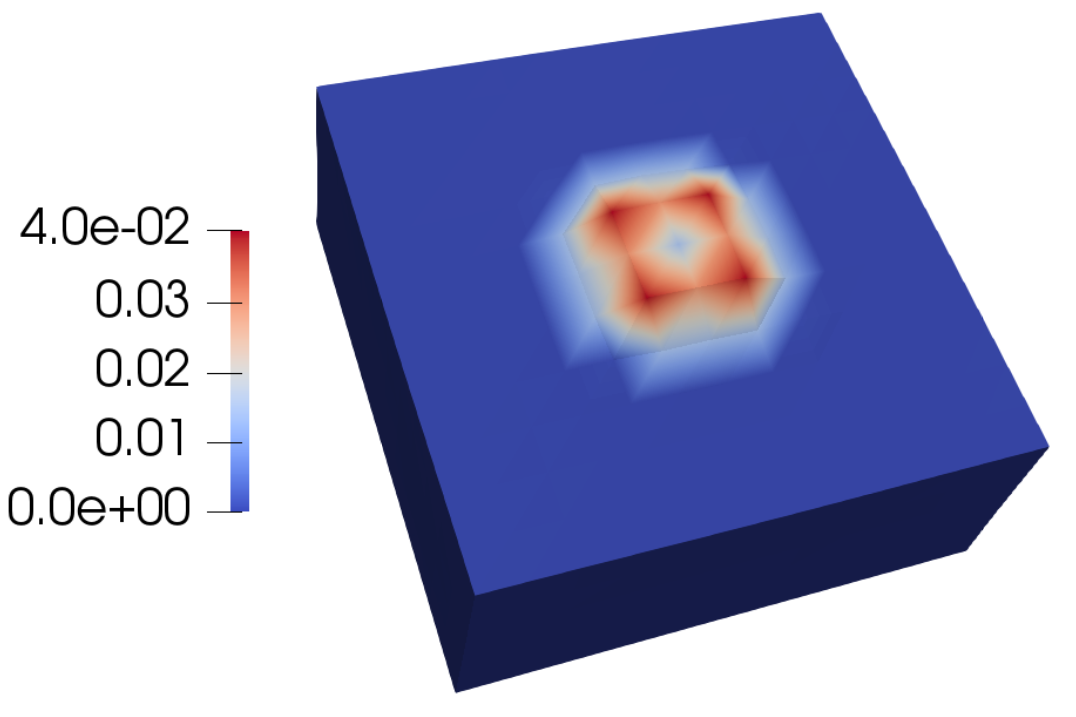}}
\subfigure[XTD, $t_5$ ]{\includegraphics[scale=0.2]{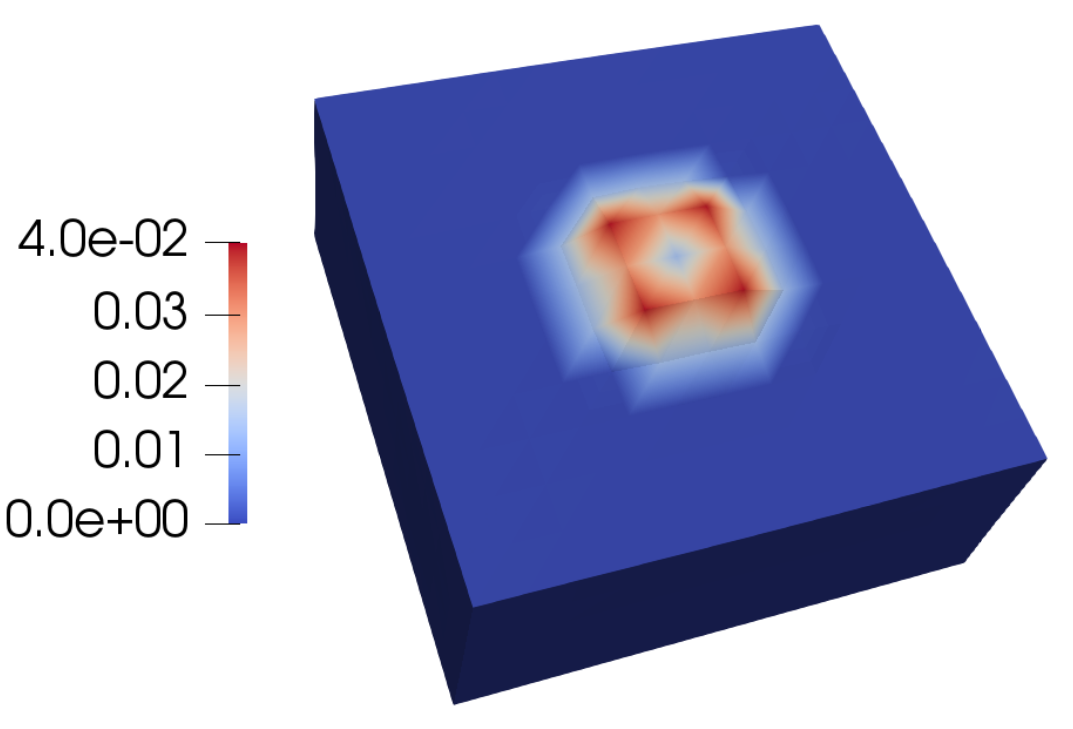}}
\caption{XTD solution versus reference FE solution for the parameter values $\sigma_y=1747$ MPa, $H=101.6$ GPa, test case 2. (a)-(d) von Mises stress, (e)-(h) equivalent plastic strain.}
\label{fig:Exmesh1compare}
\end{figure}

For further investigating the online prediction accuracy, a test database containing 10 randomly selected parameters has been constructed using FEA. As reported in \tablename~\ref{table:ex1cost}, 
the XTD predictions seem  at a well acceptable level. The solutions for the stress and plastic strain are well predicted within a mean error of $2\%$. \figurename~\ref{fig:Exmesh1compare} depicts an example of the online solution provided by XTD for the parameter values: $\sigma_y=1747$ MPa, $H=101.6$ GPa. The prediction agrees well with that obtained by FEA.

\begin{figure}[htbp]
\centering
\subfigure[Mode 1, $\boldsymbol{a}_3^{(1)}$]{\includegraphics[scale=0.2]{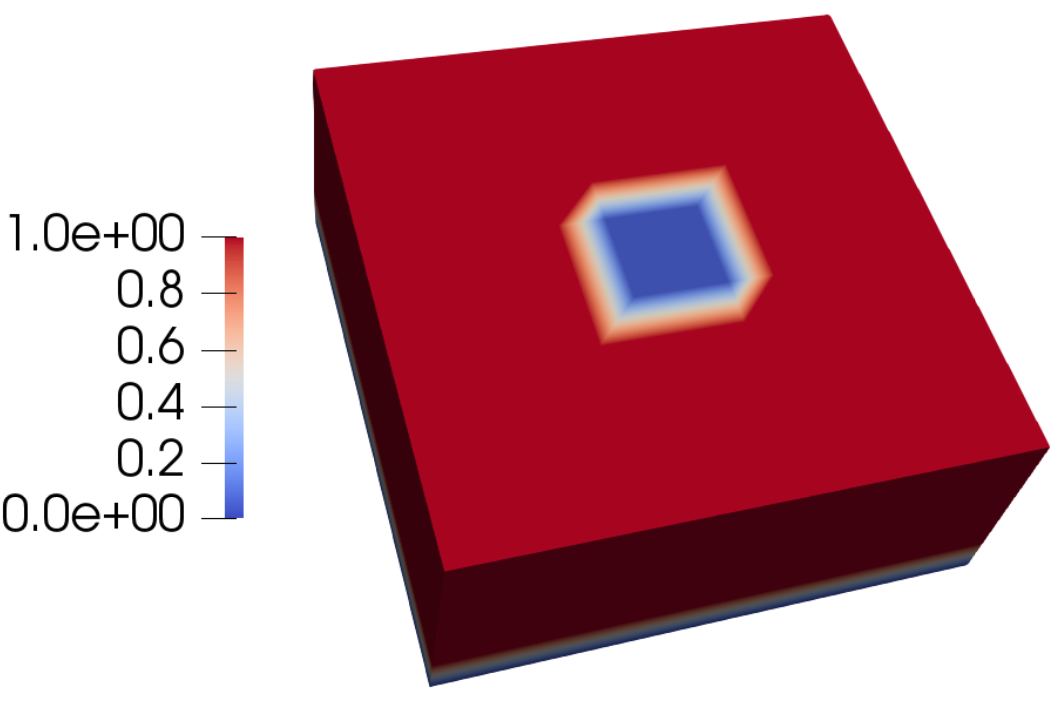}}
\subfigure[Mode 2, $\boldsymbol{a}_3^{(2)}$]{\includegraphics[scale=0.2]{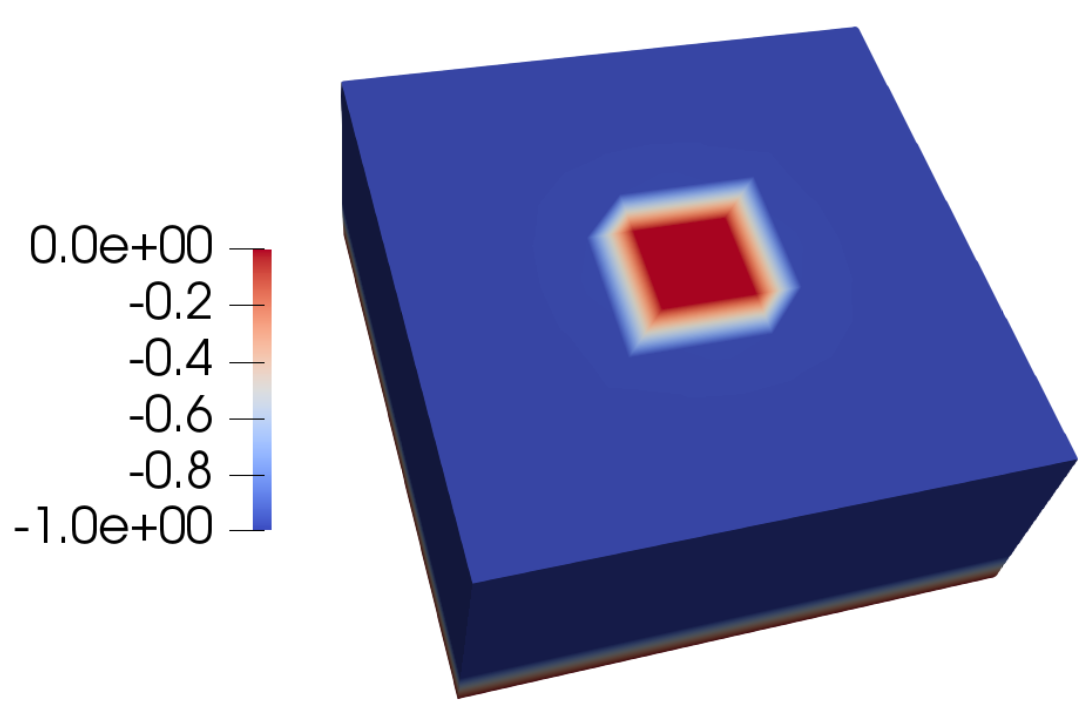}}
\subfigure[Mode 4,  $\boldsymbol{a}_3^{(4)}$]{\includegraphics[scale=0.2]{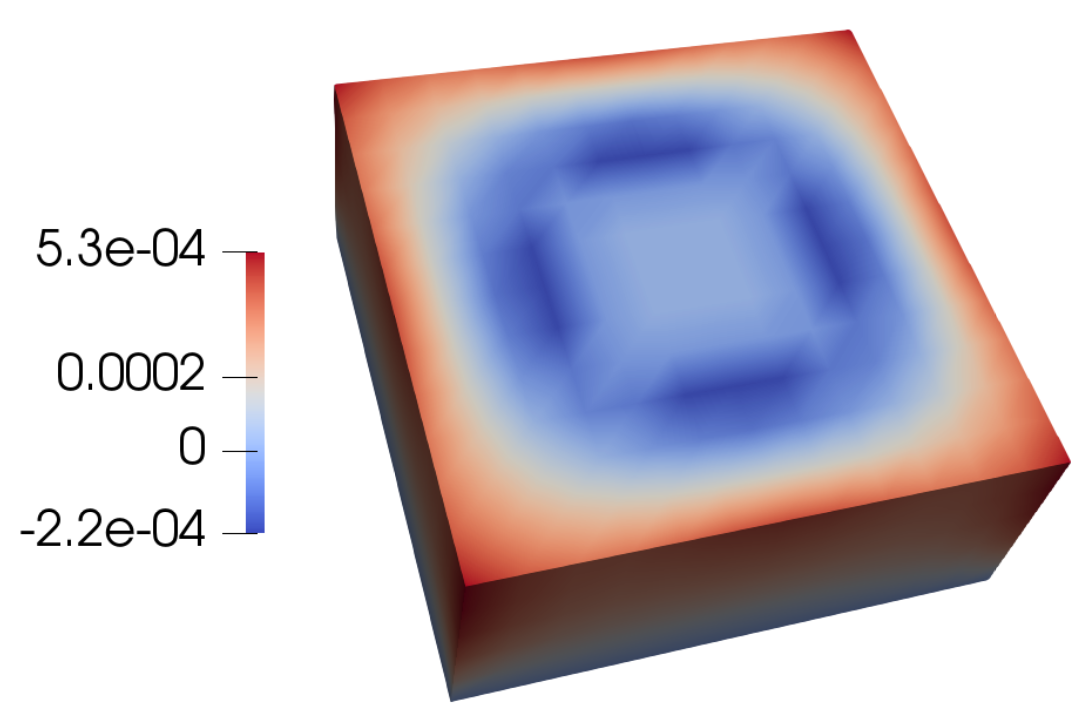}}
\subfigure[Mode 5,  $\boldsymbol{a}_3^{(5)}$ ]{\includegraphics[scale=0.2]{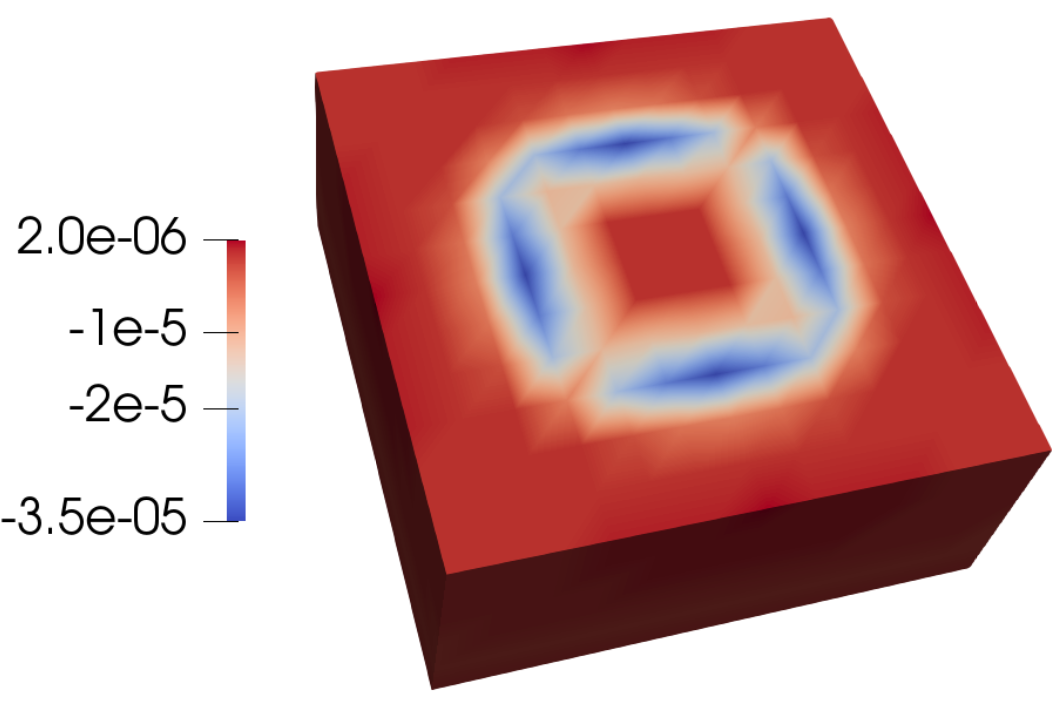}}\\
\subfigure[Mode 6,  $\boldsymbol{a}_3^{(6)}$]{\includegraphics[scale=0.2]{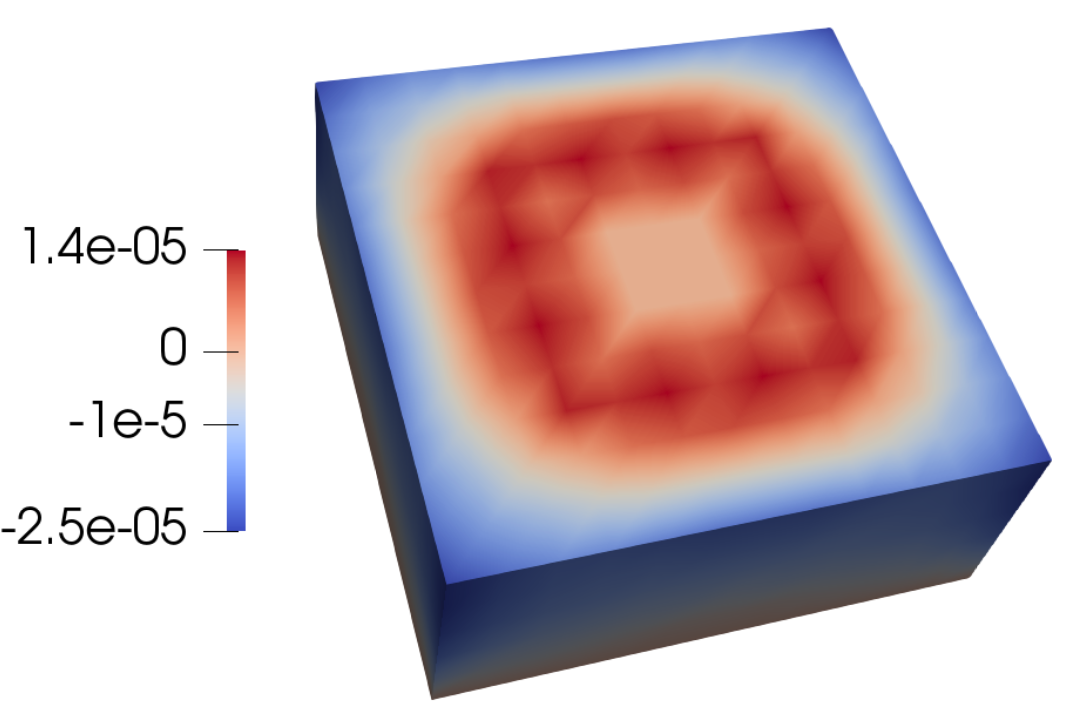}}
\subfigure[Mode 7,  $\boldsymbol{a}_3^{(7)}$]{\includegraphics[scale=0.2]{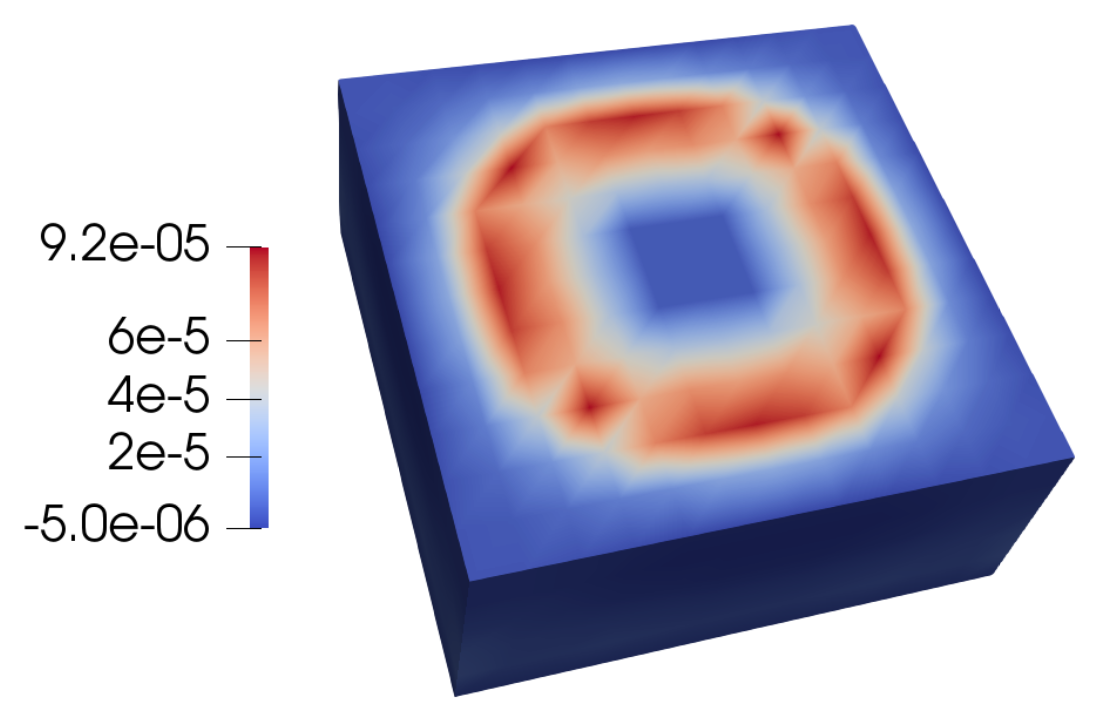}}
\subfigure[Mode 8,  $\boldsymbol{a}_3^{(8)}$]{\includegraphics[scale=0.2]{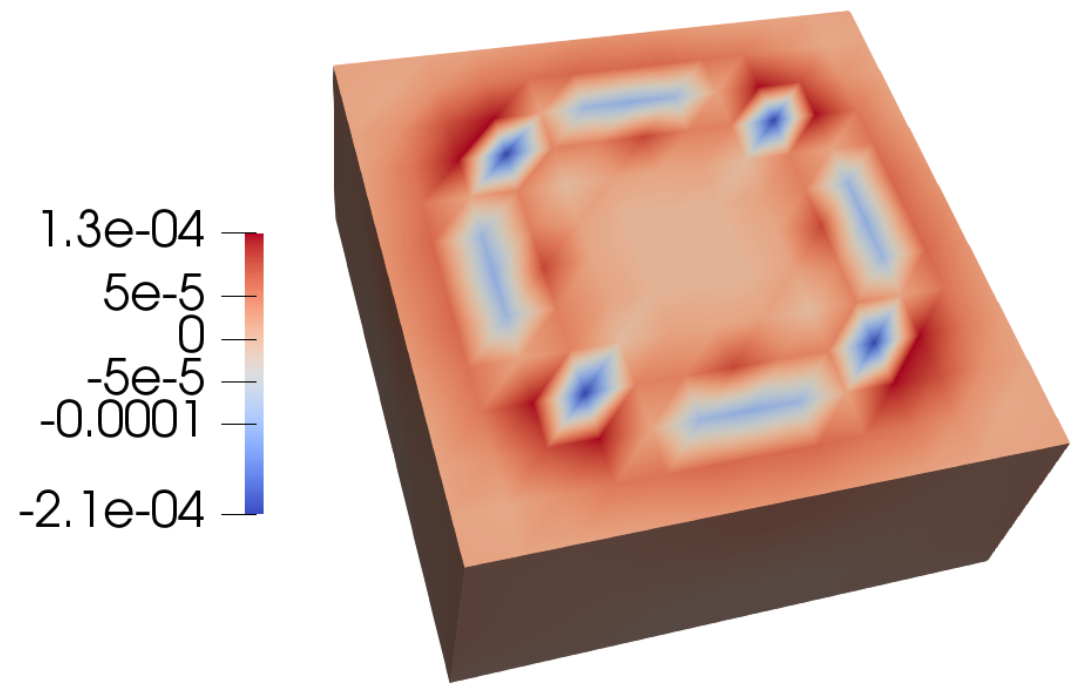}}
\subfigure[Mode 9,  $\boldsymbol{a}_3^{(9)}$]{\includegraphics[scale=0.2]{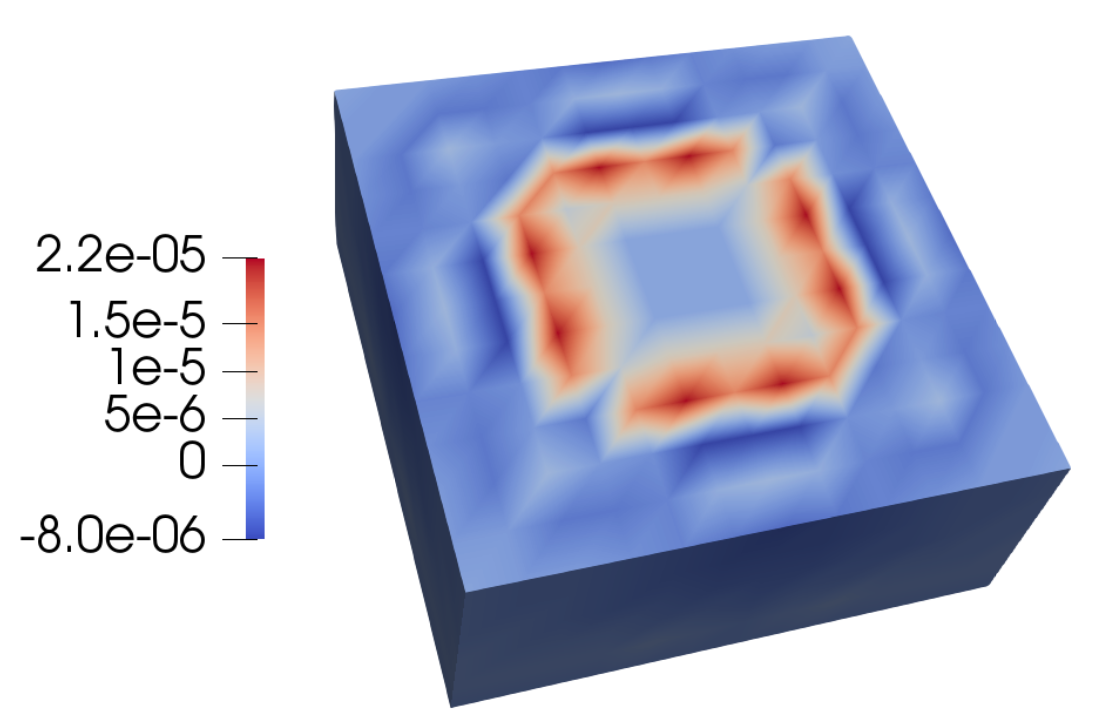}}
\caption{XTD separated modes $\boldsymbol{a}^{(m)}(\X)$ for test case 2, time step $t_5$, $z$-component $\boldsymbol{a}_3^{(m)}(\X)$.}
\label{fig:ExC2t5modesep}
\end{figure}

\begin{figure}[htbp]
\centering
\subfigure[$\Tilde{\boldsymbol{u}}_3^{(1)}(\X,360,21)$]{\includegraphics[scale=0.2]{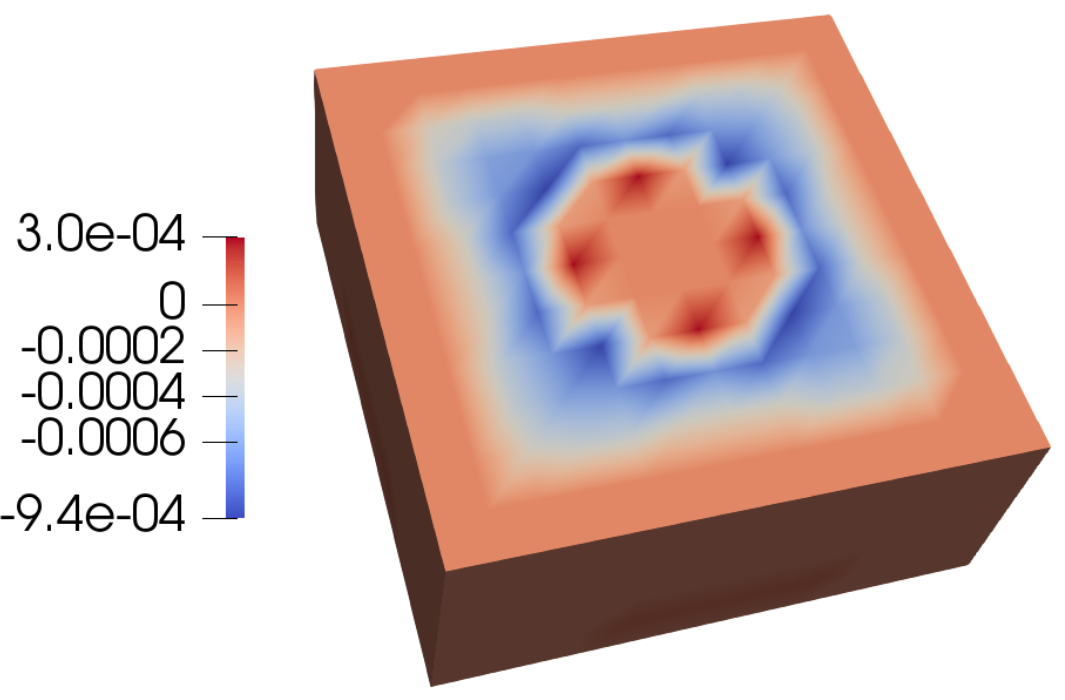}}
\subfigure[$\Tilde{\boldsymbol{u}}_3^{(1)}(\X,300,31.5)$]{\includegraphics[scale=0.2]{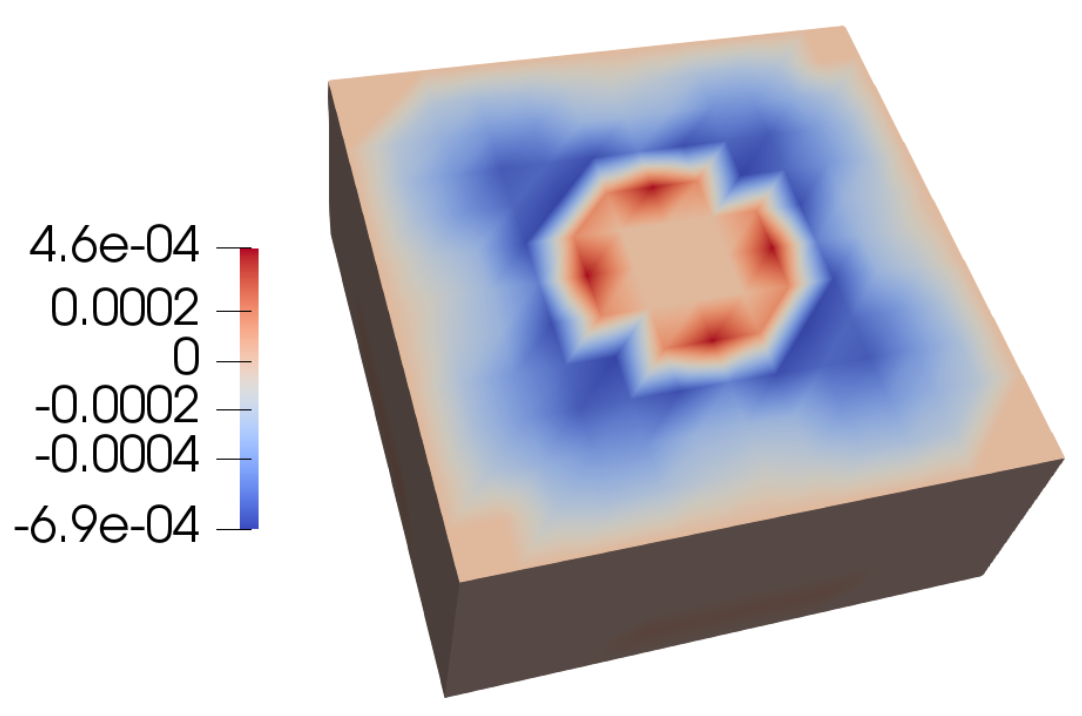}}
\subfigure[$\Tilde{\boldsymbol{u}}_3^{(1)}(\X,1380,52.5)$]{\includegraphics[scale=0.2]{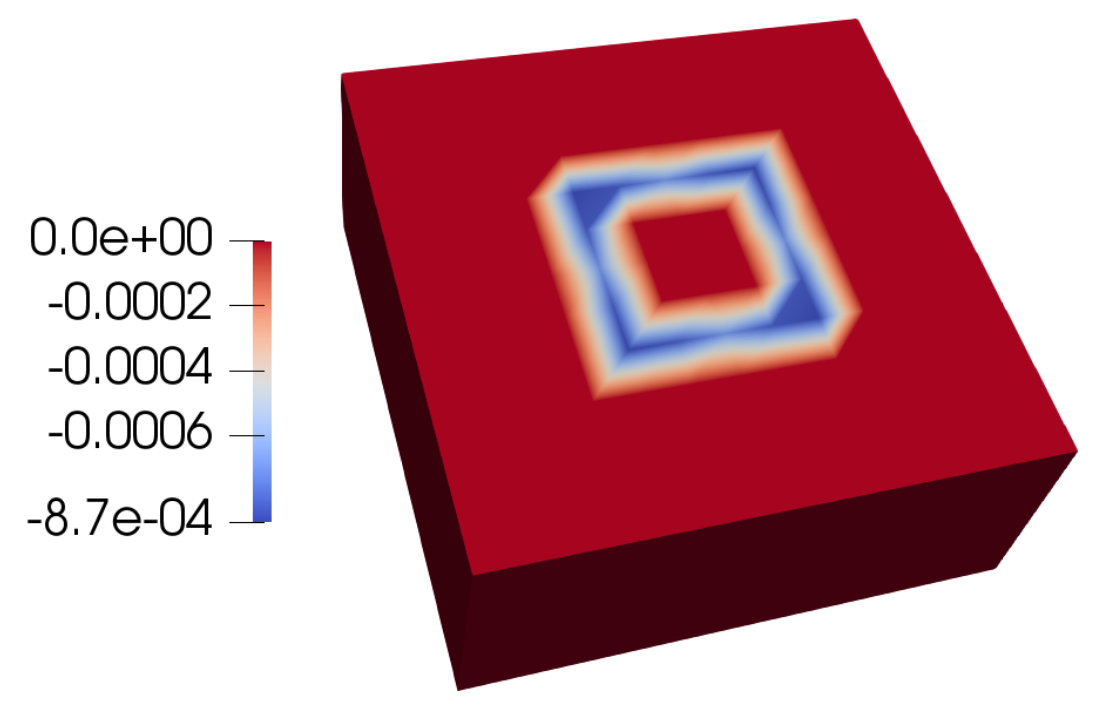}}
\subfigure[$\Tilde{\boldsymbol{u}}_3^{(1)}(\X,1140,94.5)$ ]{\includegraphics[scale=0.2]{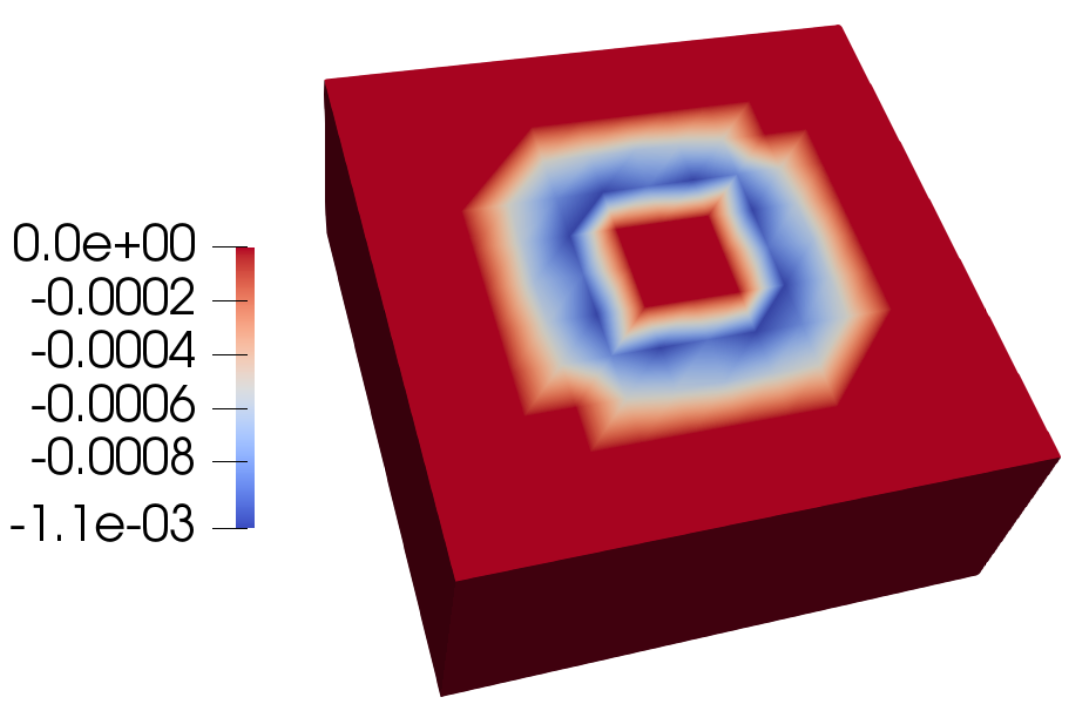}}\\
\subfigure[$\Tilde{\boldsymbol{u}}_3^{(2)}(\X,360,21)$]{\includegraphics[scale=0.2]{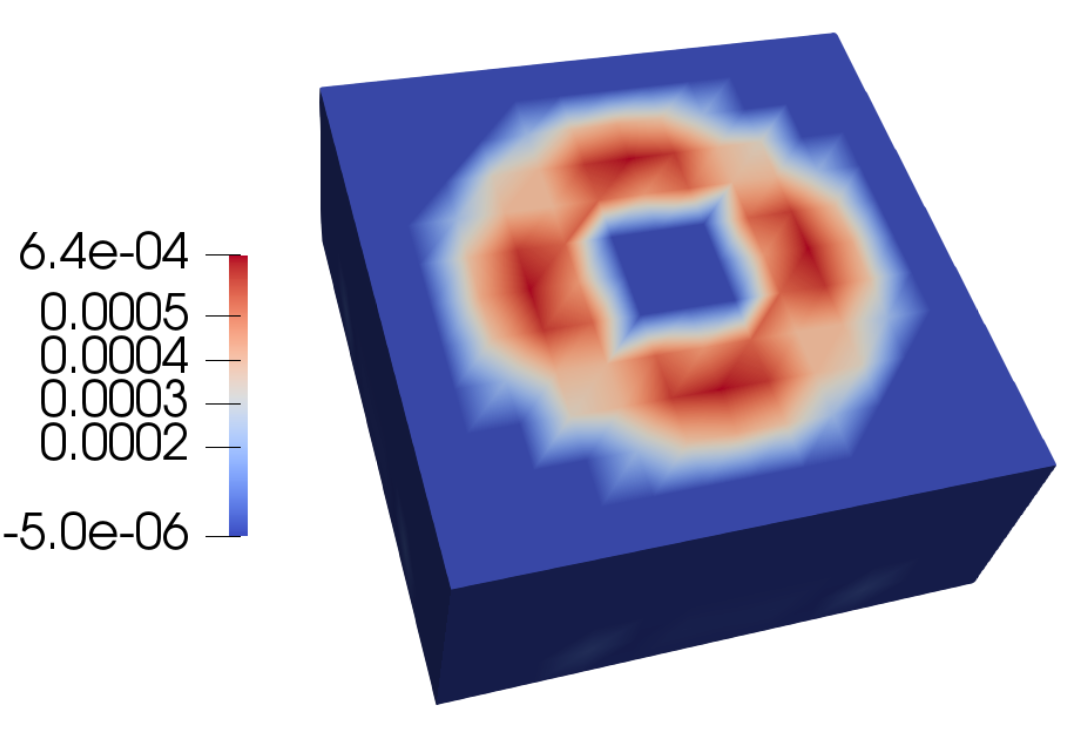}}
\subfigure[$\Tilde{\boldsymbol{u}}_3^{(2)}(\X,300,31.5)$]{\includegraphics[scale=0.2]{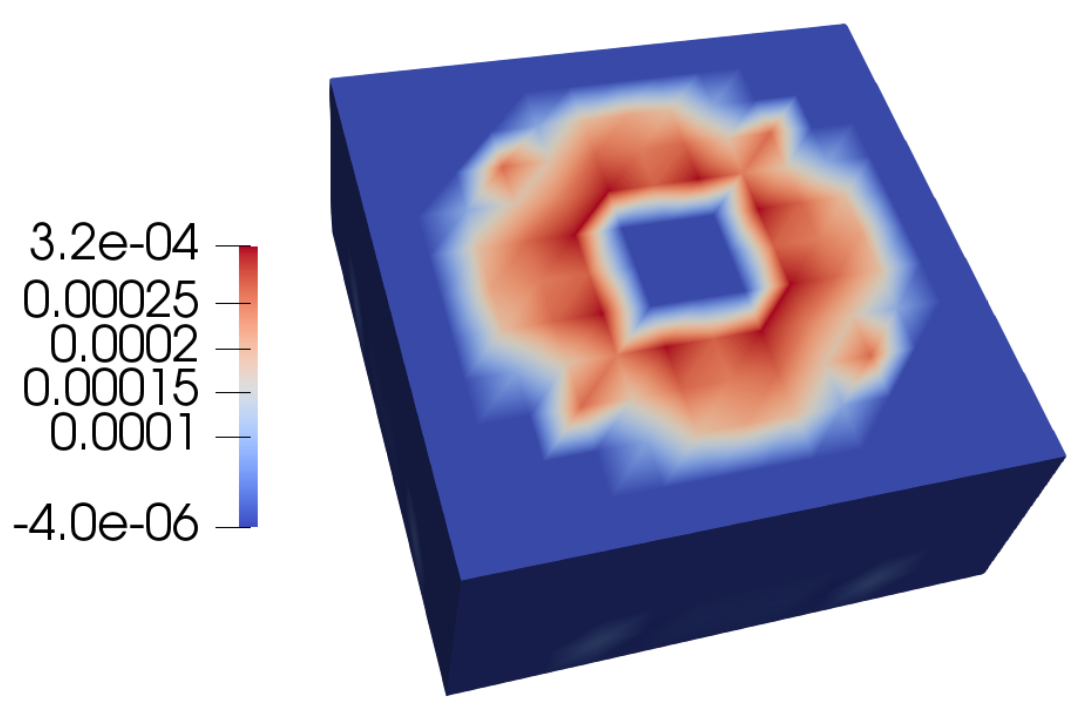}}
\subfigure[$\Tilde{\boldsymbol{u}}_3^{(2)}(\X,1380,52.5)$]{\includegraphics[scale=0.2]{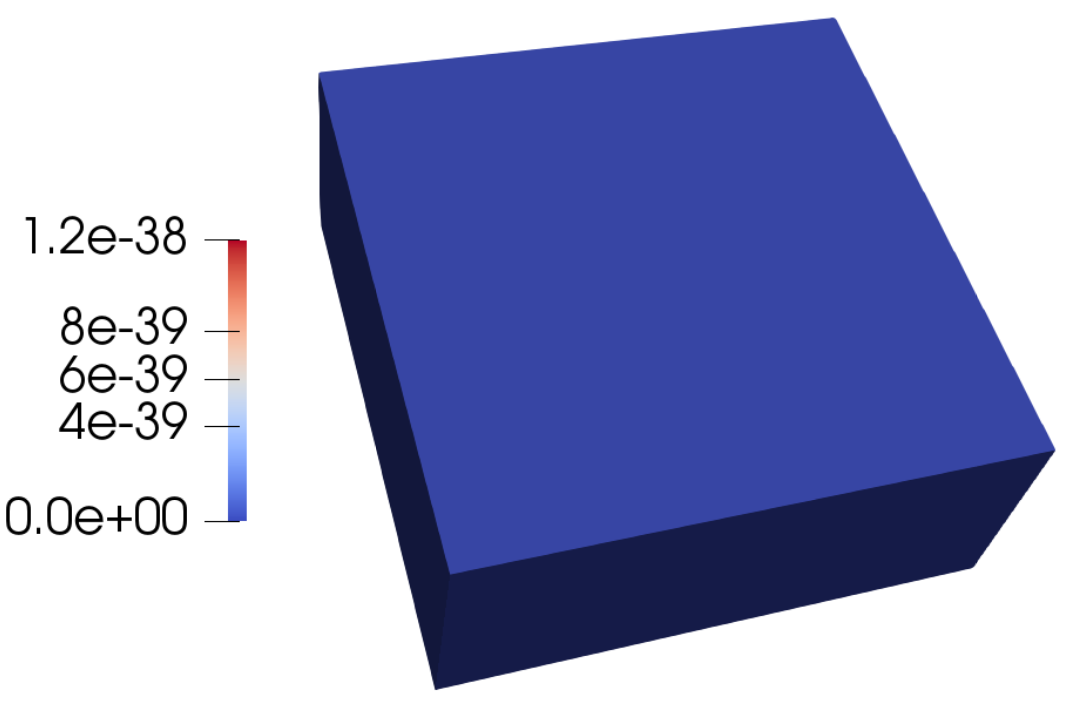}}
\subfigure[$\Tilde{\boldsymbol{u}}_3^{(2)}(\X,1140,94.5)$]{\includegraphics[scale=0.2]{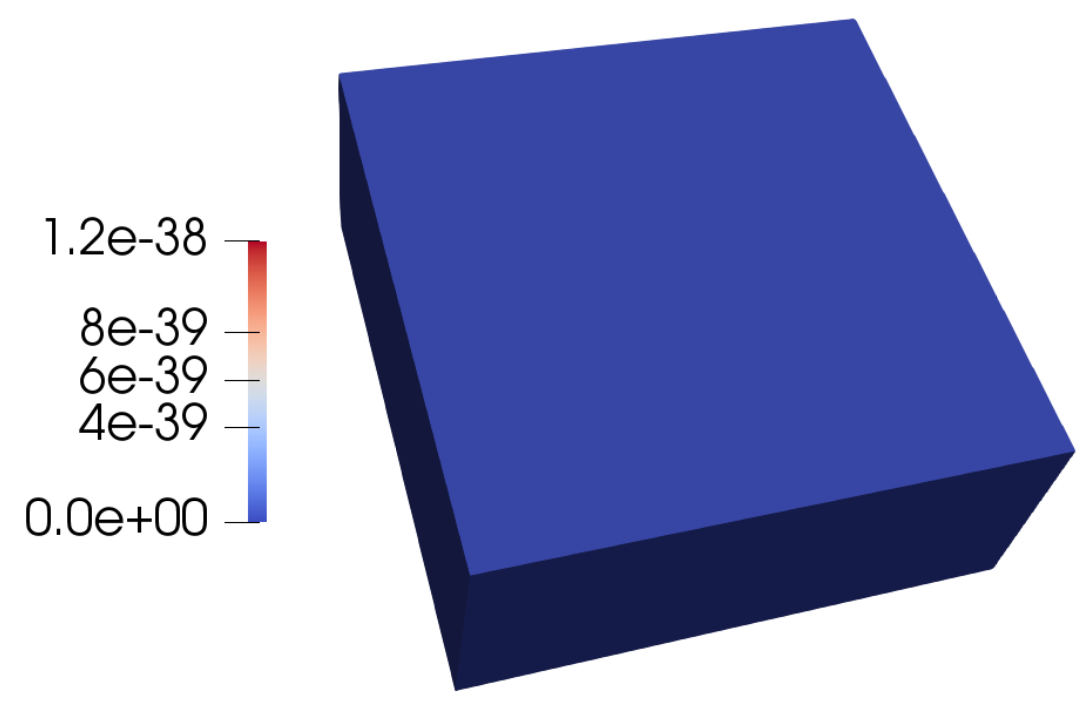}}
\caption{XTD extended modes $\Tilde{\boldsymbol{u}}^{(k)}(\X,{\sigma_y},H)$ for test case 2, time step $t_5$, $z$-component $\Tilde{\boldsymbol{u}}_3^{(k)}(\X,{\sigma_y},H)$. For example, $\Tilde{\boldsymbol{u}}_3^{(1)}(\X,360,21)$ stands for the first mode with $\sigma_y=360$ MPa, $H=21$ GPa.}
\label{fig:ExC2t5modenonsep}
\end{figure}

\begin{figure}[htbp]
\centering
\subfigure[$g_1^{(m)}(\sigma_y)$ ]{\includegraphics[scale=0.28]{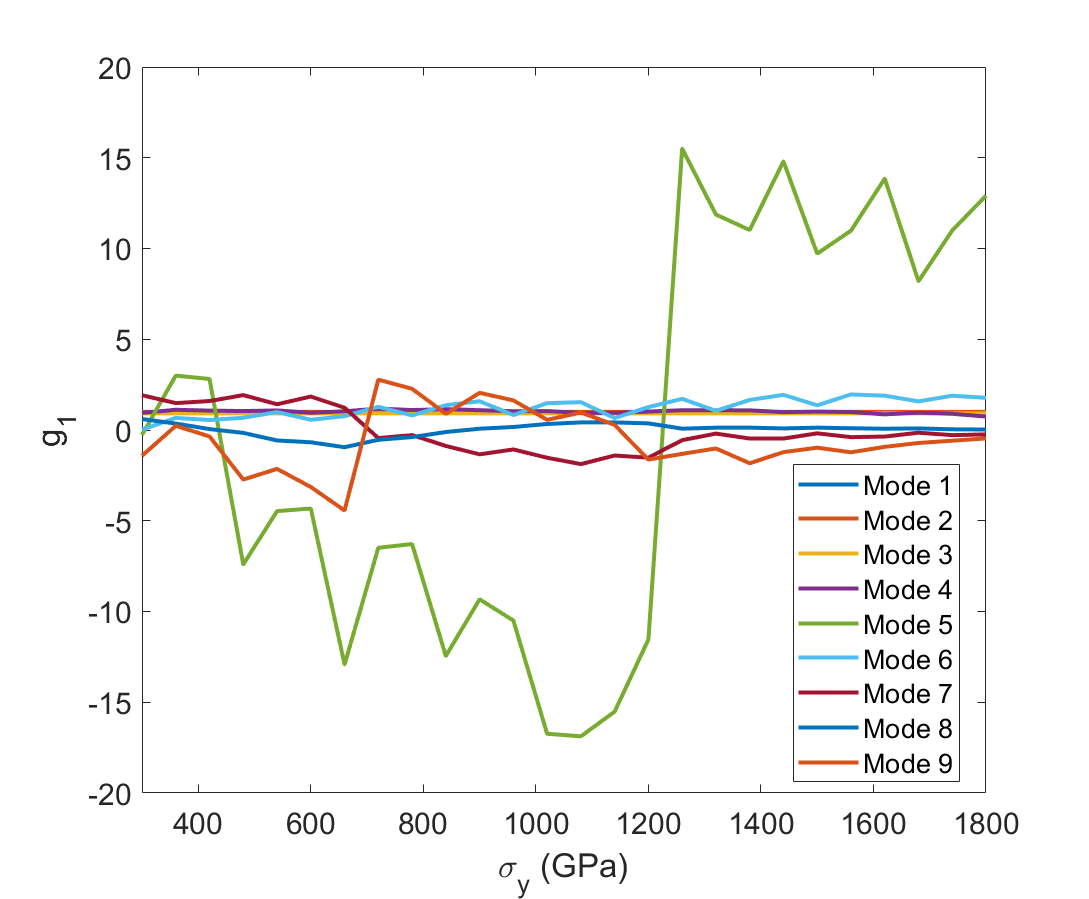}}
\subfigure[$g_2^{(m)}(H)$ ]{\includegraphics[scale=0.28]{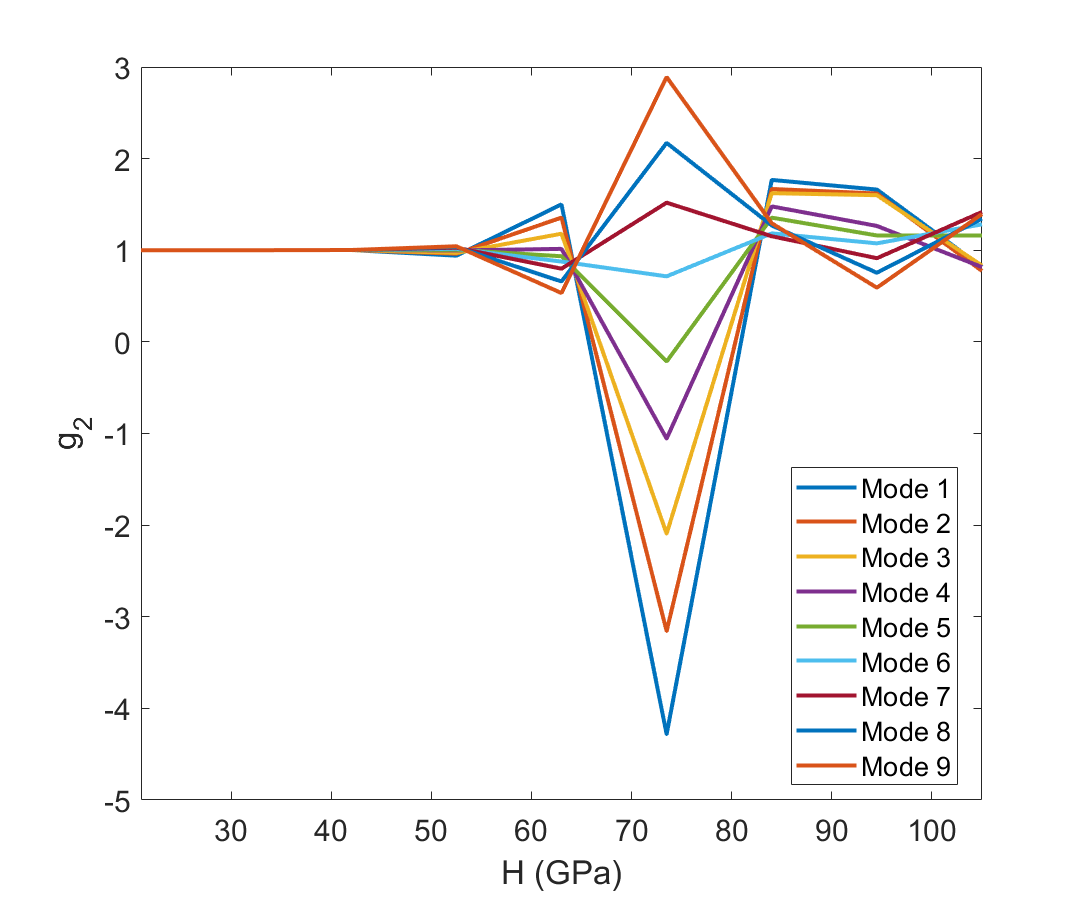}}
\caption{XTD separated modes $g_1^{(m)}(\sigma_y), \ g_2^{(m)}(H)$  for test case 2, time step $t_5$.}
\label{fig:ExC2t5modeg1g2}
\end{figure}

\begin{figure}[htbp]
\centering
{\includegraphics[scale=0.45]{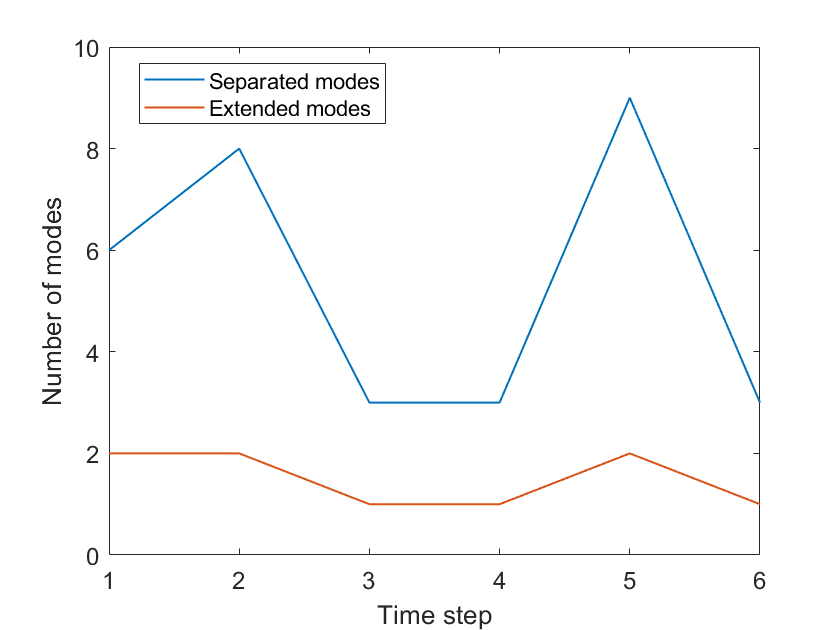}}
\caption{The number of modes for each time step in the test case 2}
\label{fig:ex1nbmode}
\end{figure}

\figurename~\ref{fig:ExC2t5modesep} illustrates the different modes of $\boldsymbol{a}^{(m)}(\X)$  obtained by XTD for the test case 2 at the time step 5. In particular, the component along the $z$-direction is plotted, which is the most important one in this example. It can be seen that the first mode is related to the boundary conditions and the loading displacement. The following modes are more and more distributed in the spatial domain.  The components of the non-separated extended modes are illustrated in \figurename~\ref{fig:ExC2t5modenonsep}. It confirms that this extended term is sparse in the global domain, as most of the region remains zero. This enables an efficient computation for the non-separated term, as mentioned previously. The other separated modes $g_1^{(m)}(\sigma_y)$ and $\ g_2^{(m)}(H)$ are plotted in \figurename~\ref{fig:ExC2t5modeg1g2}.    They are 1D functions that represent the variation of the solutions in a single parameter direction. It can be seen that the solution variation is not necessarily smooth in the parameter space, which motivated us to develop the XTD method. The total number of modes for each time step is reported in \figurename~\ref{fig:ex1nbmode}. The small number of modes for each time step confirms that the displacement increment is a good choice for the separated variable. For the 7 test cases, we observed similar results. Thus we illustrated the test case 2 as an example. 

This first numerical experiment investigated the performance of the proposed XTD method. It confirms its efficiency, in particular, for problems with localized loading. In the next example, we will study the application to problems with a localized moving load, which is relevant to  additive manufacturing processes.

\subsection{Application to thermal residual stress prediction}
\subsubsection{Offline training}
Additive manufacturing or welding is controlled by many uncertain parameters related to the materials, heat source, and boundary conditions. The parametric study for such a process fits well the aim of the overall XTD framework. In this example, we will particularly investigate the thermally induced residual stress with respect to the variation of two parameters, again the initial yield stress $\sigma_y$ and the  hardening coefficient $H$. The detailed material properties and the variation range will be described later. Assuming the temperature of each time step is given and the thermal and mechanical problems are weakly coupled, the XTD method  computes the displacement increment as below

\begin{equation}
\displaystyle
\label{eq:XTDex2}
\Delta\ud(\X,\sigma_y,H)= \sum_{m=1}^{M}\boldsymbol{a}^{(m)}(\X)g_1^{(m)}(\sigma_y)g_2^{(m)}(H) + \sum_{k=1}^{K}\Tilde{\boldsymbol{u}}^{(k)}(\X,{\sigma_y},H)
\end{equation}

\begin{figure}[htbp]
\centering
{\includegraphics[scale=0.3]{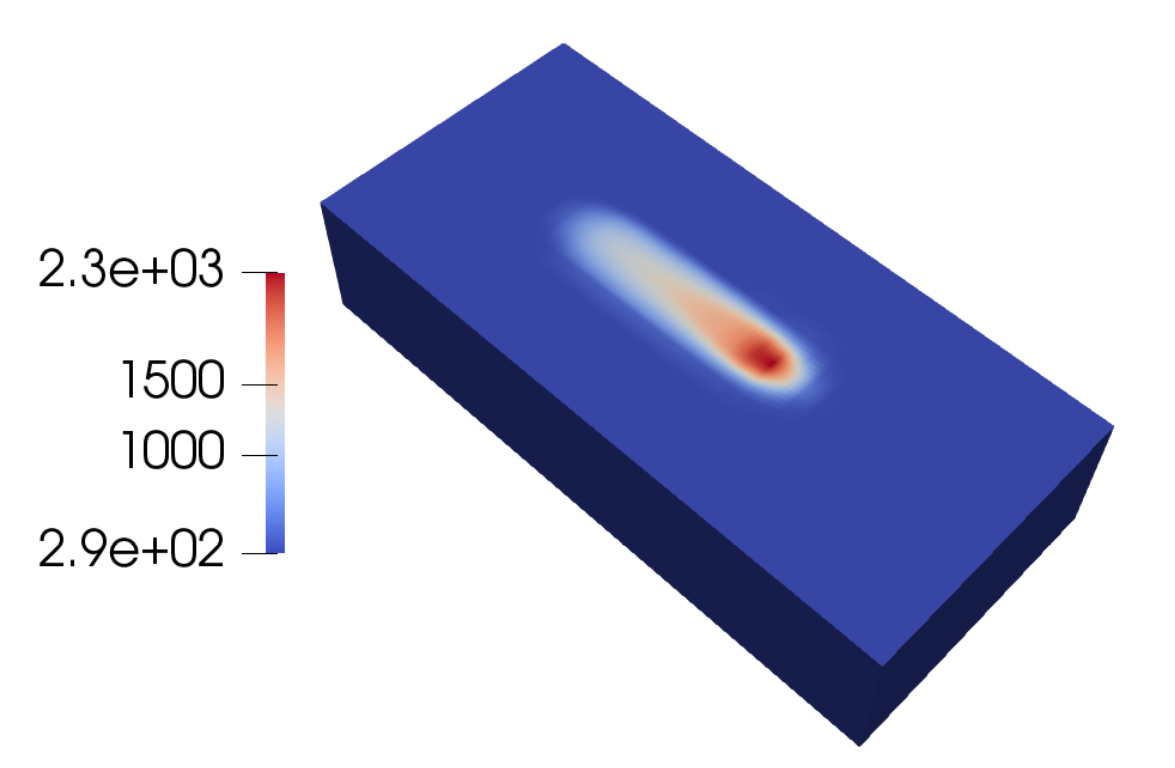}}
\caption{Snapshot of the temperature profile. Domain size: $2.7$ mm $\times 1.25$ mm $\times 0.6$ mm.}
\label{fig:ex2T}
\end{figure}

In our work, we first used the  hyper reduction method \cite{lu2020adaptive} for the thermal fluid analysis to compute the temperature evolution in a predefined domain. The details of the model can be found in \cite{lu2020adaptive}. {We remark that this choice is only for demonstration purposes. In general, any code that can provide the temperature evolution can be used without special difficulties.} Some key  parameters for the thermal computations are summarized in \tablename~\ref{table:ex2thermalprocess}. \figurename~\ref{fig:ex2T} illustrates a snapshot of the output temperature profile on the computation domain. We selected 30 output steps as the input for the subsequent mechanical analysis. The resulting stress is considered by the thermal expansion as follows
\begin{equation}
\displaystyle
\label{eq:stressT}
\sig =\textbf{D}:\epse =\textbf{D}:(\eps -\epsp-\eps_{\theta})
\end{equation}
with
\begin{equation}
\displaystyle
\label{eq:expansion}
\eps_{\theta}=\alpha\textbf{I}\theta
\end{equation}
where $\theta$ denotes the variation of  temperature and $\textbf{I}$ is the identity matrix. In real implementation, this can be performed by applying an equivalent external force induced by the thermal expansion.

\begin{table}[htbp]
\caption{Parameters used for the thermal  fluid analysis}
\centering
\begin{tabular}{|c|c|c|c|c|c|}
\hline
Power & Scan speed &  Step size & Total time steps& Output steps& Materials\\ \hline
195 W & 0.8 m/s& $10^{-5}$ s & 145 &30& IN625\\ \hline
\end{tabular}
\label{table:ex2thermalprocess}
\end{table}

In the mechanical analysis, the bottom surface is fixed as the boundary condition. The elastic-plastic materials properties are considered temperature-dependant, as shown in \tablename~\ref{table:ex2mate} and reported in \cite{denlinger2015thermo}. Considering these reported properties as reference values, the variations of the temperature dependant yield stress $\sigma_y (\theta)$ and the  hardening coefficient $H(\theta)$ can be  considered by two amplitude parameters $A_1$ and $A_2$, i.e. $\sigma_y (\theta)\leftarrow A_1 \sigma_y (\theta) $ with $A_1\in [1,\ 2]$ and $H (\theta)\leftarrow A_2 H (\theta) $ with $A_2\in [1,\ 3]$. 

\begin{table}[htbp]
\caption{Temperature-dependant mechanical properties}
\centering
\begin{tabular}{|c|c|c|c|c|c|}
\hline
Temperature (K) & $\alpha$ ($\times 10^{-6}$ K$^{-1}$) &  $\sigma_y$ (MPa) & $H$ (GPa)& Young's modulus (GPa)\\ \hline
273 & 12.8& 493 & 20.8 &208\\ \hline
293 & 12.8& 493 & 20.8 &208\\ \hline
366 & 12.8& 479 & 20.4 &204\\ \hline
478 & 13.1& 443 & 19.8 &198\\ \hline
588 & 13.3& 430 & 19.2 &192\\ \hline
698 & 13.7& 424 & 18.6 &186\\ \hline
813 & 14.0& 423 & 17.9 &179\\ \hline
923 & 14.8& 422 & 17.0 &170\\ \hline
1033 & 15.3& 415 & 16.1 &161\\ \hline
1143 & 15.8& 386 & 14.8 &148\\ \hline
$>$ 1143 & 15.8& 386 & 14.8 &148\\ \hline
\end{tabular}
\label{table:ex2mate}
\end{table}

Concerning the discretization, a mesh of size $36\times22\times9$ is used for the physical (spatial) domain. The total number of time steps is 30. For the parameter domain, the $A_1$ is discretized by 6 equally distancing points and  $A_2$ is discretized by 5 points. Hence, if the FEA is used for computing the space-time-parametric solutions, we need to run 30 ($6\times 5$) times the nonlinear space-time simulation. Each of this FEA takes around 1100 s. The estimated cost for the 30 runs is then 33000 s.

\begin{figure}[htbp]
\centering
\subfigure[FEA, $t_{30}$ ]{\includegraphics[scale=0.28]{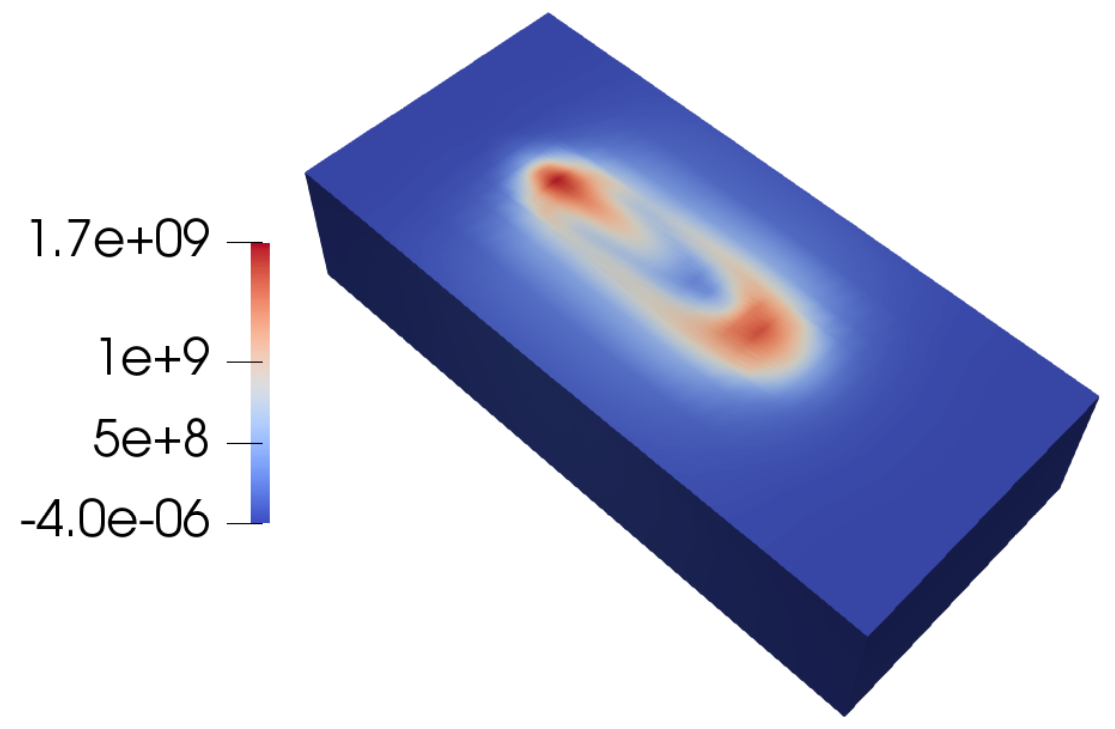}}\quad\quad
\subfigure[XTD, $t_{30}$ ]{\includegraphics[scale=0.28]{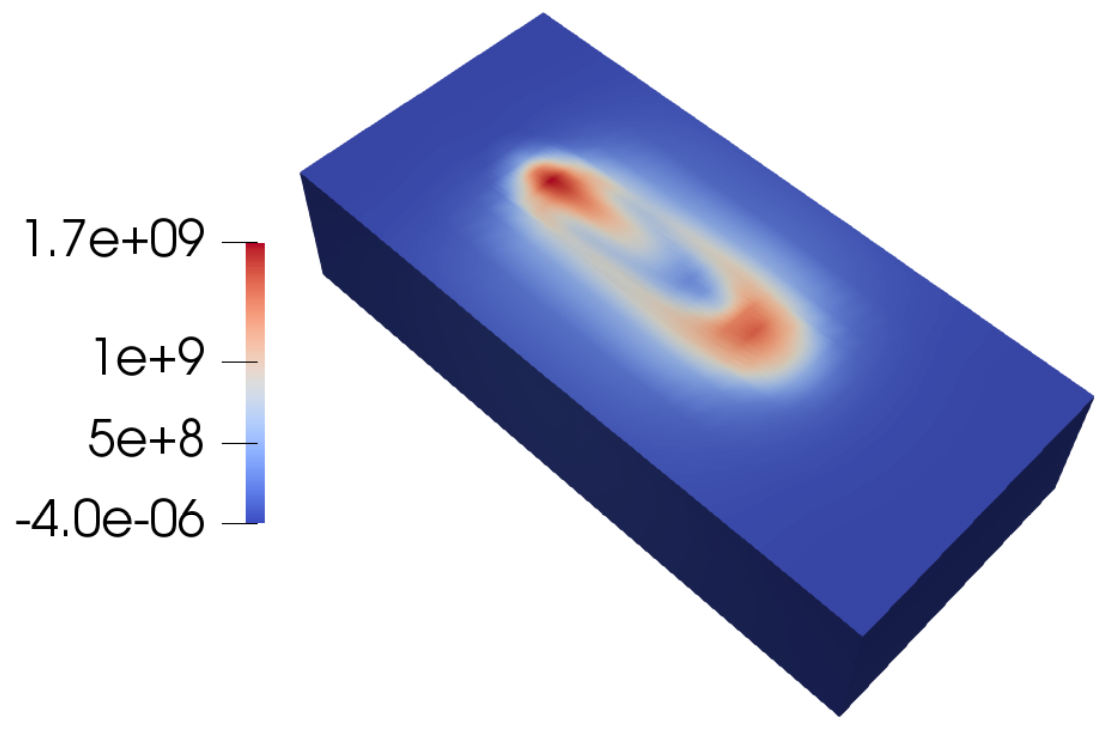}}\\
\subfigure[FEA, $t_{30}$ ]{\includegraphics[scale=0.28]{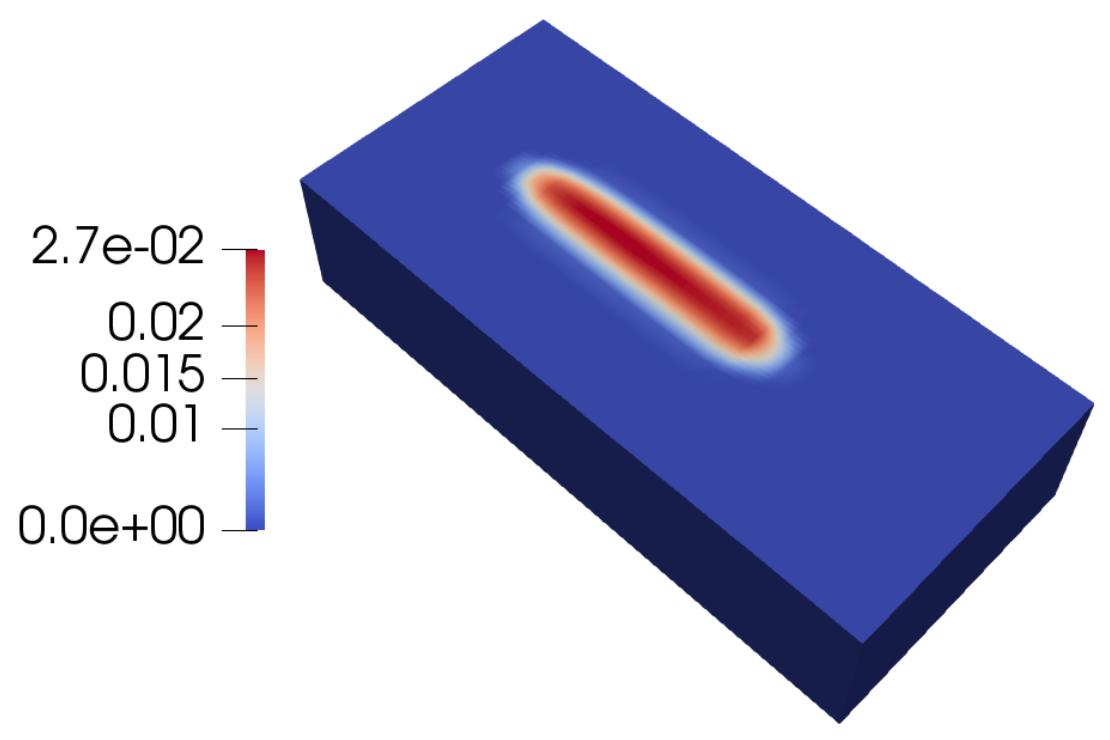}}
\quad\quad\subfigure[XTD, $t_{30}$ ]{\includegraphics[scale=0.28]{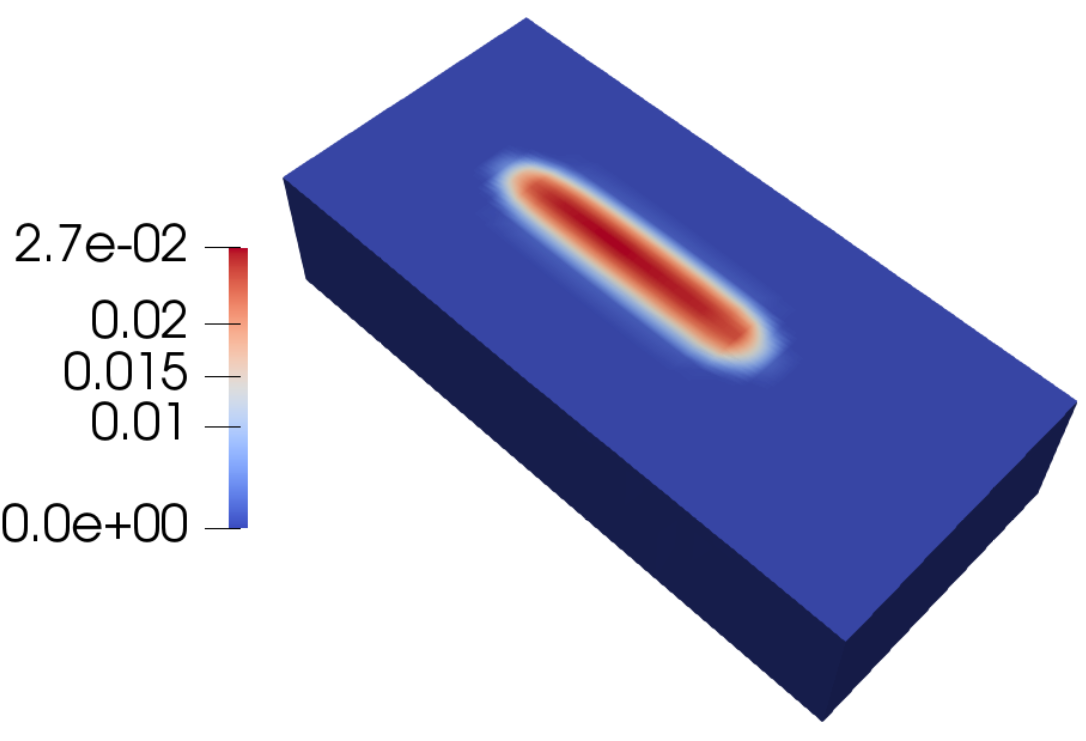}}
\caption{Comparison of FEA and XTD solutions for the parameters: $A_1=1.9$, $A_2=2.1$. (a)(b) von Mises stress stress, (c)(d) equivalent plastic strain.}
\label{fig:Ex2FEXTD}
\end{figure}

The XTD model \eqref{eq:XTDex2} is used for computing the parametric solution. The total cost is 10480 s, which is significantly less than that of FEA. This confirms the efficiency of the method.  Again, this is the offline computational cost, the online prediction is very fast ($< 0.01$ s). \figurename~\ref{fig:Ex2FEXTD} compares the online XTD prediction to the FEA reference solution for the parameters $A_1=1.9$, $A_2=2.1$. They show a very good agreement in terms of stress and plastic strain. The relative  $L^2$ norm difference is less than $2\%$.

\begin{figure}[htbp]
\centering
\subfigure[XTD, $t_{5}$ ]{\includegraphics[scale=0.2]{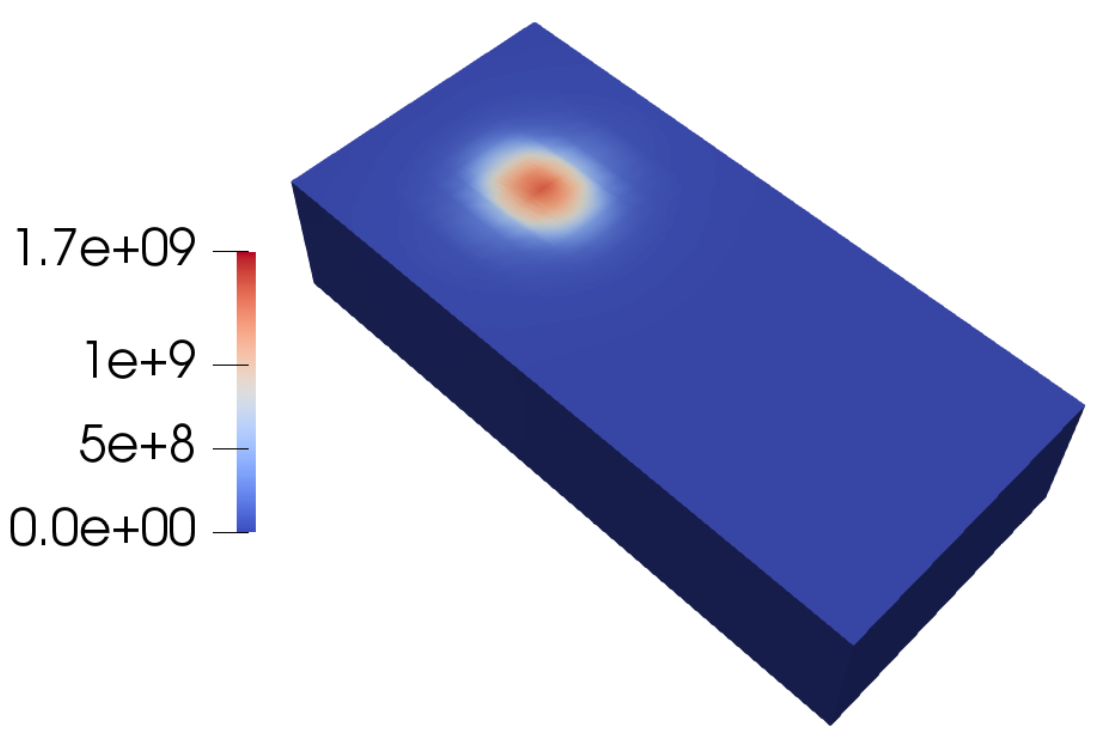}}
\subfigure[XTD, $t_{10}$ ]{\includegraphics[scale=0.2]{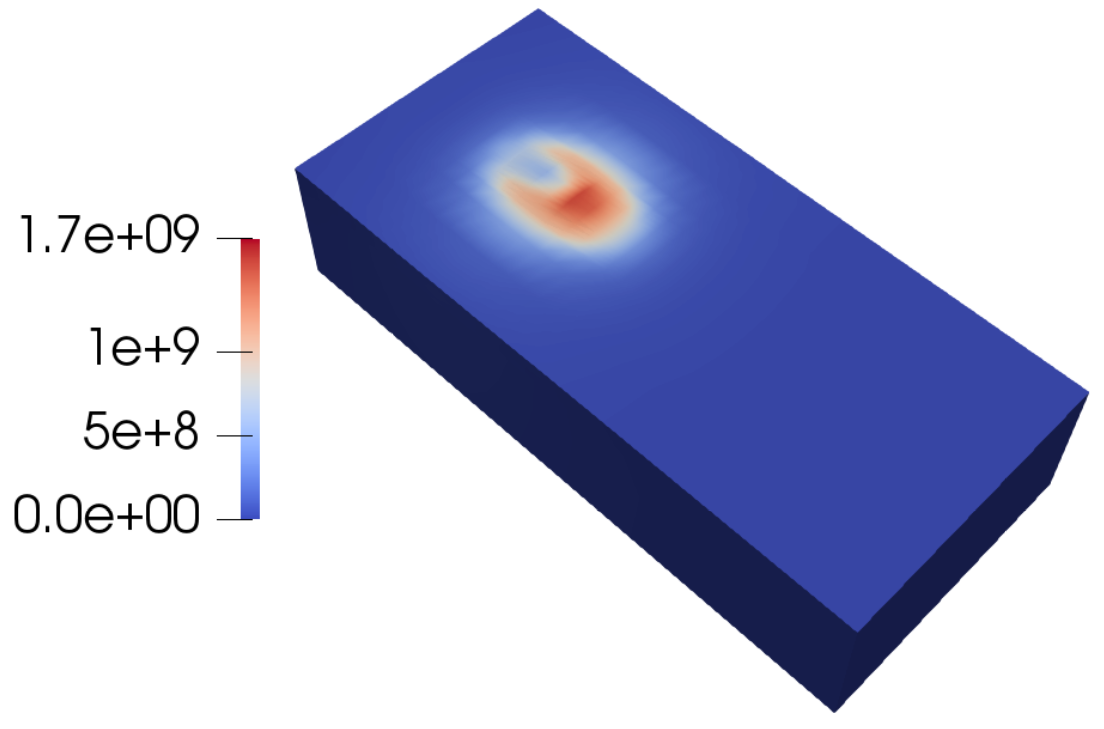}}
\subfigure[XTD, $t_{20}$ ]{\includegraphics[scale=0.2]{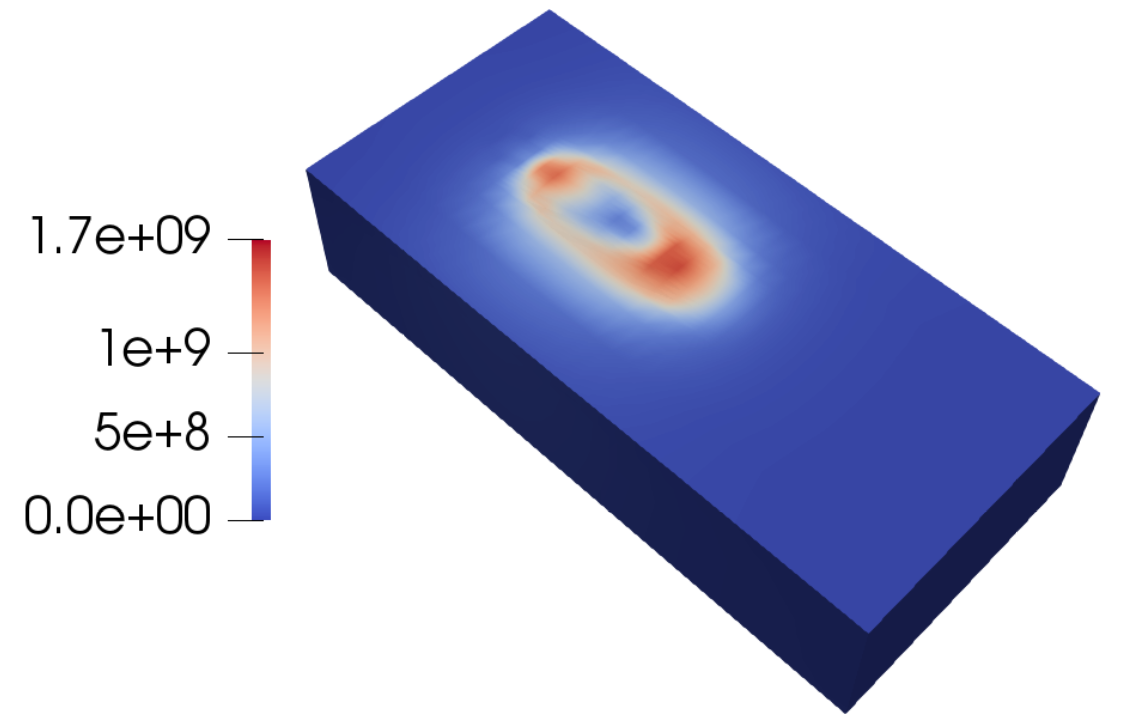}}
\subfigure[XTD, $t_{30}$ ]{\includegraphics[scale=0.2]{img/Ex2XTDVMt30.png}}\\
\subfigure[XTD, $t_{5}$ ]{\includegraphics[scale=0.2]{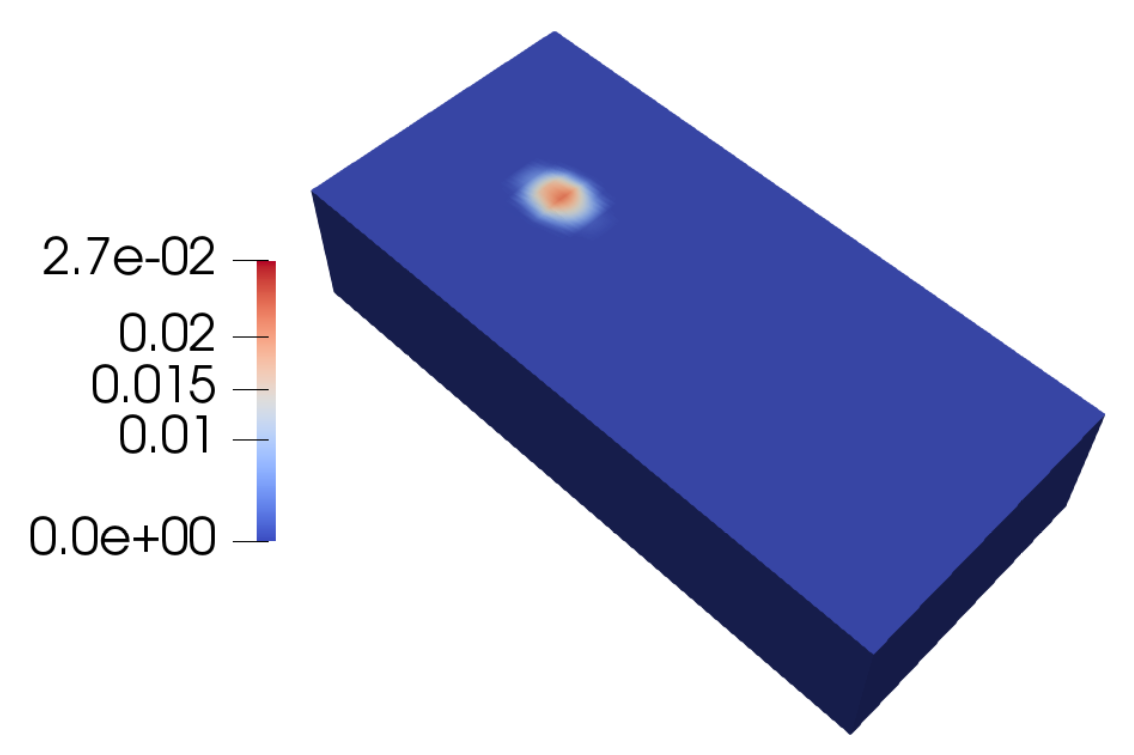}}
\subfigure[XTD, $t_{10}$ ]{\includegraphics[scale=0.2]{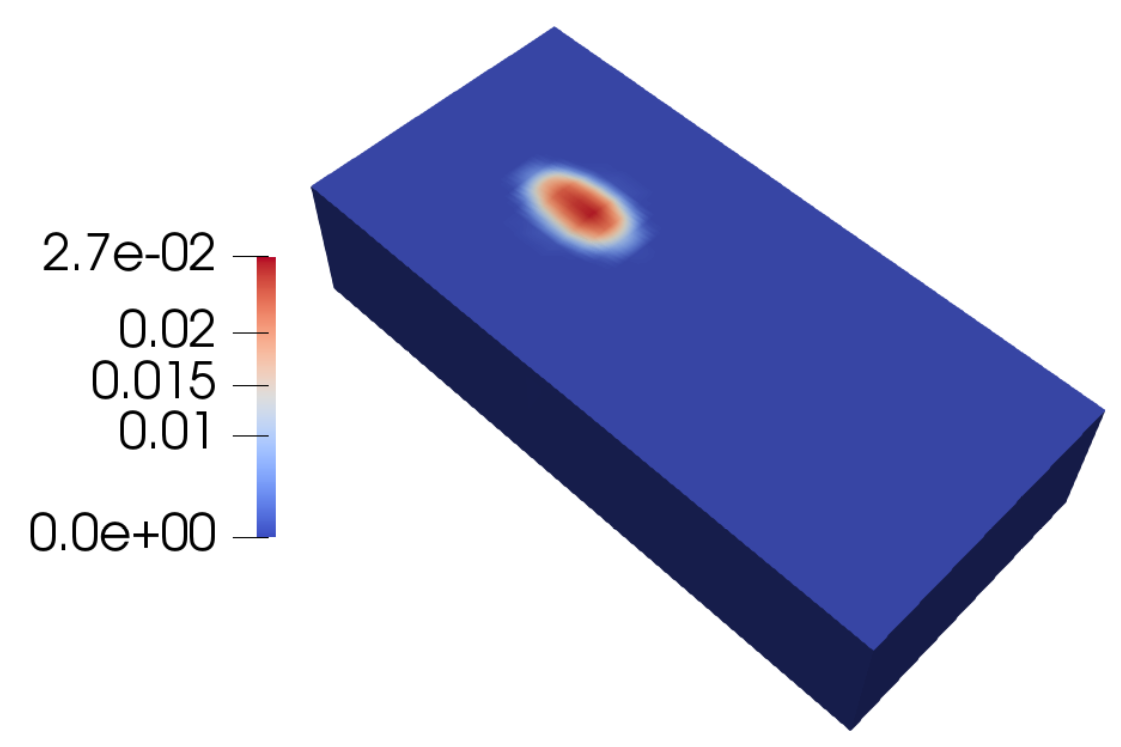}}
\subfigure[XTD, $t_{20}$ ]{\includegraphics[scale=0.2]{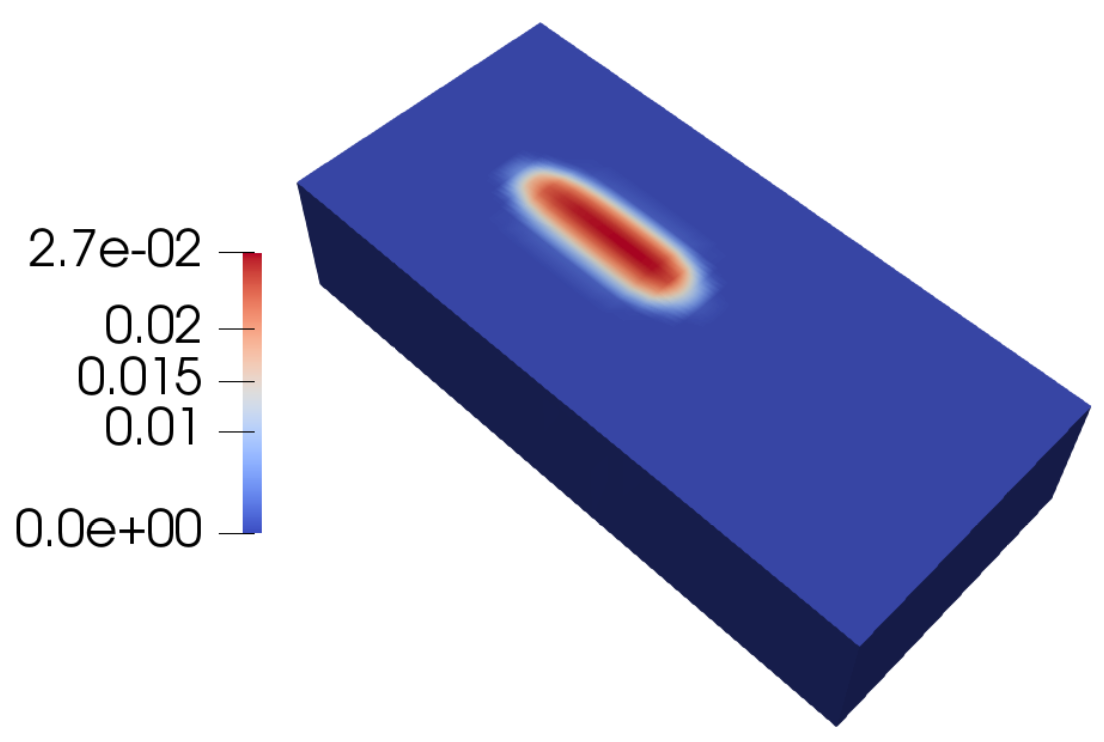}}
\subfigure[XTD, $t_{30}$ ]{\includegraphics[scale=0.2]{img/Ex2XTDpt30.png}}
\caption{Evolution of stress and strain fields for the parameters: $A_1=1.9$, $A_2=2.1$. (a)-(d) von Mises stress stress, (e)-(h) equivalent plastic strain.}
\label{fig:Ex2XTD}
\end{figure}

\begin{figure}[htbp]
\centering
\subfigure[Mode 1 ]{\includegraphics[scale=0.2]{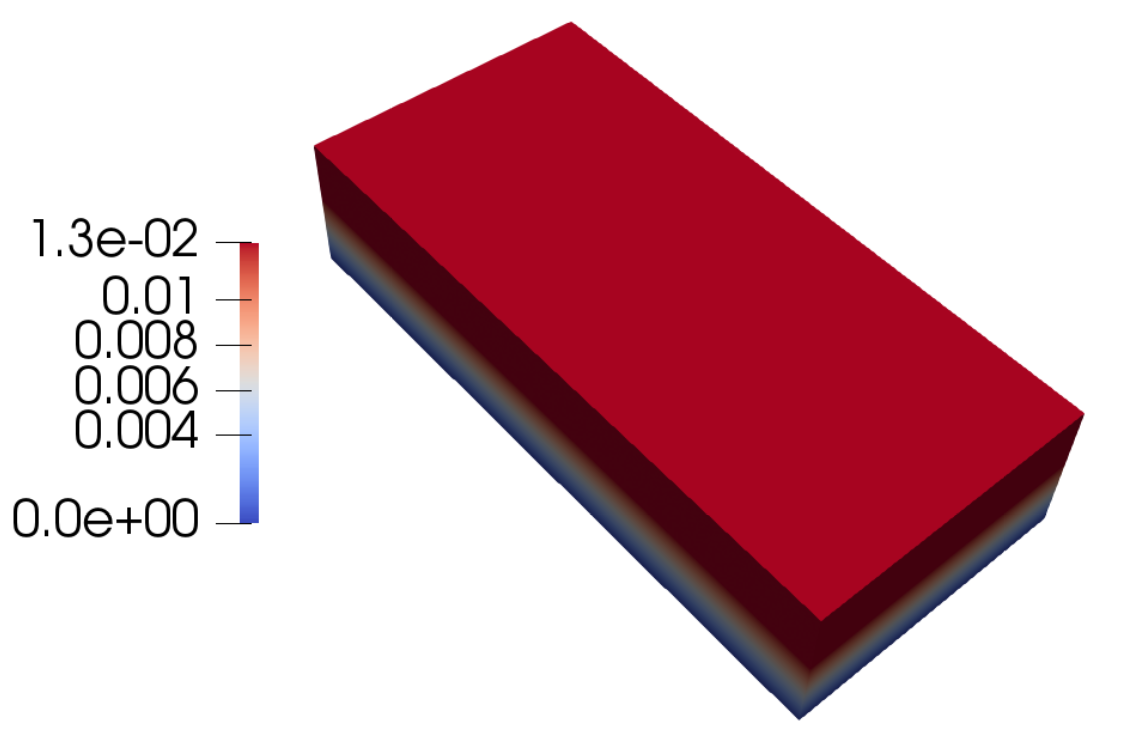}}
\subfigure[Mode 2 ]{\includegraphics[scale=0.2]{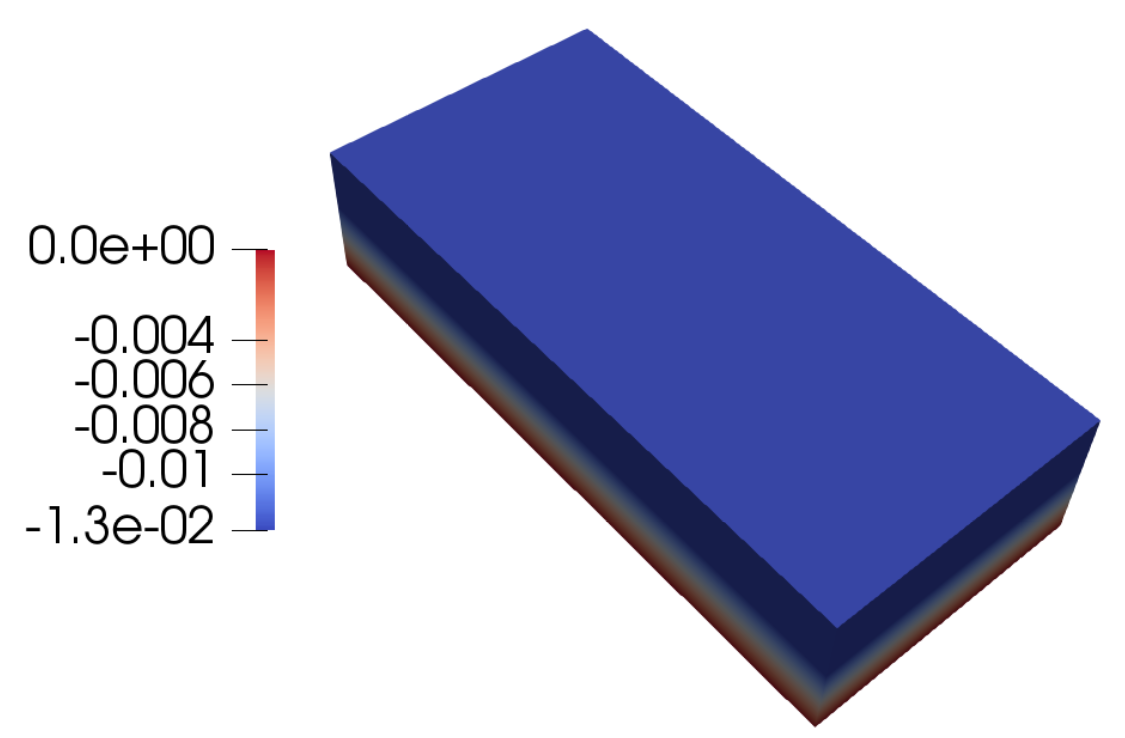}}
\subfigure[Mode 3 ]{\includegraphics[scale=0.2]{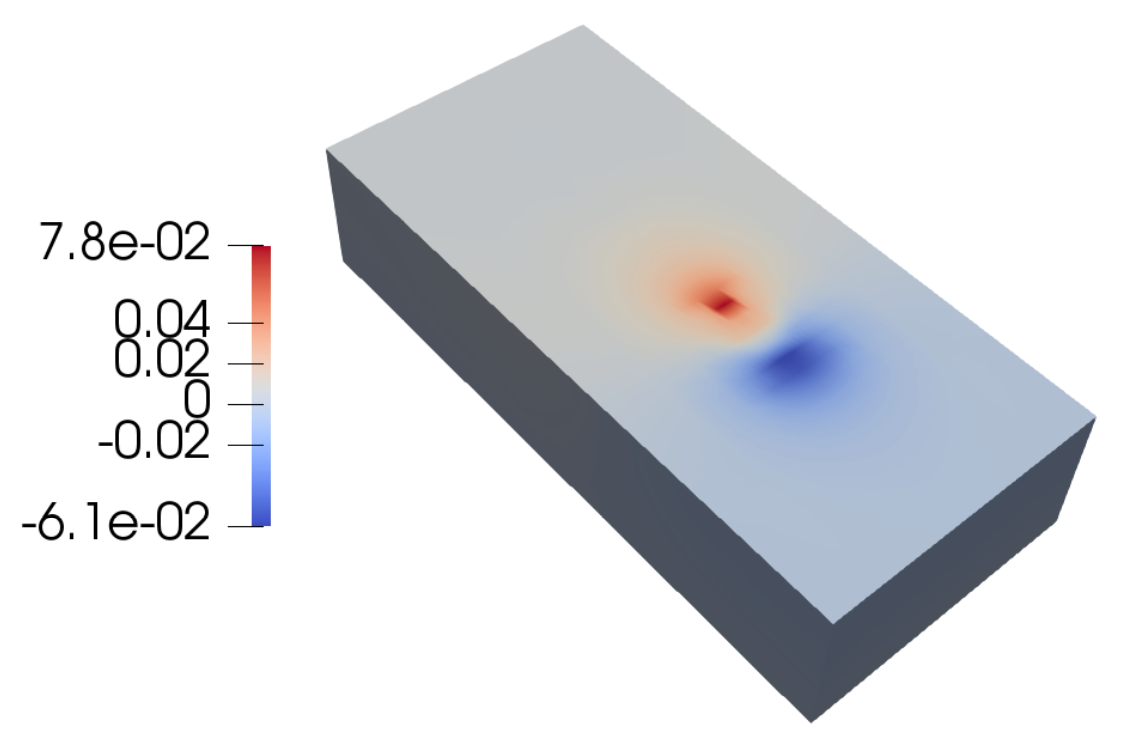}}
\subfigure[Mode 4 ]{\includegraphics[scale=0.2]{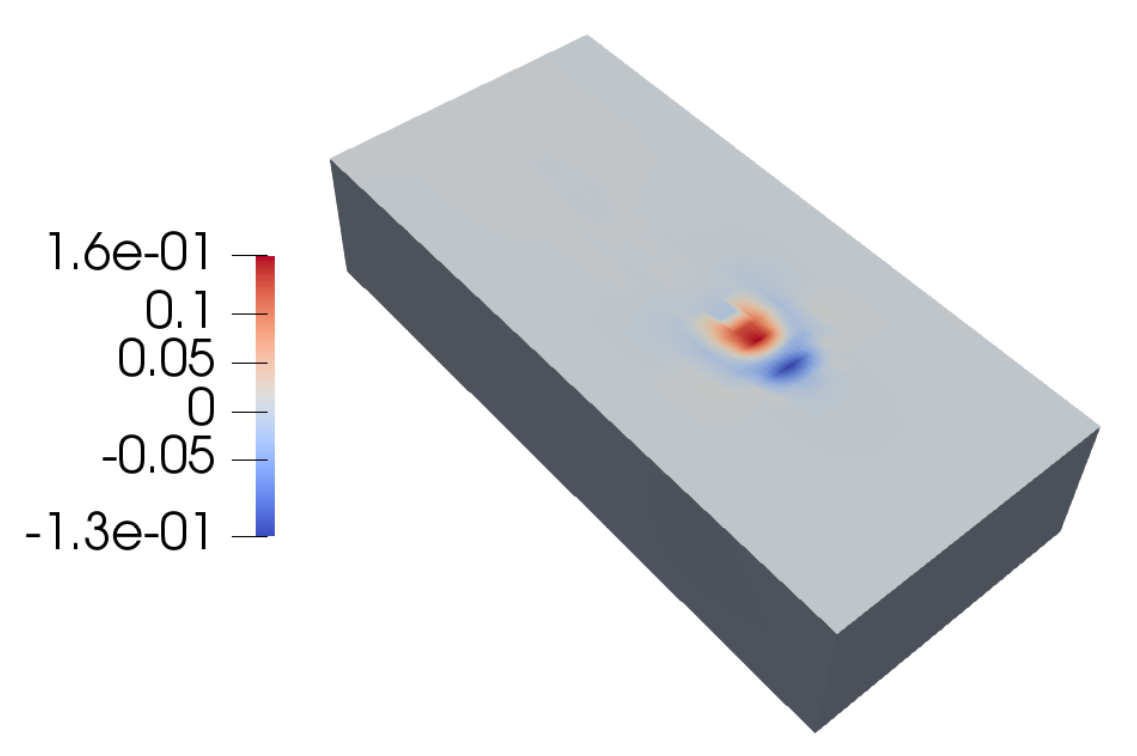}}\\
\subfigure[Mode 1 ]{\includegraphics[scale=0.2]{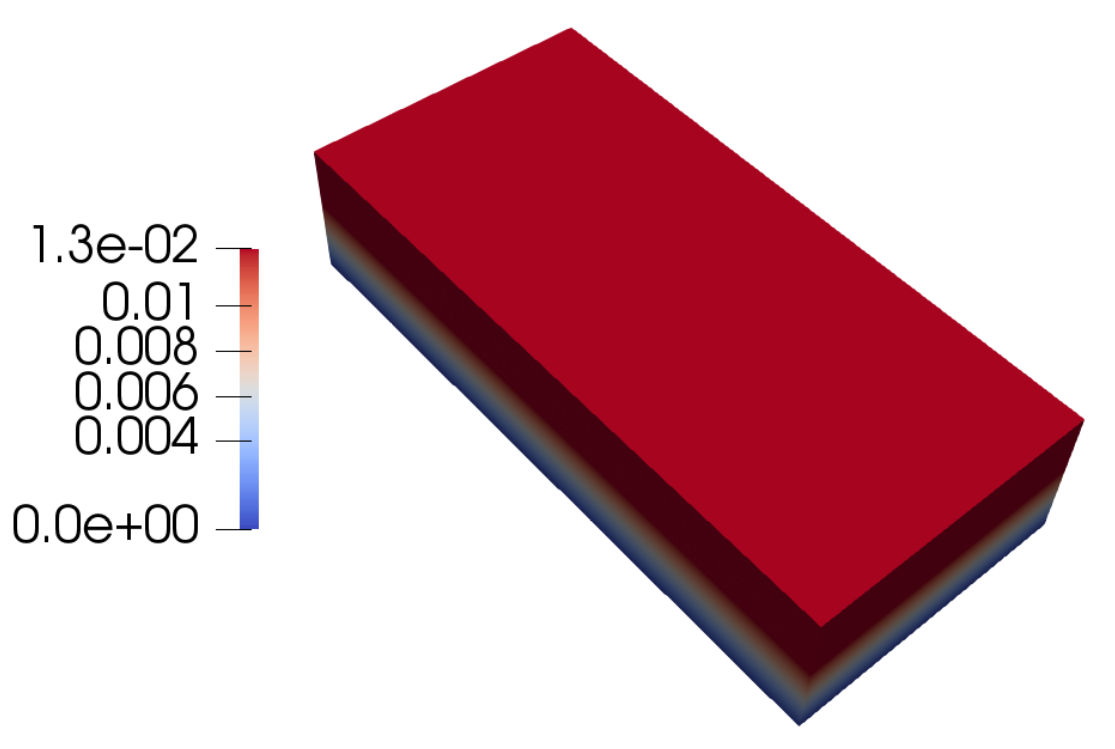}}
\subfigure[Mode 2 ]{\includegraphics[scale=0.2]{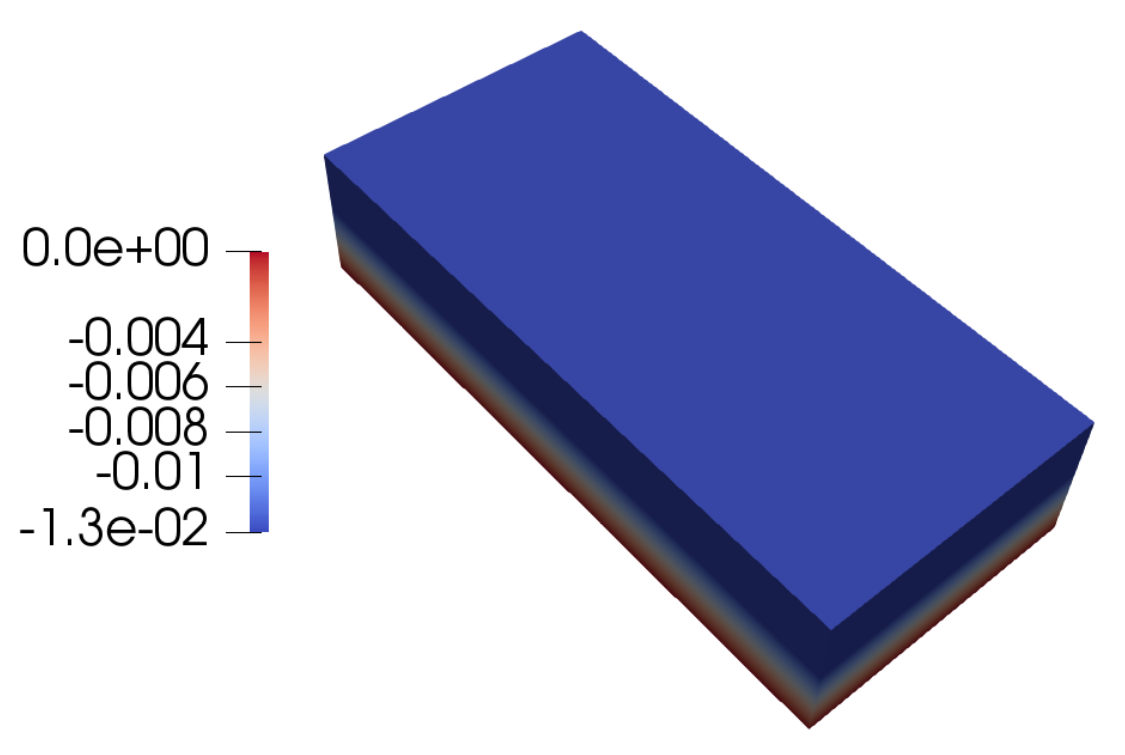}}
\subfigure[Mode 3 ]{\includegraphics[scale=0.2]{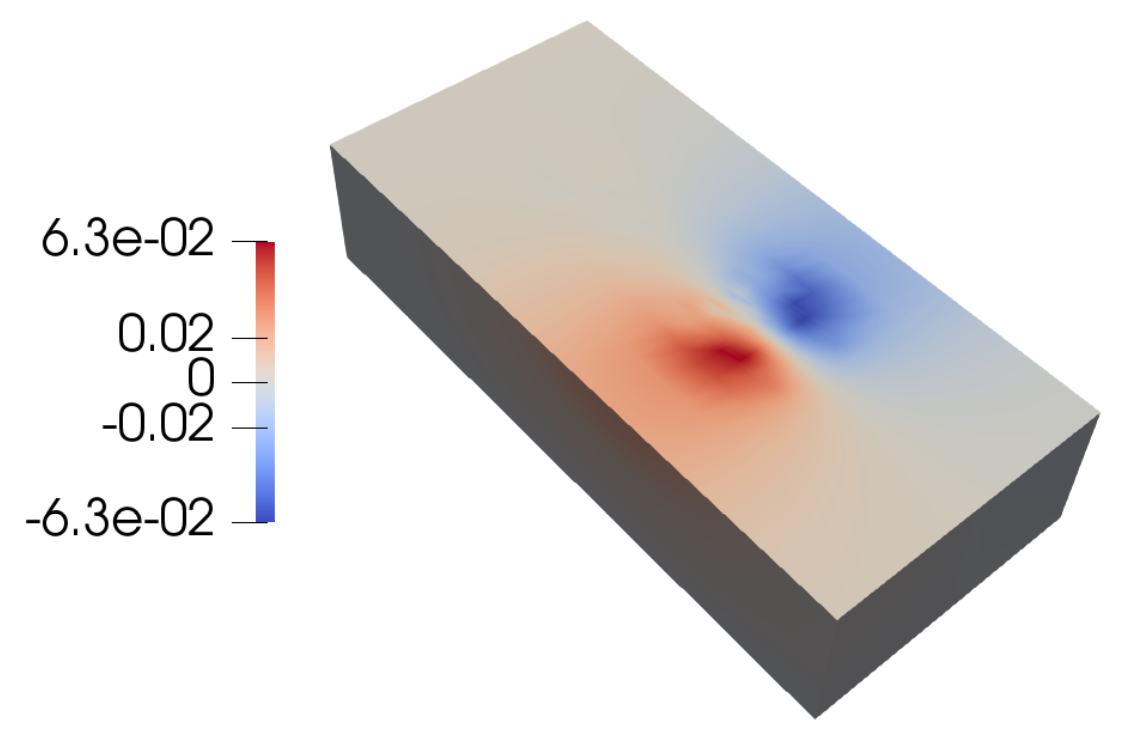}}
\subfigure[Mode 4 ]{\includegraphics[scale=0.2]{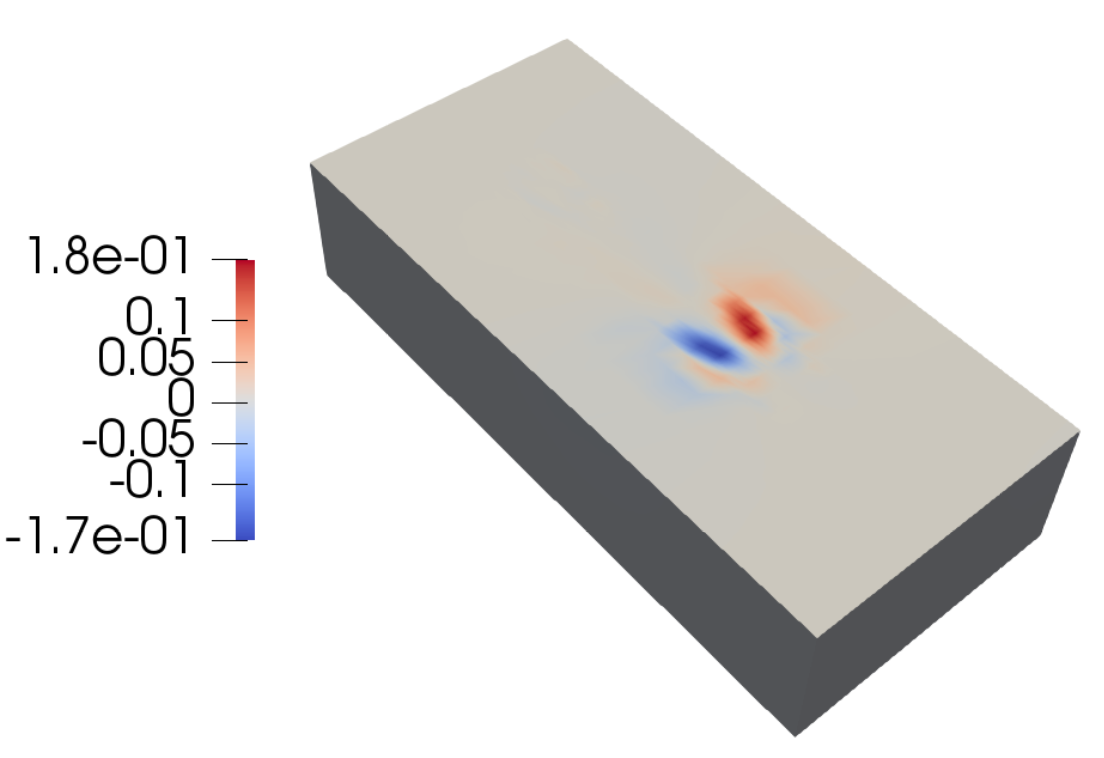}}\\
\subfigure[Mode 1 ]{\includegraphics[scale=0.2]{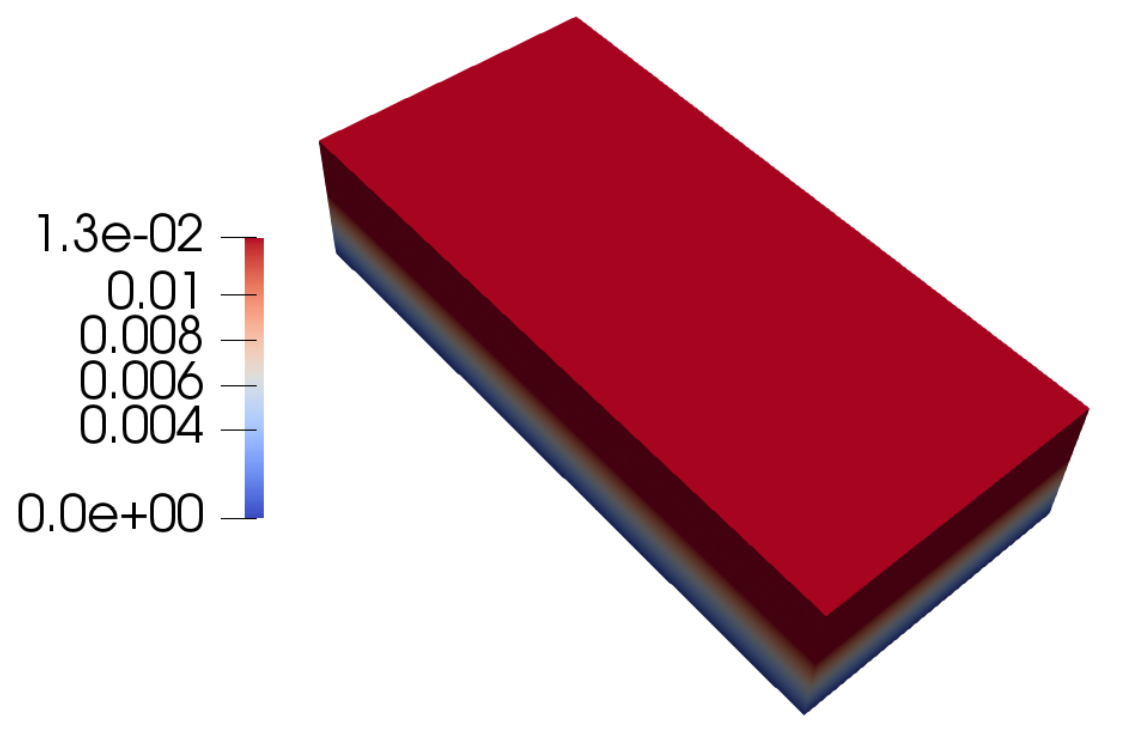}}
\subfigure[Mode 2 ]{\includegraphics[scale=0.2]{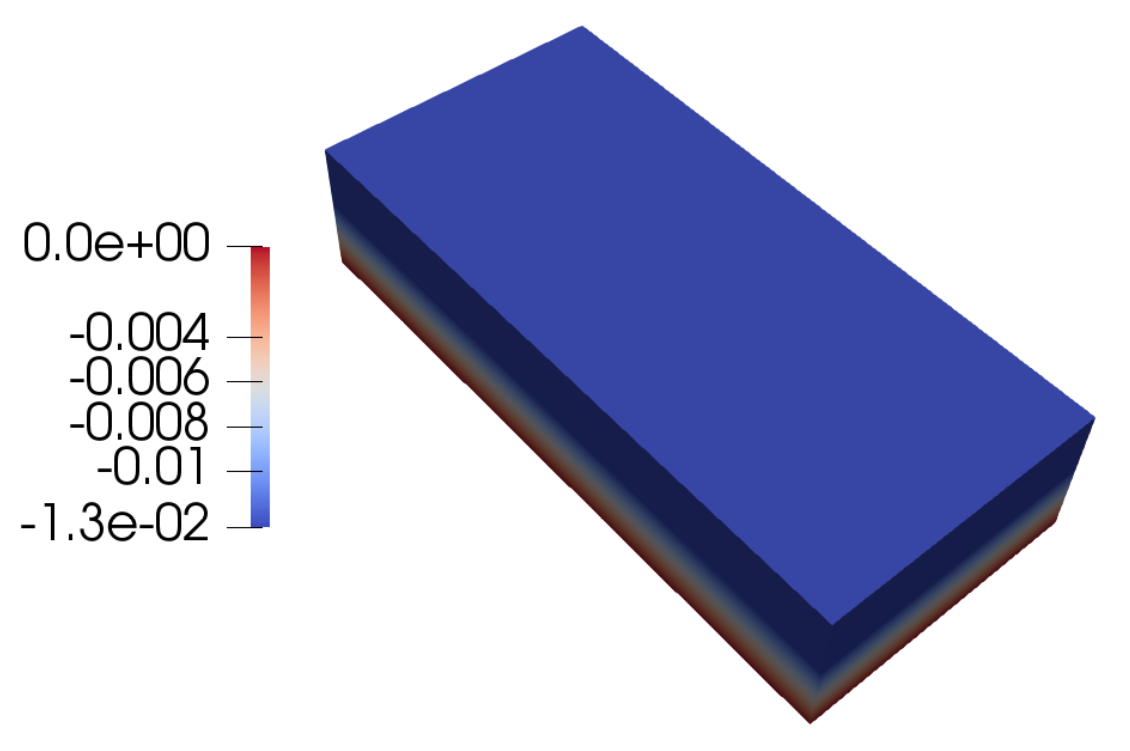}}
\subfigure[Mode 3 ]{\includegraphics[scale=0.2]{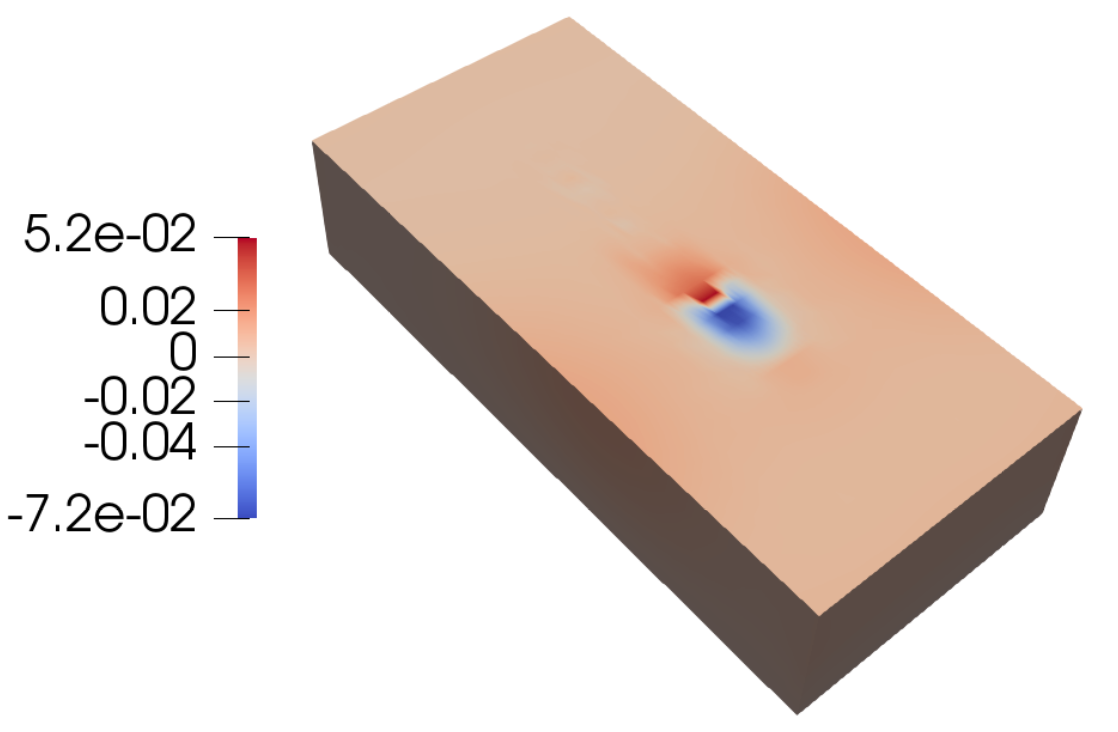}}
\subfigure[Mode 4 ]{\includegraphics[scale=0.2]{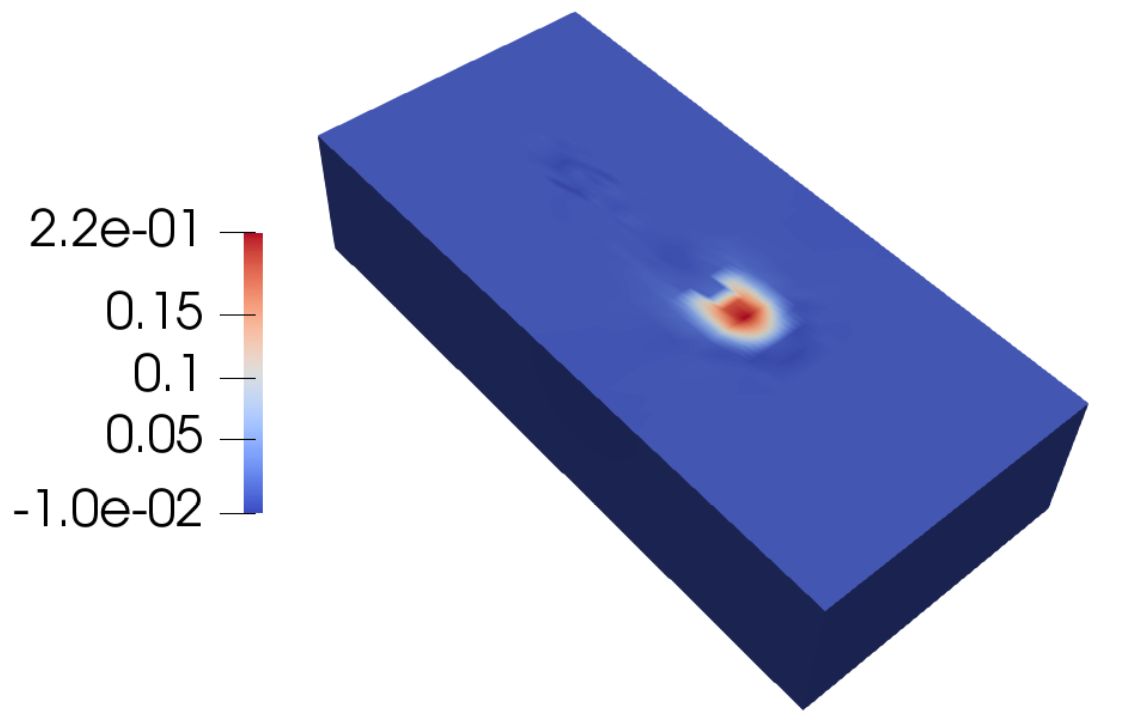}}\\
\caption{XTD  separated modes (normalized) $\boldsymbol{a}^{(m)}(\X)$, time step $t_{30}$.
(a)-(d) $x$-component $\boldsymbol{a}_1^{(m)}/\|\boldsymbol{a}_1^{(m)}\|$,
(e)-(h) $y$-component $\boldsymbol{a}_2^{(m)}/\|\boldsymbol{a}_2^{(m)}\|$,
(i)-(l) $z$-component $\boldsymbol{a}_3^{(m)}/\|\boldsymbol{a}_3^{(m)}\|$.}
\label{fig:Ex2XTDmode}
\end{figure}

\begin{figure}[htbp]
\centering
\subfigure[$A_1= 1.4,\ A_2= 1$ ]{\includegraphics[scale=0.2]{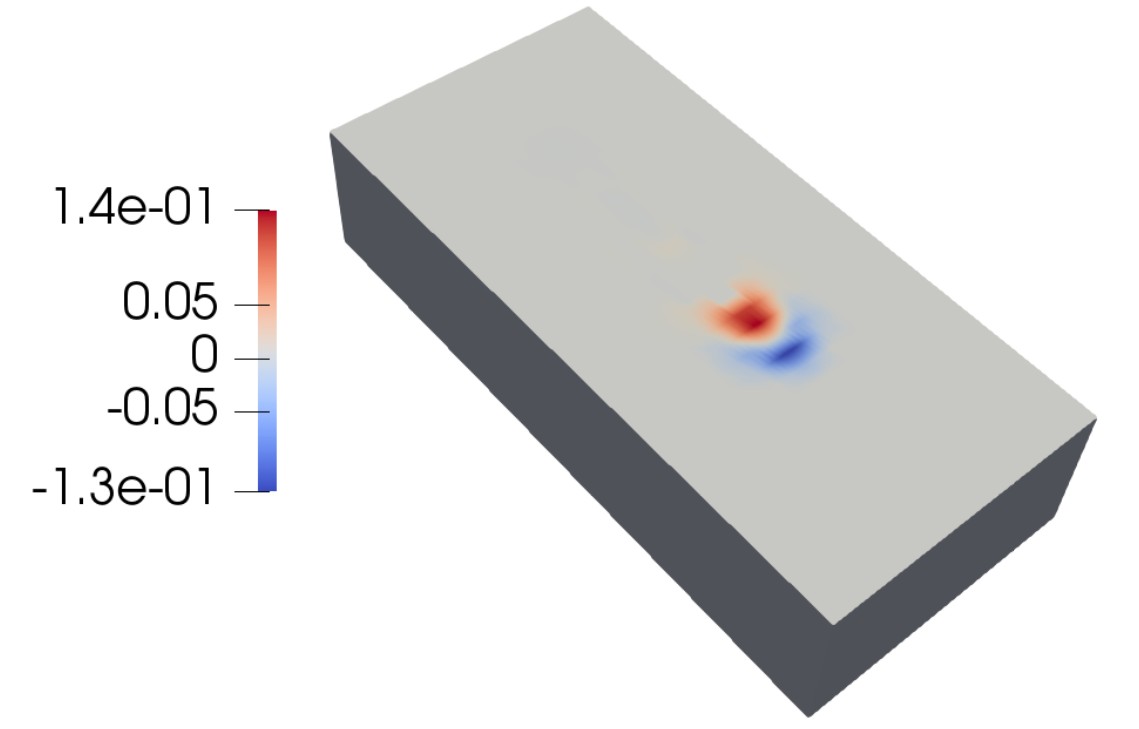}}
\subfigure[$A_1=1,\ A_2=1.5$ ]{\includegraphics[scale=0.2]{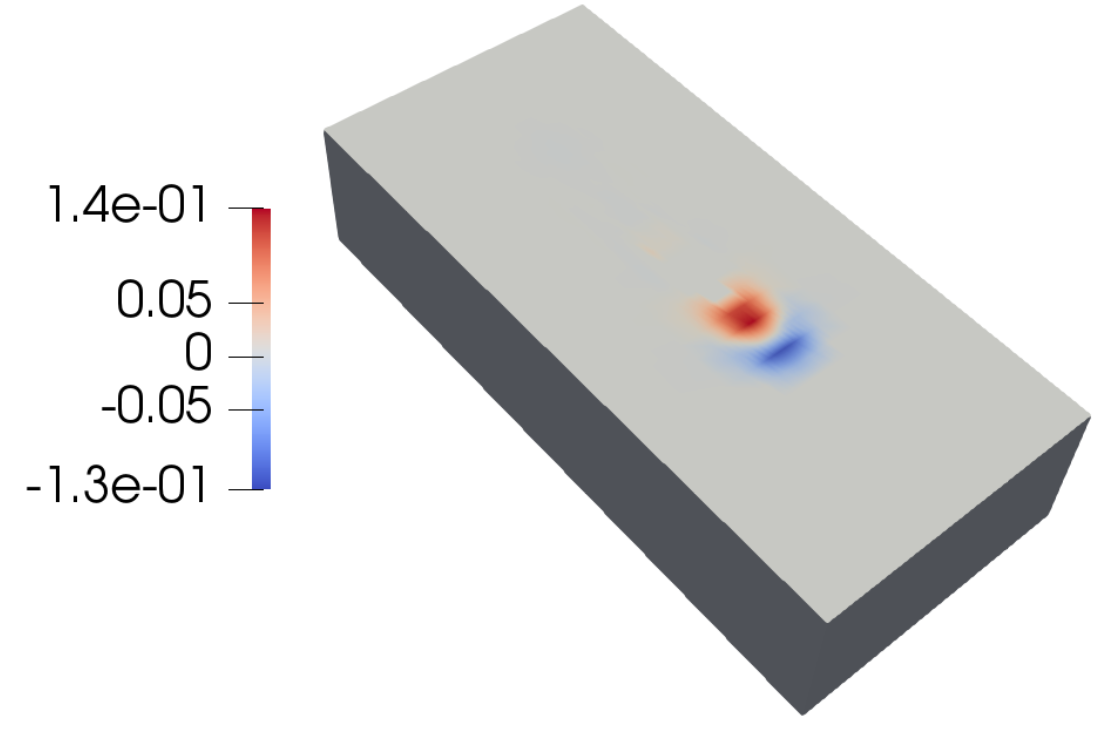}}
\subfigure[$A_1=2,\ A_2=1.5$ ]{\includegraphics[scale=0.2]{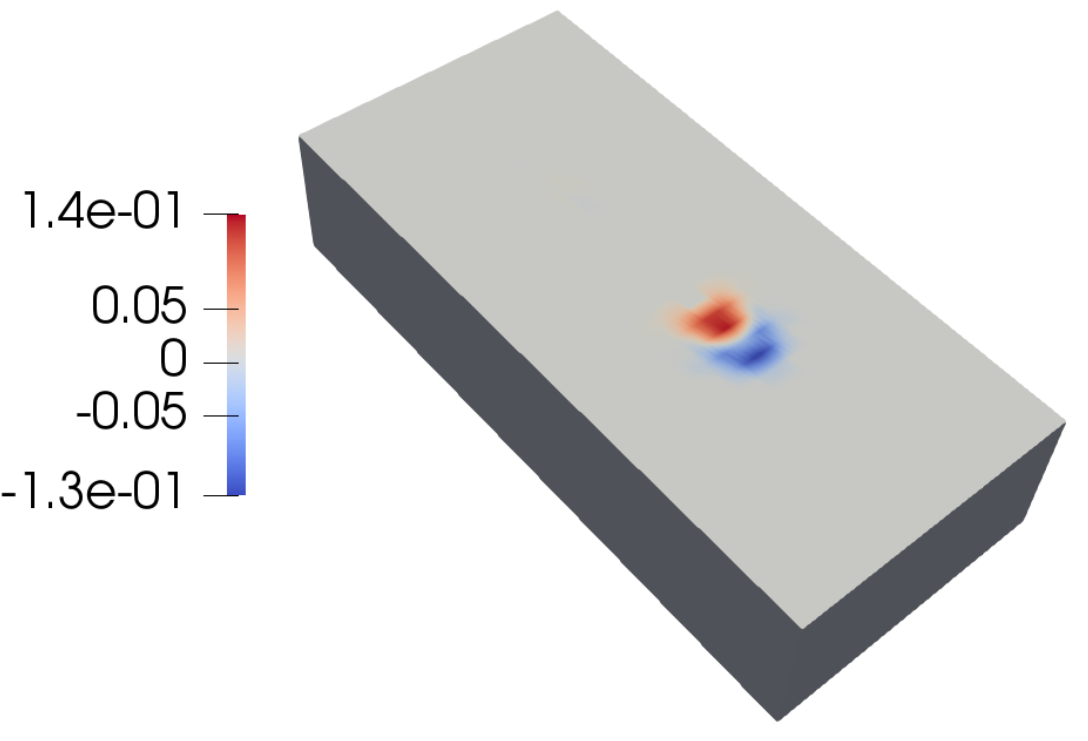}}
\subfigure[$A_1=1.6,\ A_2=2.5$ ]{\includegraphics[scale=0.2]{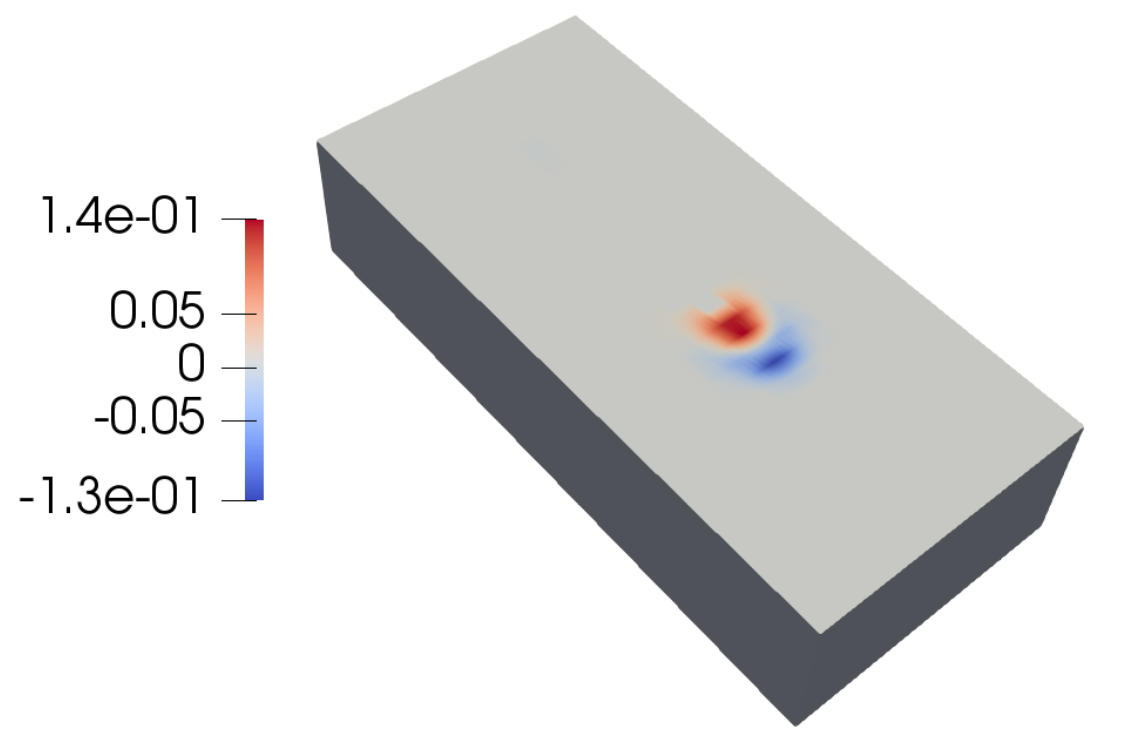}}\\
\subfigure[$A_1= 1.4,\ A_2= 1$ ]{\includegraphics[scale=0.2]{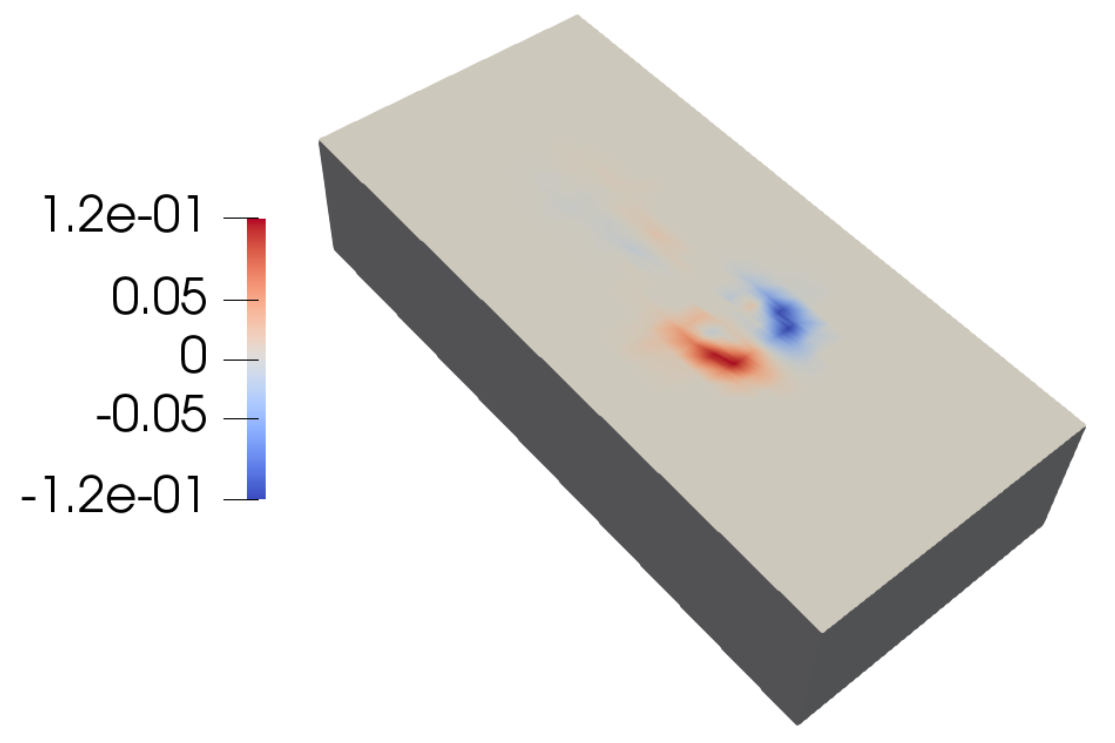}}
\subfigure[$A_1=1$, $A_2=1.5$ ]{\includegraphics[scale=0.2]{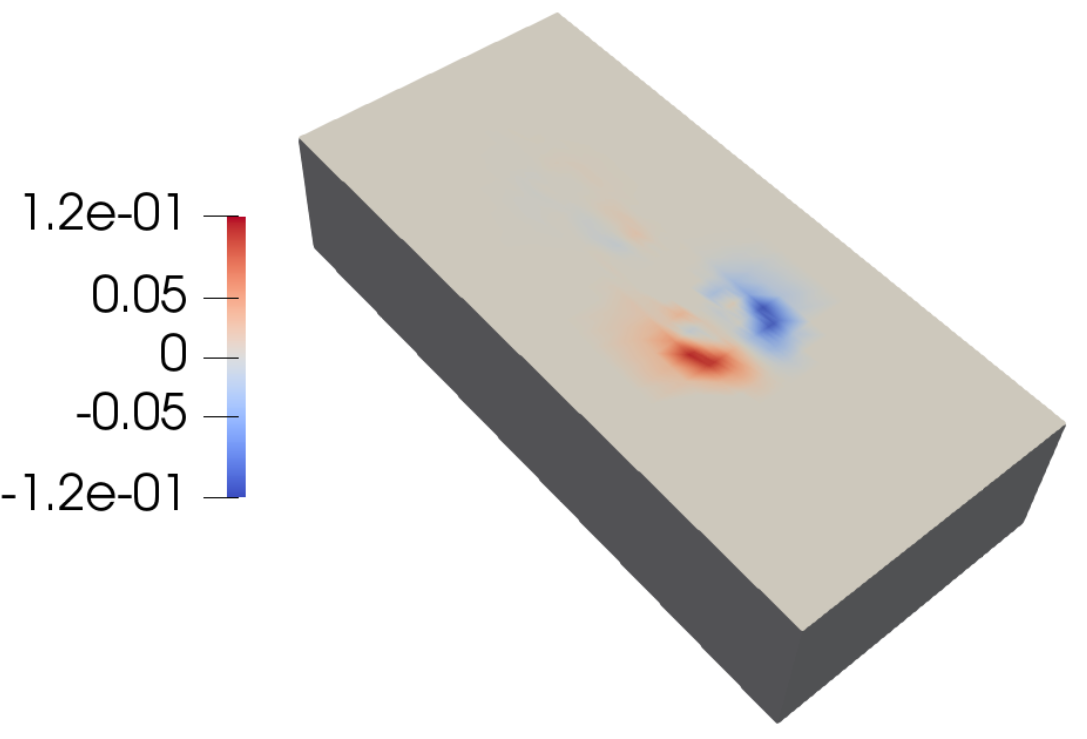}}
\subfigure[$A_1=2,\ A_2=1.5$ ]{\includegraphics[scale=0.2]{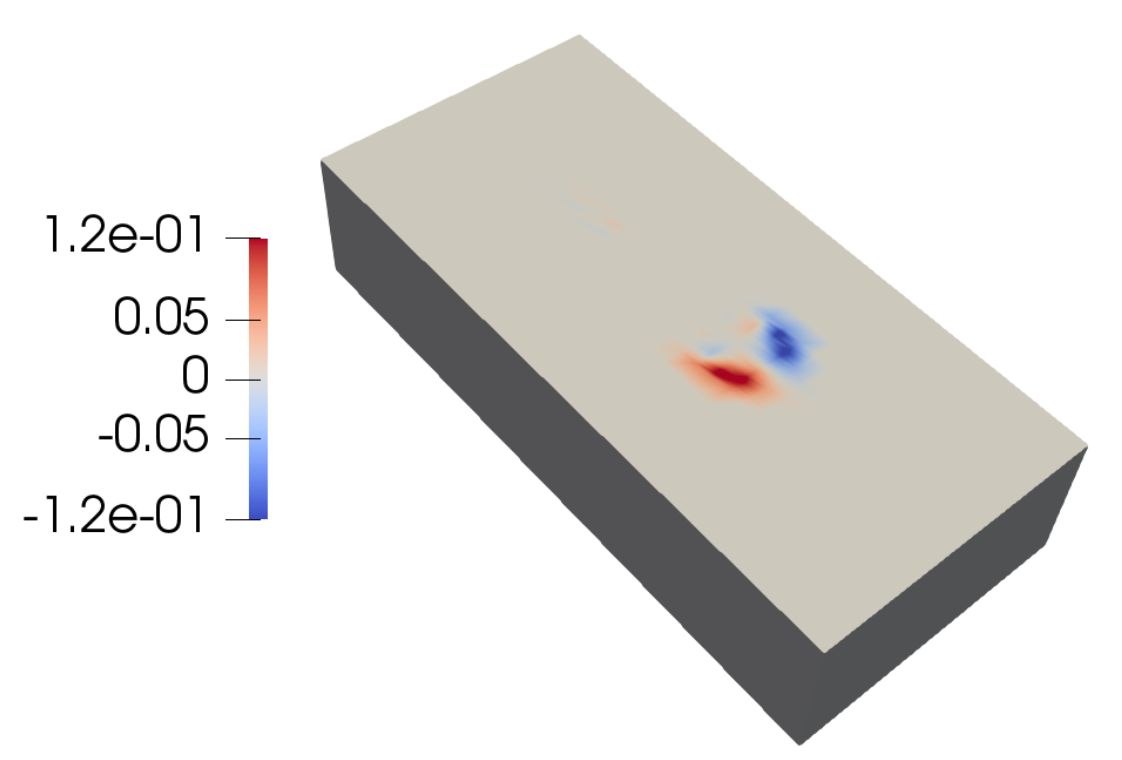}}
\subfigure[$A_1=1.6,\ A_2=2.5$ ]{\includegraphics[scale=0.2]{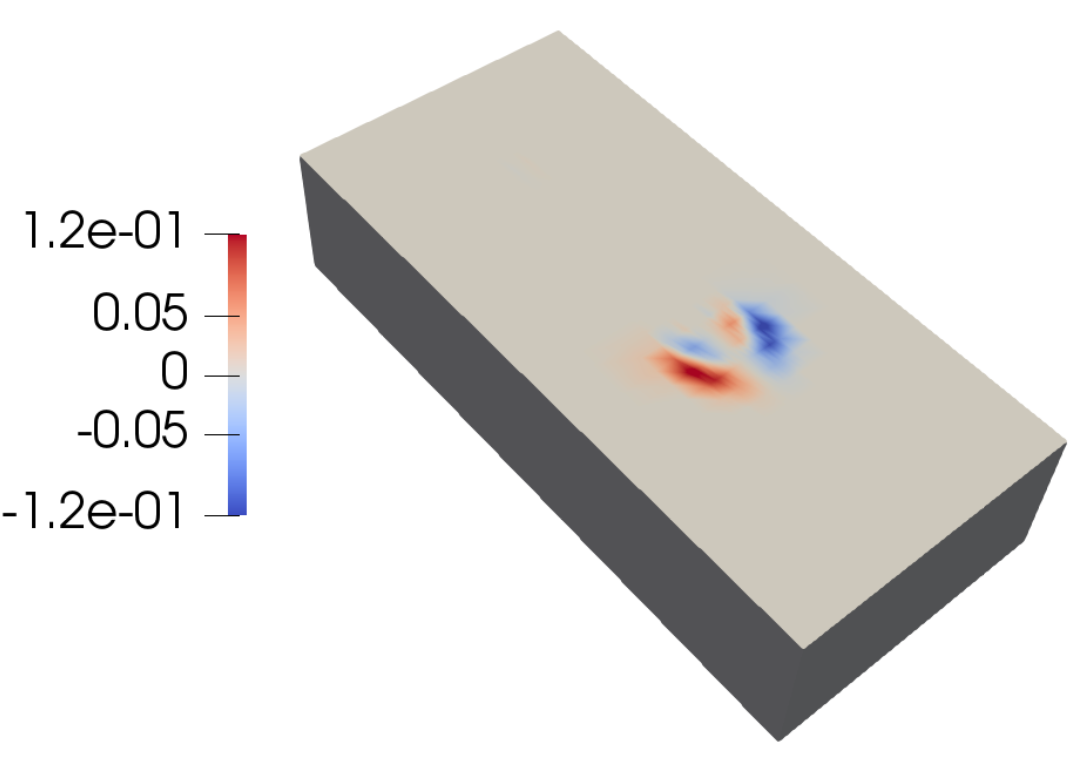}}\\
\subfigure[$A_1=1.4,\ A_2= 1$ ]{\includegraphics[scale=0.2]{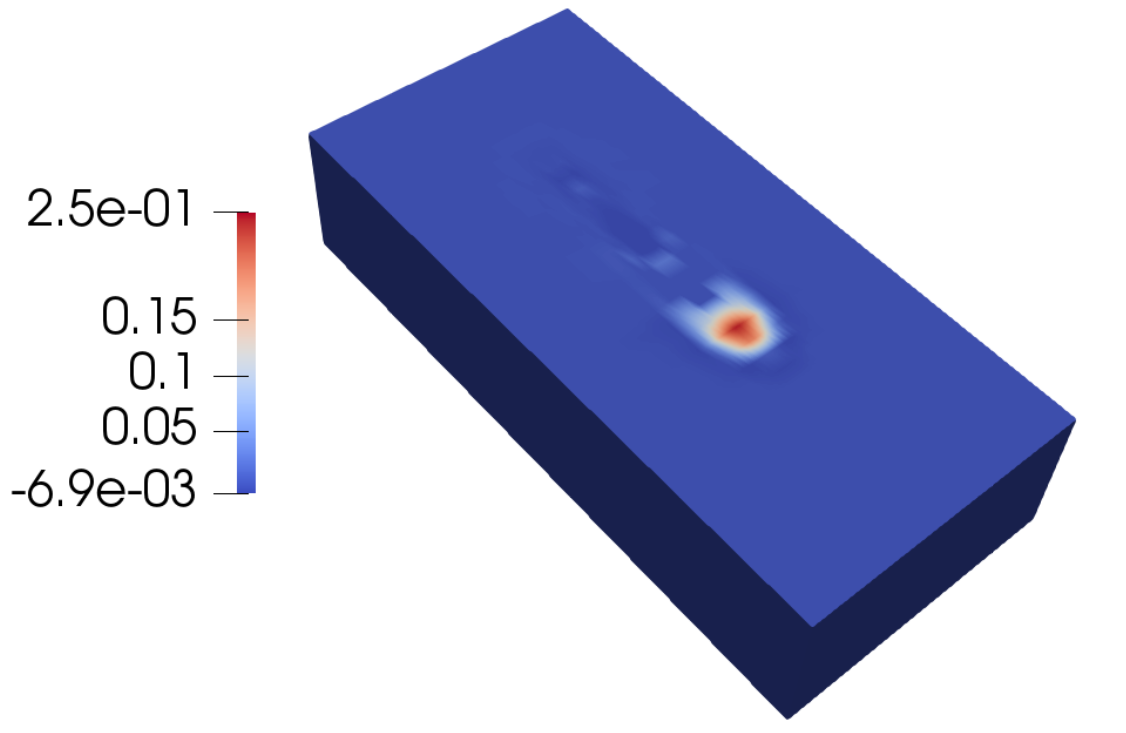}}
\subfigure[$A_1=1,\ A_2=1.5$ ]{\includegraphics[scale=0.2]{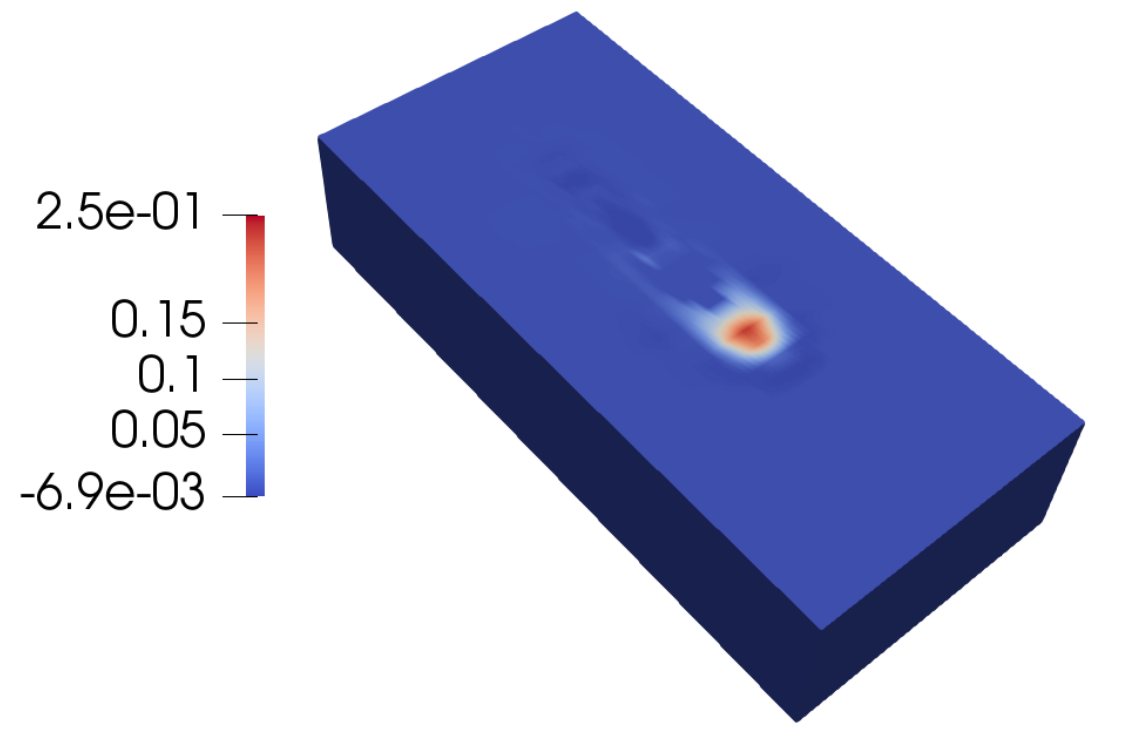}}
\subfigure[$A_1=2,\ A_2=1.5$ ]{\includegraphics[scale=0.2]{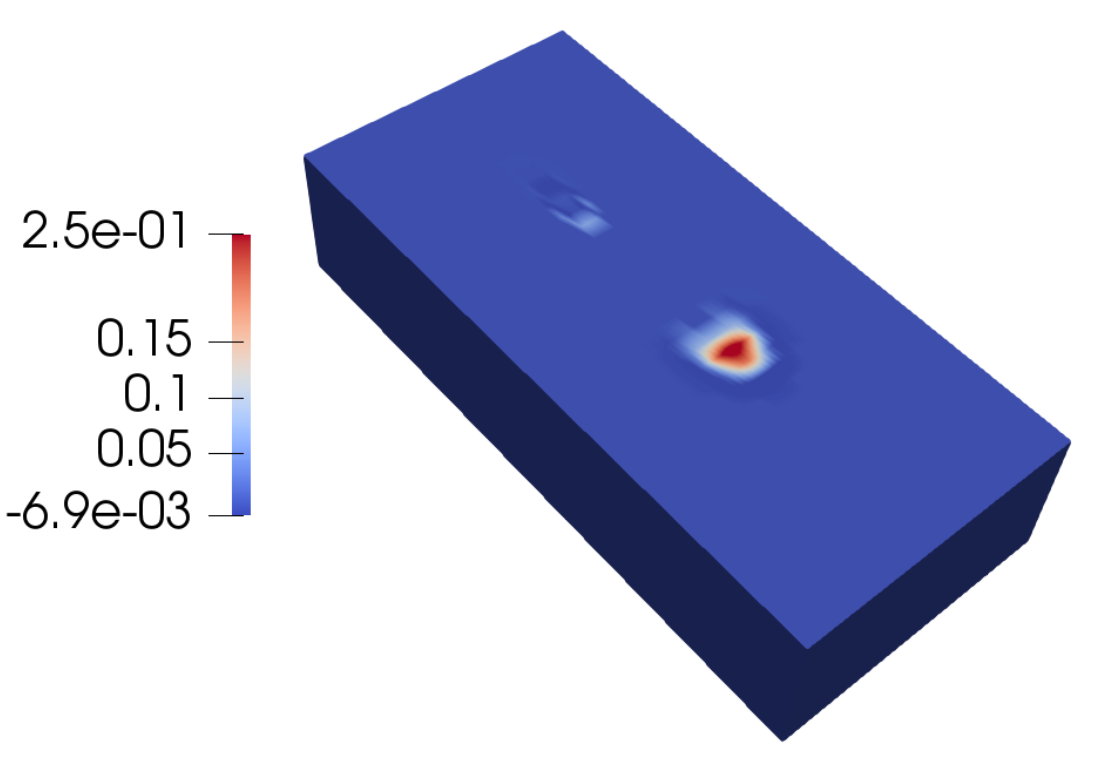}}
\subfigure[$A_1=1.6,\ A_2=2.5$ ]{\includegraphics[scale=0.2]{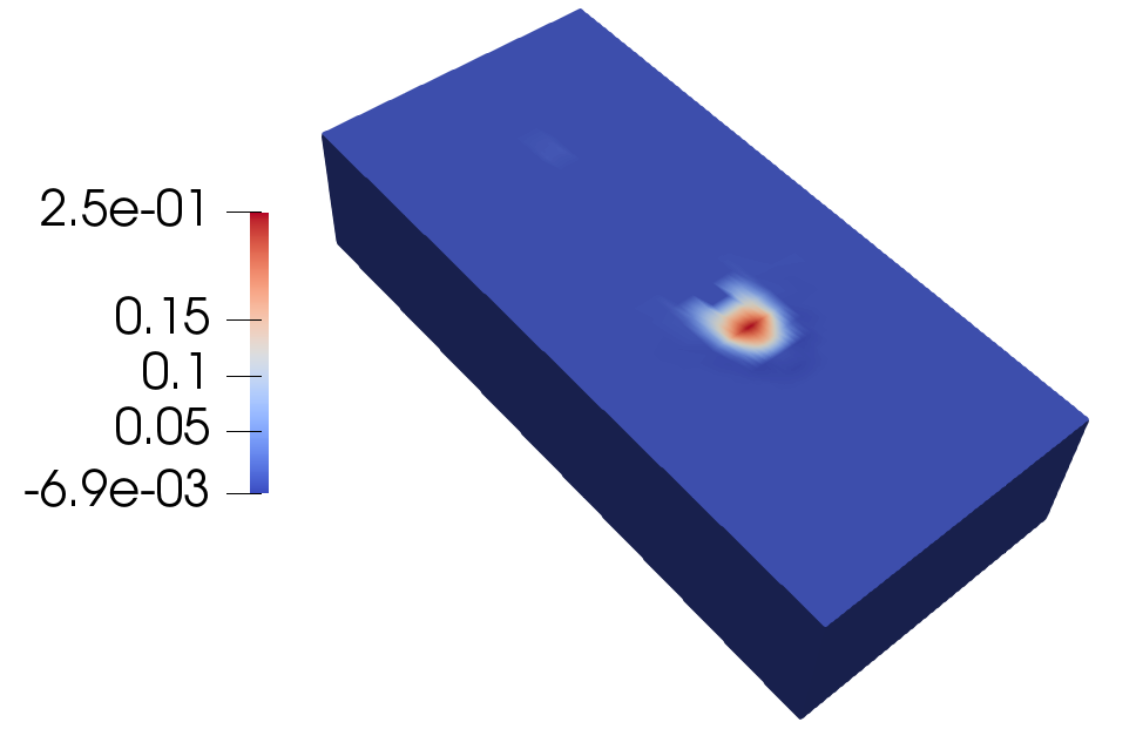}}\\
\caption{XTD  extended mode (normalized) $\Tilde{\boldsymbol{u}}^{(1)}$, time step $t_{30}$.
(a)-(d) $x$-component $\Tilde{\boldsymbol{u}}_1^{(1)}/\|\Tilde{\boldsymbol{u}}_1^{(1)}\|$,
(e)-(h) $y$-component $\Tilde{\boldsymbol{u}}_2^{(1)}/\|\Tilde{\boldsymbol{u}}_2^{(1)}\|$,
(i)-(l) $z$-component $\Tilde{\boldsymbol{u}}_3^{(1)}/\|\Tilde{\boldsymbol{u}}_3^{(1)}\|$.}
\label{fig:Ex2XTDmodeext}
\end{figure}

\begin{figure}[htbp]
\centering
{\includegraphics[scale=0.45]{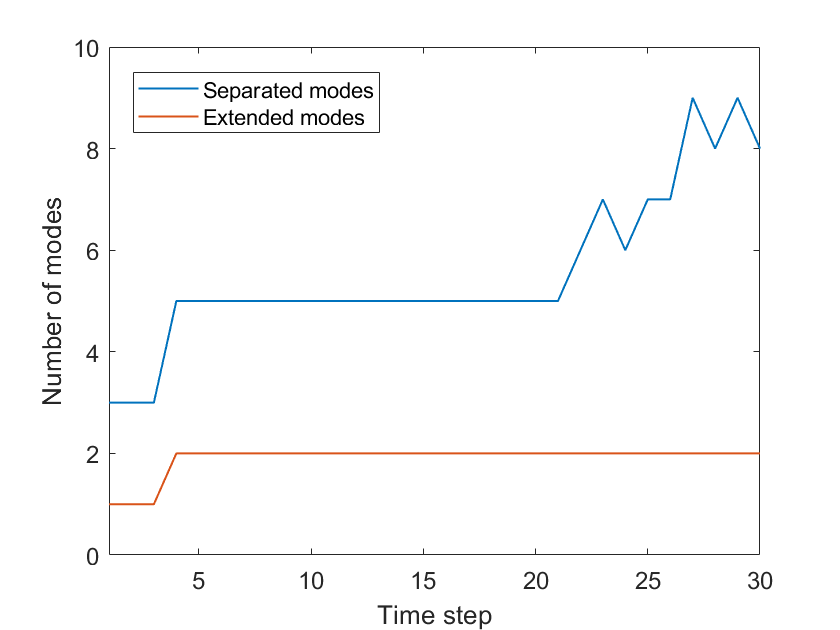}}
\caption{The number of modes for each time step in the second example}
\label{fig:ex2nbmode}
\end{figure}

For visualizing the space-time solution, the evolution of stress and plastic strain fields is depicted in \figurename~\ref{fig:Ex2XTD}. We can see that the overall solution field is evolving significantly over the time range, which challenges the commonly used model reduction methods in the literature. To some extend, the success of the present method in this challenging problem relies on the good choice of the separated variable, i.e. displacement increment. As shown in \figurename~\ref{fig:Ex2XTDmode} and \figurename~\ref{fig:Ex2XTDmodeext}, both the separated and extended modes are concentrated around the region of the heat source. This makes the problem reducible and enables the efficiency of the proposed solution method. The total number of modes at each time step for the XTD model is shown in \figurename~\ref{fig:ex2nbmode}. The overall number is limited within a range of 10. It increases only slightly with time, especially for the extended modes. 

This is the offline training stage for the XTD model. The online stage concerns more about the application (or usage) of these parametric solutions. In the following subsections, we present two applications concerning the online stage.

\subsubsection{Online fast uncertainty quantification}
The first online application concerns the uncertainty quantification (UQ) using the computed XTD solutions. One of the computationally costly problems in UQ is the uncertainty propagation. Given the material parameter uncertainty, the UQ problem consists in propagating the uncertainty to the model output, like residual stresses or distortions. The Monte Carlo (MC) simulation is usually adopted to this end, which is however very expensive as the numerical model (FEA) has to be run many times. In our work, since the parametric solution has been known by offline XTD training, the MC procedure becomes very fast and even in real time.

\begin{figure}[htbp]
\centering 
\def\svgscale{0.3}
{%% Creator: Inkscape inkscape 0.92.4, www.inkscape.org
%% Accompanies image file 'Ex2lineAB.pdf' (pdf, eps, ps)
\begingroup%
  \makeatletter%
  \providecommand\color[2][]{%
    \errmessage{(Inkscape) Color is used for the text in Inkscape, but the package 'color.sty' is not loaded}%
    \renewcommand\color[2][]{}%
  }%
  \providecommand\transparent[1]{%
    \errmessage{(Inkscape) Transparency is used (non-zero) for the text in Inkscape, but the package 'transparent.sty' is not loaded}%
    \renewcommand\transparent[1]{}%
  }%
  \providecommand\rotatebox[2]{#2}%
  \newcommand*\fsize{\dimexpr\f@size pt\relax}%
  \newcommand*\lineheight[1]{\fontsize{\fsize}{#1\fsize}\selectfont}%
  \ifx\svgwidth\undefined%
    \setlength{\unitlength}{557.00000192bp}%
    \ifx\svgscale\undefined%
      \relax%
    \else%
      \setlength{\unitlength}{\unitlength * \real{\svgscale}}%
    \fi%
  \else%
    \setlength{\unitlength}{\svgwidth}%
  \fi%
  \global\let\svgwidth\undefined%
  \global\let\svgscale\undefined%
  \makeatother%
  \begin{picture}(1,0.65798925)%
    \lineheight{1}%
    \setlength\tabcolsep{0pt}%
    \put(0,0){\includegraphics[width=\unitlength,page=1]{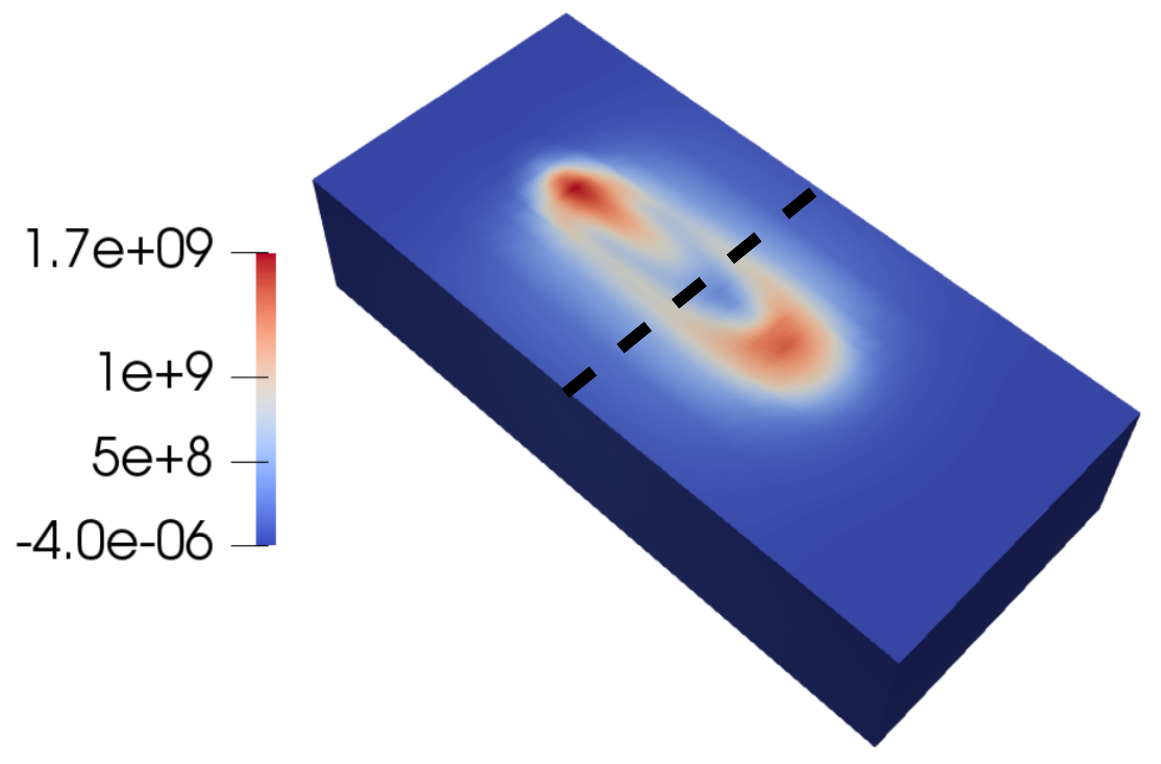}}%
    \put(0.40382146,0.22127471){\color[rgb]{0,0,0}\makebox(0,0)[lt]{\lineheight{1.25}\smash{\begin{tabular}[t]{l}A\end{tabular}}}}%
    \put(0.73165492,0.50946652){\color[rgb]{0,0,0}\makebox(0,0)[lt]{\lineheight{1.25}\smash{\begin{tabular}[t]{l}B\end{tabular}}}}%
  \end{picture}%
\endgroup%
}
\caption{Central line A-B for the residual stress observation}
\label{fig:Ex2lineAB}
\end{figure}

\begin{figure}[htbp]
\centering
\subfigure[$\sigma_{xx}$ over the central line A-B ]{\includegraphics[scale=0.25]{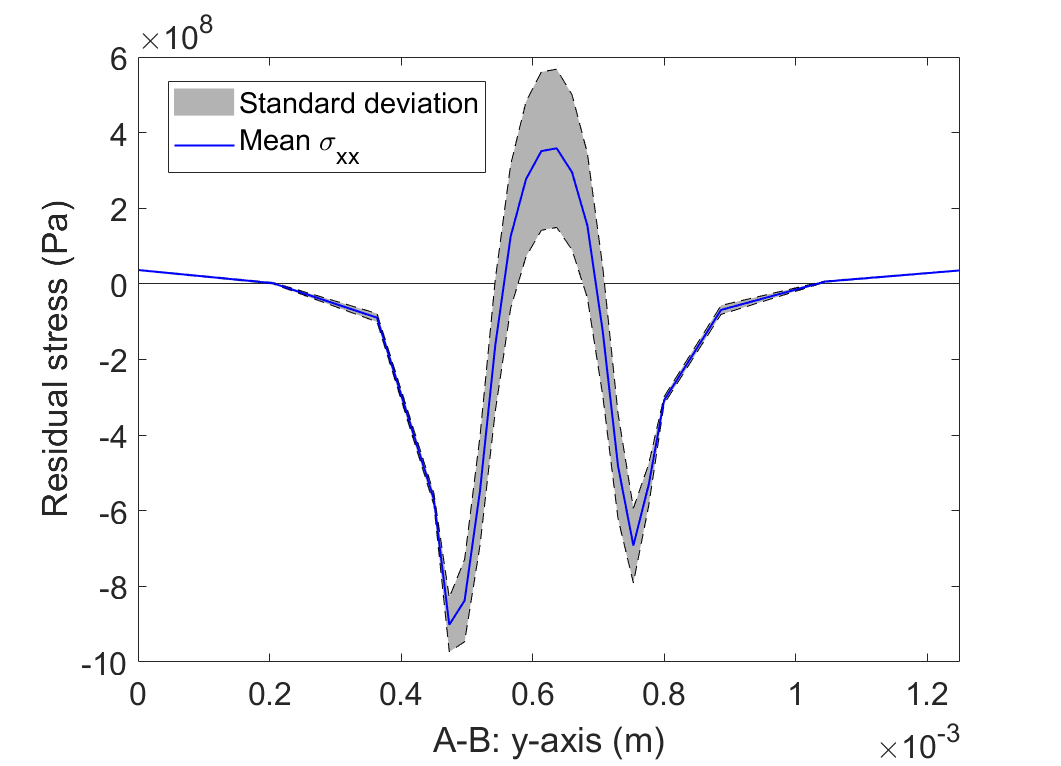}}\quad\quad
\subfigure[$\sigma_{yy}$ over the central line A-B ]{\includegraphics[scale=0.25]{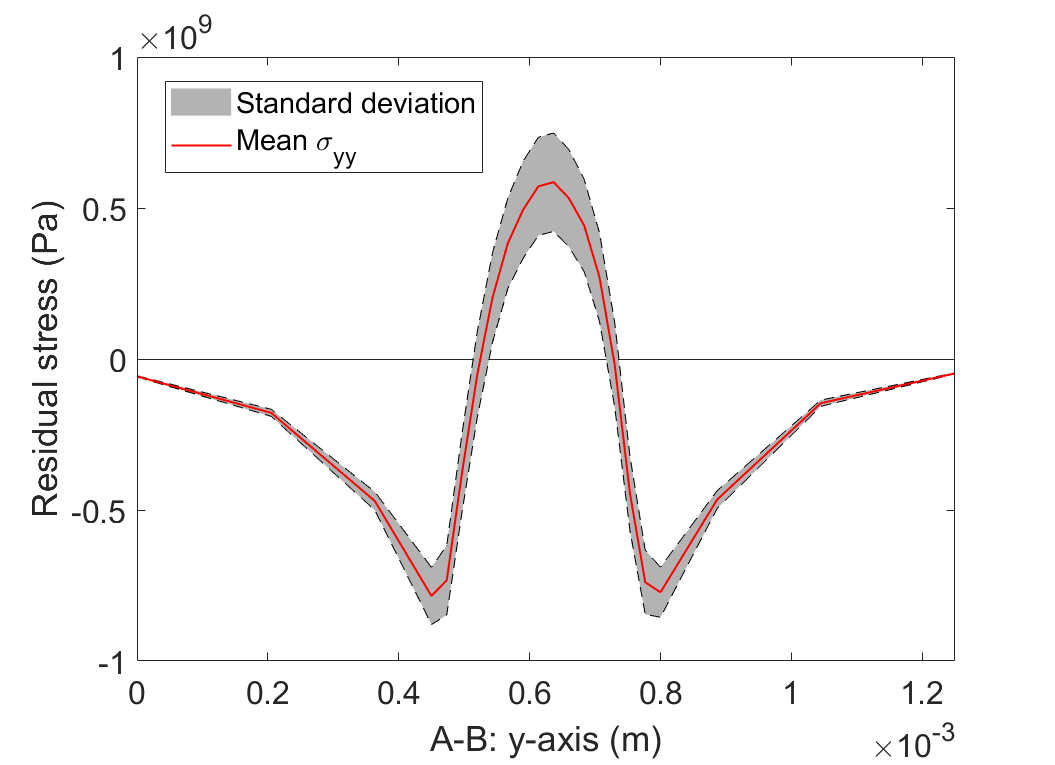}}
\caption{Online uncertainty quantification for residual stresses, cost $<0.2$ s}
\label{fig:Ex2onlinestd}
\end{figure}

Let us assume that the material properties $\sigma_y$ and $H$ are two independent random variables and follow two normal distributions, i.e.
\begin{equation}
\displaystyle
\label{eq:distribution}
\sigma_y \sim \mathcal{N}(1.5\sigma^{ref}_y,0.5\sigma^{ref}_y), \quad H \sim \mathcal{N}(2.1 H^{ref},0.5H^{ref})
\end{equation}
where $\sigma^{ref}_y$ and $H^{ref}$ are the reference parameter values shown in \tablename~\ref{table:ex2mate}. We can now study the variability of the residual stress with respect to the material randomness. As an example, the central line A-B (\figurename~\ref{fig:Ex2lineAB}) is chosen to observe the stress variation. Then, 63 random sets of material properties $\sigma_y$ and $H$ are generated according to the distribution \eqref{eq:distribution}. Using the XTD solution, we can easily obtain the propagated uncertainty on the residual stress, as shown in \figurename~\ref{fig:Ex2onlinestd}. Since there is no expensive computations at the online stage, the MC procedure only takes less than 0.2 s. If the parallel computing is used, this can be in real time.  The estimated cost for the MC simulation using FEA  for 65 runs is 69300 s, which seems unaffordable for real applications.

\subsubsection{Online fast material design/characterization under uncertainty}
In practice, the stochastic models of material properties are a priori unknown and need to be identified against experiments. Thus, we can consider an inverse  problem as below

\begin{equation}
\displaystyle
\label{eq:inverse}
\begin{cases}
\begin{aligned}
&\text{Find the hyper-parameters of the stochastic models of}\ \sigma_y,H\\
&\text{Minimize}\ J(\sig^{\text{XTD}}(\X,\sigma_y,H),\sig^{\text{target}})
\end{aligned}
\end{cases}
\end{equation}
where the hyper-parameters are the mean and standard deviation of the probabilistic distribution for $\sigma_y$ and $H$. $\sig^{\text{target}}$ can be the experimental measurement of the stress. This inverse problem consists in finding out the best stochastic models of materials that fit the experimental observation in terms of both means and variations. In this case, the objective function reads

\begin{equation}
\displaystyle
\label{eq:inverseJ}
J= w_1\|\text{mean}\ \sig^{\text{XTD}}-\text{mean}\ \sig^{\text{target}}\|+w_2\|\text{std}\ \sig^{\text{XTD}}-\text{std}\ \sig^{\text{target}}\|
\end{equation}
where std stands for the standard deviation. $w_1, w_2$ are the weights. For evaluating the objective function $J$, the statistics (mean and std)  of simulation outputs are needed. Hence, for a given set of hyper-parameters, a MC simulation is required. This makes the optimization extremely expensive as the hyper-parameters need to be modified over and over again. Fortunately, we can use XTD model to do the optimization. It takes only 5 minutes for this optimization problem at the online stage. This is reasonable as the online MC simulation is very cheap (less than 0.2 s), as shown previously.

\begin{figure}[htbp]
\centering
\subfigure[$\sigma_{xx}$ over the central line A-B ]{\includegraphics[scale=0.35]{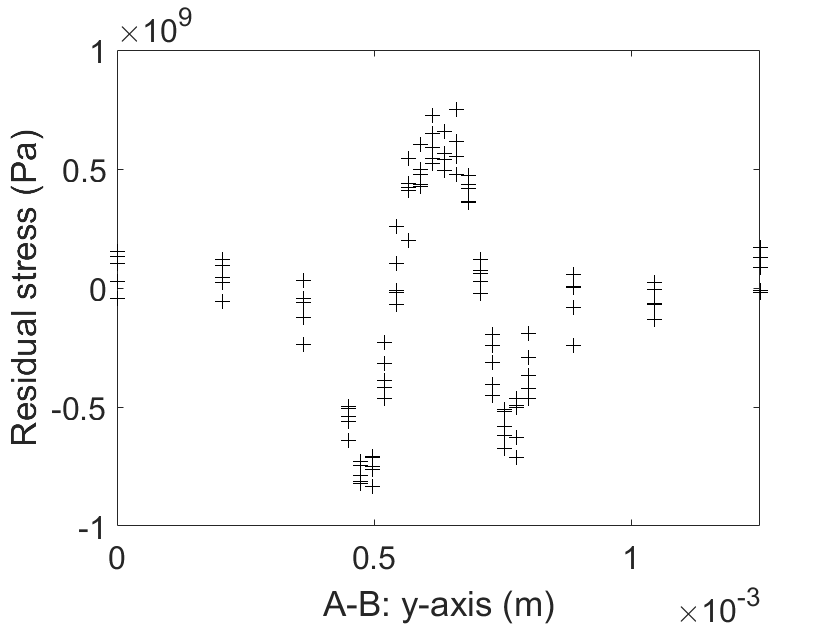}}\quad\quad
\subfigure[$\sigma_{yy}$ over the central line A-B ]{\includegraphics[scale=0.35]{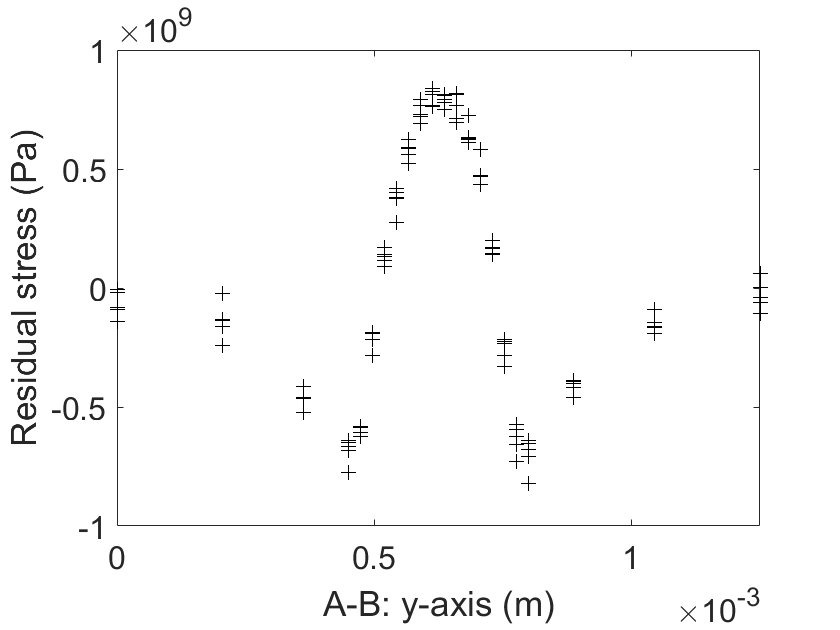}}
\caption{Experimental (target) residual stress }
\label{fig:Ex2exp}
\end{figure}

\begin{figure}[htbp]
\centering
\subfigure[$\sigma_{xx}$ over the central line A-B ]{\includegraphics[scale=0.29]{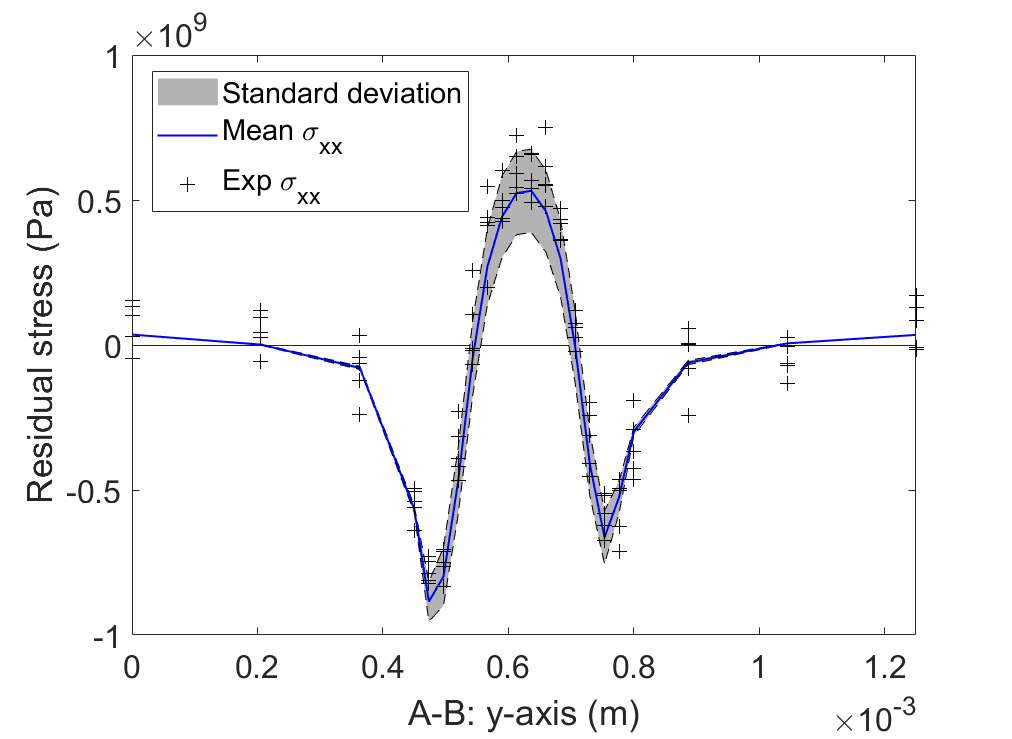}}\quad\quad
\subfigure[$\sigma_{yy}$ over the central line A-B ]{\includegraphics[scale=0.29]{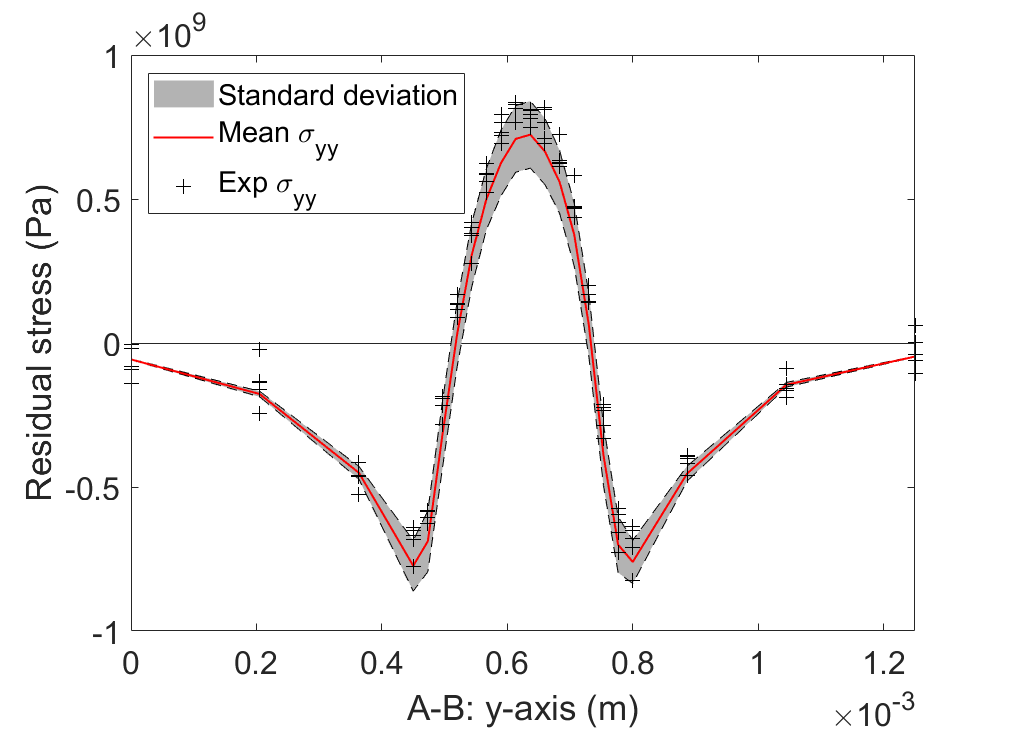}}
\caption{Identified model out against experimental (target) residual stress }
\label{fig:Ex2expstd}
\end{figure}

An example of the target experimental stress data is illustrated in \figurename~\ref{fig:Ex2exp}. We assume the two random material parameters follow two normal distributions independently. The  identified probabilistic models are given in \eqref{eq:distributionexp}. As shown in \figurename~\ref{fig:Ex2expstd}, the identified model shows a good agreement to the experimental measurement.

\begin{equation}
\displaystyle
\label{eq:distributionexp}
\sigma_y \sim \mathcal{N}(1.61\sigma^{ref}_y,0.38\sigma^{ref}_y), \quad H \sim \mathcal{N}(1.25 H^{ref},0.42H^{ref})
\end{equation}

\subsection{Coupled XTD-SCA strategy for multi-scale multi-parametric modeling of materials}
This example considers a two-scale coupled modeling for composite materials in the proposed XTD framework. In particular, the self-consistent clustering analysis (SCA) \cite{liu2016self} is used for the microscale simulation. Other options like FEA and the FFT (fast Fourier transform) method \cite{moulinec1998numerical} can also be used for the microscale problem. The SCA is adopted because of its efficiency. 

We summarize here the two-scale coupled solution procedure. In the macroscale, the XTD provides first an estimate of strain based on the nodal displacement increment. The macro-strain applies then as boundary conditions to the microscale representative volume element (RVE) associated to each macroscopic integration point. The SCA will provide a fast response of the micro-stress and micro-strain and return the homogenized stress for the macroscale problem. In other words, the SCA serves as the replacement to {the} constitutive material {laws} in the XTD framework. This allows to directly consider the anisotropic effect of the microstructure without making hypotheses on the constitutive laws. The cost to pay is the microscale simulations. Hence, the two-scale coupled solution is usually very expensive. The {presented} XTD-SCA method offers a way to enable real-time two scale simulations.

\begin{figure}[htbp]
\centering 
\def\svgscale{0.3}
{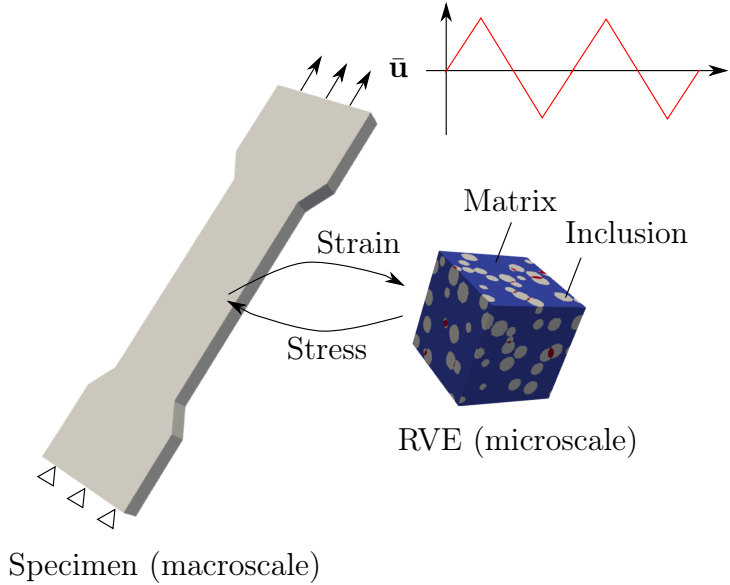}
\caption{Two-scale coupled composite modeling}
\label{fig:Ex3specimen}
\end{figure}

In this example, we consider {that} the RVE is the same for all the macroscopic integration points. The geometry of the specimen is shown in \figurename~\ref{fig:Ex3specimen}. The RVE has two material phases, matrix and inclusion. The inclusion is considered purely elastic, whereas the matrix materials are elastoplastic with isotropic hardening.  The yield stress of the RVE matrix is described by the following equation
\begin{equation}
\displaystyle
\label{eq:XTDex3yield}
\sigma_Y=\sigma_y+H \text{p}^n
\end{equation}
where $\sigma_y$ is the initial yield stress, $H$ and $n$ are the hardening coefficients, \text{p} is the equivalent plastic strain. 

The specimen is subject to a cycle load. The increment of displacement loading is given as follows: $\Delta\bar{u}_y=0.00002$ m for the first two steps, then $\Delta\bar{u}_y=0.00001$ m for the following four steps, $\Delta\bar{u}_y=-0.00002$ m for the remaining four steps in a half cycle. The second half cycle follows the same step increments and switches to a compression mode. This complete cycle is repeated twice in this example. There is in total 40 time steps.

\begin{figure}[htbp]
\centering
\subfigure[Clusters in the matrix phase ]{\includegraphics[scale=0.32]{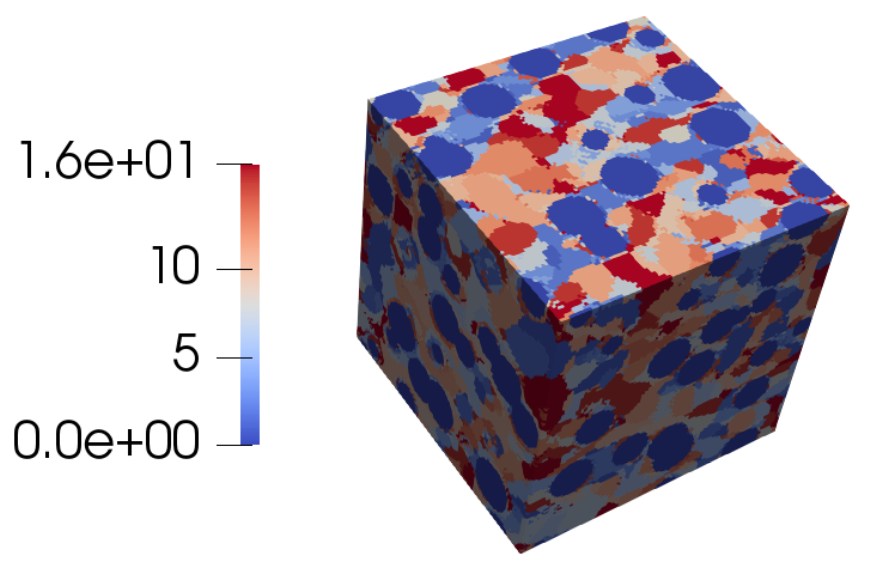}}\quad\quad\quad
\subfigure[Clusters in the inclusion phase ]{\includegraphics[scale=0.32]{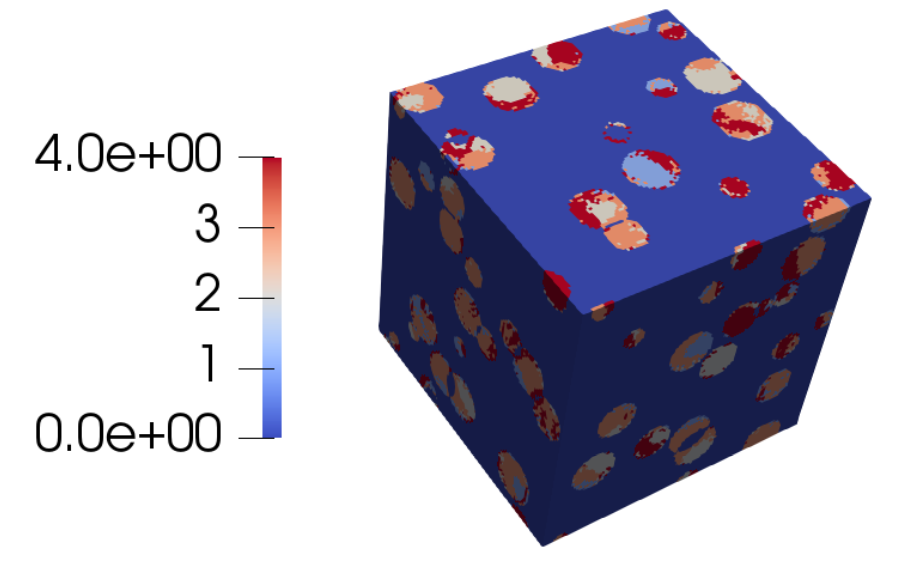}}
\caption{Clusters for the micro-structural RVE (legends are cluster ID)}
\label{fig:Ex3cluster}
\end{figure}

The specimen is discretized by a mesh of 160 cubic elements with 8 integration points. The RVE is discretized by a mesh of $100\times100\times100$. For SCA, the matrix domain is divided by 16 clusters, whereas the inclusion is divided into 4 clusters. The clustering result is depicted in \figurename~\ref{fig:Ex3cluster}.

For the XTD model, we consider two variable parameters, $\sigma_y\in [\sigma_y^{\min},\ \sigma_y^{\max}]$ and $H\in [H^{\min},\ H^{\max}]$. The following parametric solutions are computed for each time step.
\begin{equation}
\displaystyle
\label{eq:XTDex3}
\Delta\ud(\X,\sigma_y,H)= \sum_{m=1}^{M}\boldsymbol{a}^{(m)}(\X)g_1^{(m)}(\sigma_y)g_2^{(m)}(H) + \sum_{k=1}^{K}\Tilde{\boldsymbol{u}}^{(k)}(\X,{\sigma_y},H)
\end{equation}
The material properties for the RVE are illustrated in \tablename~\ref{table:ex3mate}. The initial yield stress $\sigma_y$ is discretized using 11 points. The hardening coefficient $H$ is discretized using 3 points. 

\begin{figure}[htbp]
\centering
\subfigure[$t_6$ ]{\includegraphics[scale=0.32]{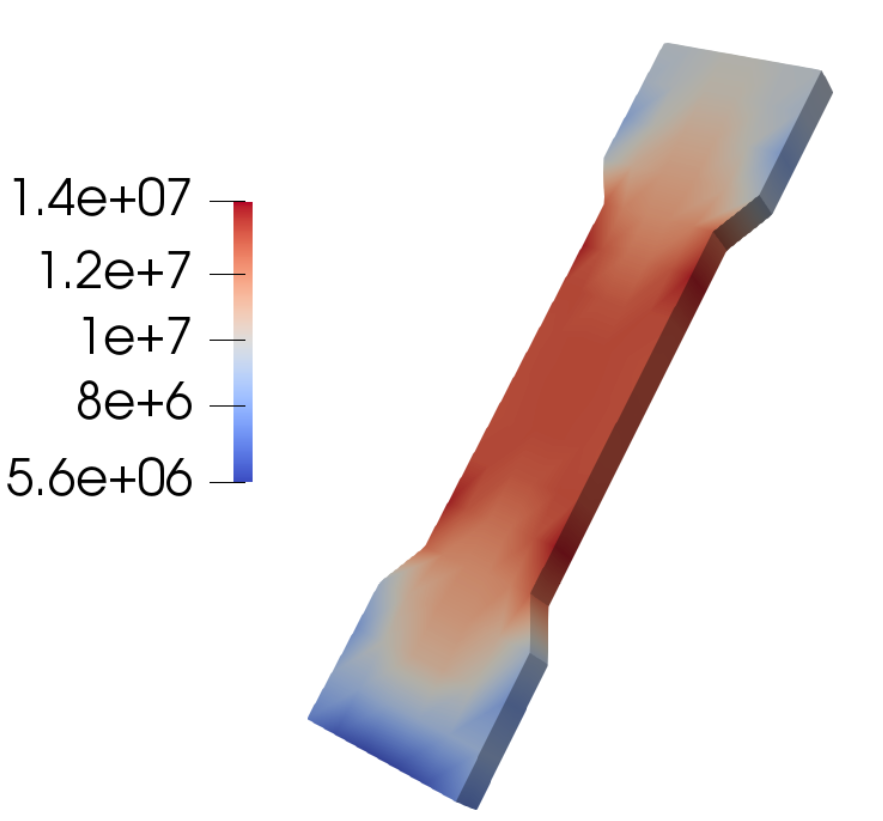}}
\subfigure[$t_{17}$ ]{\includegraphics[scale=0.32]{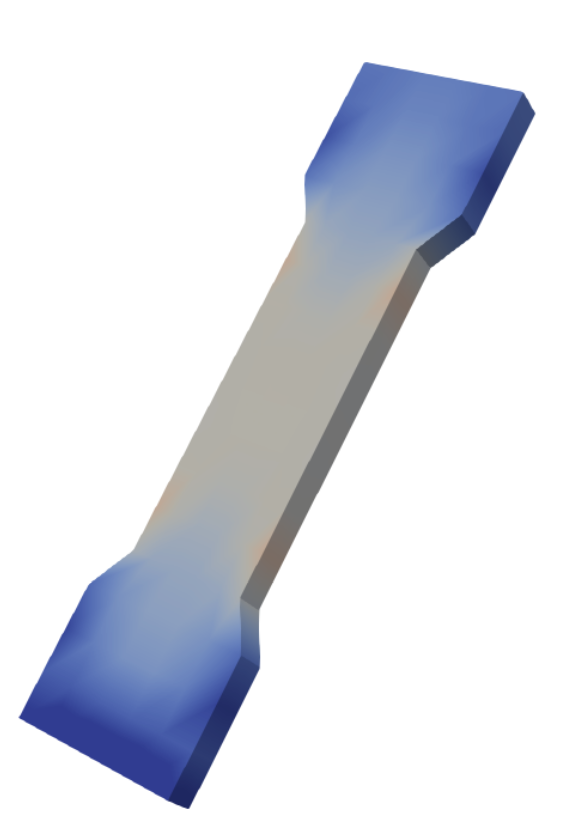}}
\subfigure[$t_{26}$ ]{\includegraphics[scale=0.32]{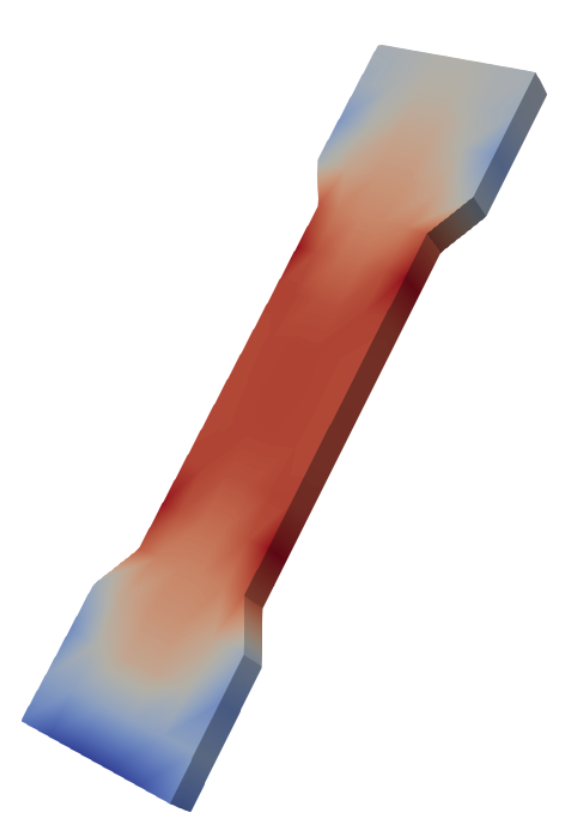}}
\subfigure[$t_{36}$ ]{\includegraphics[scale=0.32]{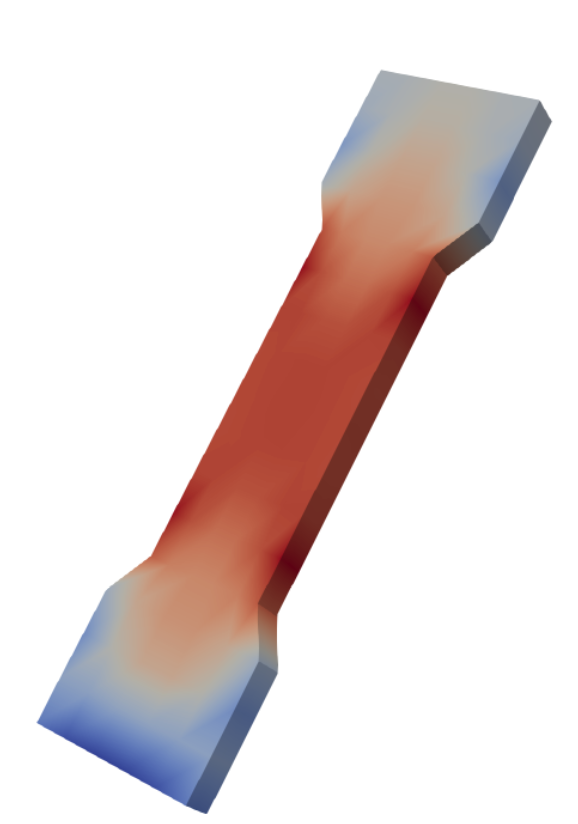}}
\subfigure[$t_{6}$
]{\includegraphics[scale=0.28]{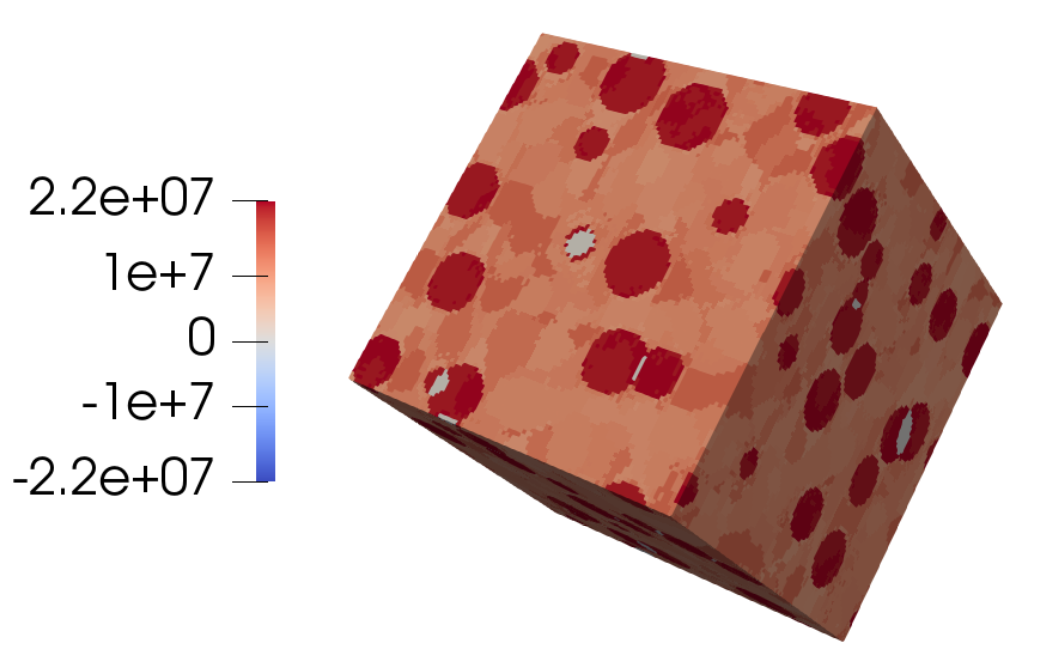}}
\subfigure[$t_{17}$ ]{\includegraphics[scale=0.28]{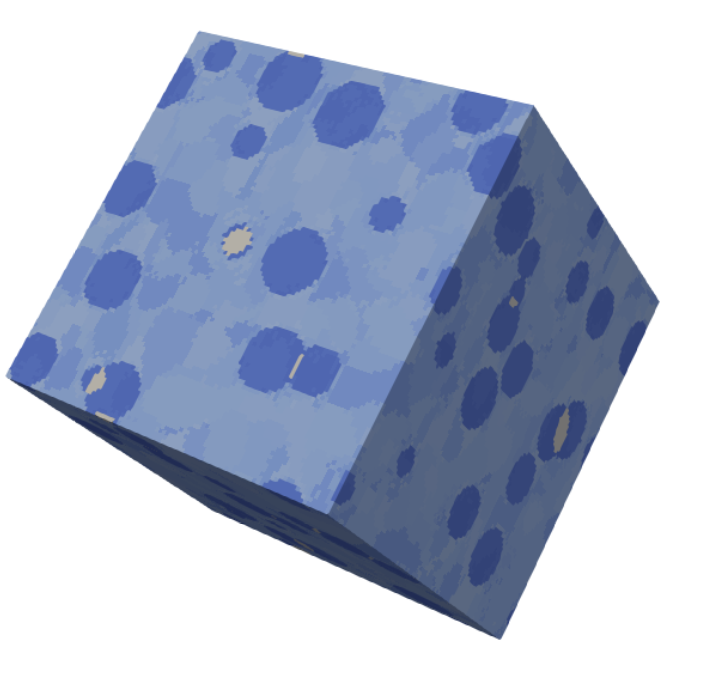}}
\subfigure[$t_{26}$ ]{\includegraphics[scale=0.28]{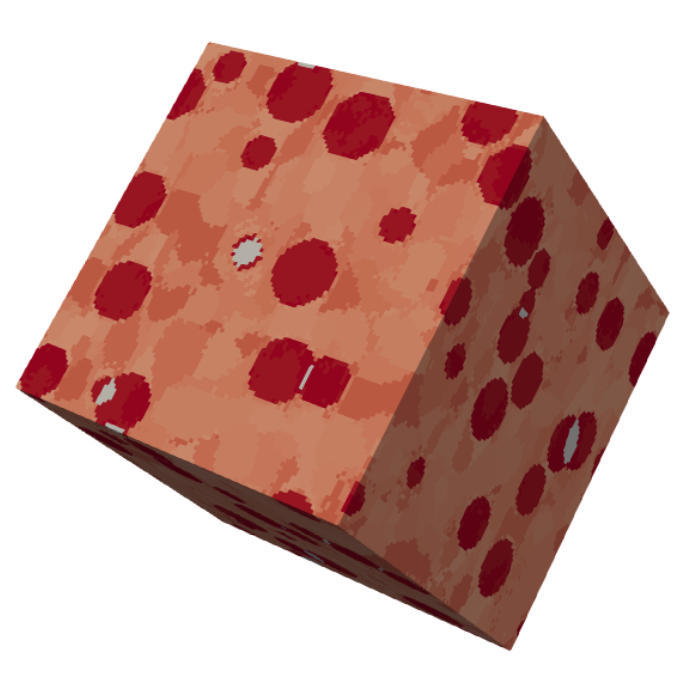}}
\subfigure[$t_{36}$ ]{\includegraphics[scale=0.28]{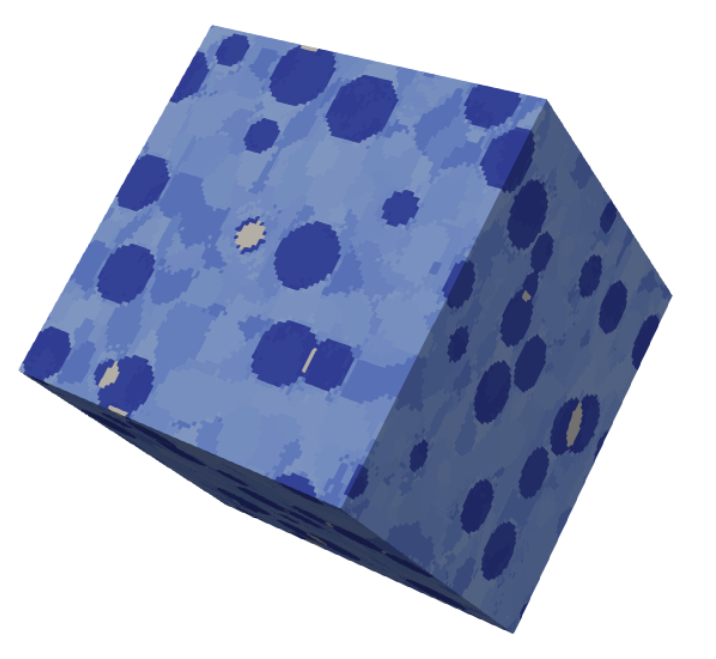}}
\caption{Two scale response of materials to the cyclic load. (a)-(d) macroscopic von Mises stress, (e)-(h) microscopic von Mises stress.}
\label{fig:Ex3XTD}
\end{figure}

The offline computation of the two-scale coupled XTD-SCA model takes 37633 s. For comparison, the FEA-SCA would take about 105044 s for the 33 runs. If a larger mesh size is used, the speedup is even more significant. Again, once the offline computation is finished, the online XTD prediction for a given parameter set is in real time ($<0.01$ s). \figurename~\ref{fig:Ex3XTD} illustrates the two-scale von Mises stress prediction for the parameters: $\sigma_y=0.013$ (GPa) and $H=0.045$ (GPa). As applications, {the} XTD-SCA enables real-time monitoring of the two-scale material response and can be used to monitor the effect of microscale defects or the initiation of micro-cracks of materials. It should be noted that this is intractable with conventional FEA as one prediction from FEA-SCA would take more than 3000 s. If the FE$^2$ approach is used, the computational cost will be higher than 30000 s.

\begin{table}[htbp]
\caption{Material properties of RVE}
\centering
\begin{tabular}{|c|c|c|c|c|c|c|c|}
\hline
Phase & $E$ & $\nu$  & $\sigma_y^{\min}$& $\sigma_y^{\max}$ & $H^{\min}$ &$H^{\max}$&$n$\\ \hline
Inclusion & 66 GPa &  0.19 & - & - &- &- &- \\ \hline
Matrix & 3 GPa & 0.4 & 0.01 GPa & 0.02 GPa &0.03 GPa &0.06  GPa& 0.4 \\ \hline
\end{tabular}\\
$E$: Young's modulus, $\nu$: Poisson's ratio
\label{table:ex3mate}\\
\end{table}

{
\subsection{Discussion on the application of XTD to other types of problems}
One key idea that enables the efficiency of the XTD method is the sparsity of the non-separated enrichment in \eqref{eq:XTDdefcontinu}. Choosing appropriate local enrichment regions becomes essential for having the desired performance in terms of efficiency and accuracy. The criterion used for determining the local (sparse) enrichment region is specially designed for the presented problems. To apply the method to other problems, a specific criterion needs to be well defined for ensuring the sparsity of enrichment. Some assumptions might be needed. For example, for computational fluid mechanics, the nonlinear convection terms can be assumed to be dominant only in some local regions, at least for an incremental time step. In this case, the criterion for the local enrichment region might be related to the velocity field. As the velocity is usually non-uniform, the sparsity of the enrichment might be achieved with minimum compromising of accuracy, if the criterion is well chosen.  Hence, the method should be applicable to various problems where potential singularity and nonlinear strong coupling phenomena can be assumed to be localized in part of the global space-time-parameter domain.}

{
Another key point that makes the method efficient is the existence of separable modes. For problems that involve hardly separable variables that are not localized, such as finite elastoplasticity problems, further investigations are needed on how to efficiently compute the separated modes and choose the enrichment domain. In machine learning context, the architectures of neural networks that connect different variables can be very complex to have the universal approximation capability. If separation of variables is used to simplify machine learning structures, the resulting approximation capability may be reduced, to some extent. However, the advantage is the reduced number of DoFs and therefore reduced training costs. A balance between approximation capability and efficiency can be made with adaptive separation of variables.
}

\section{Conclusion}
This paper {presented} a novel model reduction method, referred to as XTD model reduction. The basic idea is to introduce an extra sparse term in the tensor decomposition, such as the canonical tensor decomposition. This extra-term is expected to improve the approximation accuracy and reduce the number of modes when encountering highly nonlinear or even singular problems. The overall XTD framework is split into two stages: offline and online. The offline training consists in solving the problem by considering that the solution has the XTD format. The online prediction is real time and can be used for many engineering problems in which parametric solutions are inquired many times, like design and uncertainty quantification.

The proposed method has been applied to nonlinear elastic-plastic problems, including challenging problems like additive manufacturing residual stress prediction and multi-scale composite material modeling. The method has shown significant speedups to conventional FEA, especially when the mesh size is very large. Applications concerning uncertainty quantification and inverse design and calibration have been considered in our work, which show a great potential to enable a wide range of applications of the overall framework.

Other applications can be considered in terms of multi-scale system monitoring and control for additive manufacturing or other material forming processes. The XTD method can also be combined with other machine learning techniques to further improve efficiency.

\section*{Acknowledgement}
The authors would like to acknowledge the support of the National Science Foundation under Grant No. CMMI-1762035 and CMMI-1934367.

\appendix
{\section{XTD based data learning and regression}\label{XTDdata}}
Without loss of generality, let us consider a scalar-valued multidimensional function $u(\mu_1,\dots,\mu_n)$ and the data of it is known for the general parameters: $(\mu_1,\dots,\mu_n) \in \Omega:=\Omega_{\mu_1}\times\cdots\times\Omega_{\mu_n}$. The training problem \eqref{eq:XTDdata} becomes
\begin{equation}
\displaystyle
\label{eq:XTDdatascalar}
\int_{\Omega}\delta u^{\text{XTD}}(u-u^{\text{XTD}})\  d\mu_1\cdots d\mu_{n}=0
\end{equation}
where the XTD representation of $u$ reads
\begin{equation}
\displaystyle
\label{eq:XTDdefcontinuscalar}
u^{\text{XTD}}({\mu_1},\dots,\mu_n)= {\sum_{m=1}^{M}g_1^{(m)}(\mu_1)\cdots g_n^{(m)}(\mu_n) + \sum_{k=1}^{K}\Tilde{{u}}^{(k)}( {\mu_1},\dots,\mu_n)}
\end{equation}
or in a discrete form
\begin{equation}
\displaystyle
\label{eq:XTDdefdiscretescalar} \boldsymbol{\mathcal{U}}^{\text{XTD}}=\sum_{m=1}^{M}
\boldsymbol{g}_1^{(m)}\otimes\cdots\otimes\boldsymbol{g}_n^{(m)}+\sum_{k=1}^{K}\Tilde{\boldsymbol{\mathcal{U}}}^{(k)}
\end{equation}
with the $g_i^{(m)}$  and $\Tilde{{u}}^{(k)}$  denoting respectively the separated modes and the extended (non-separated) modes. Similar to PGD and HOPGD, the XTD uses an incremental scheme to compute the modes.

\subsection{Compute separated modes with $\Tilde{{u}}^{(1)}=0$}
Starting from $M=1,K=1$, then
\begin{equation}
\displaystyle
\label{eq:XTDonemodeu}
u^{\text{XTD}}=g_1^{(1)}\cdots g_n^{(1)} + \Tilde{{u}}^{(1)}
\end{equation}
\begin{equation}
\displaystyle
\label{eq:XTDdeltau}
\delta u^{\text{XTD}}=\delta g_1^{(1)}\cdots g_n^{(1)}+\cdots+g_1^{(1)}\cdots\delta g_n^{(1)} + \delta \Tilde{{u}}^{(1)}
\end{equation}

We can now compute an estimate of $g^{(1)}_1$ by giving initial values to $g^{(1)}_i,\ \forall i\neq 1$ and with an initial guess $\Tilde{{u}}^{(1)}=0$. This also leads to $\delta g^{(1)}_i=0,\ \forall i\neq 1$ and $\delta \Tilde{{u}}^{(1)}=0$. Hence, 
\begin{equation}
\displaystyle
\label{eq:uxtdonemode}
 u^{\text{XTD}}= g_1^{(1)}\cdots g_n^{(1)}
\end{equation}
\begin{equation}
\displaystyle
\label{eq:deltaxtdonemode}
\delta u^{\text{XTD}}=\delta g_1^{(1)}\cdots g_n^{(1)}
\end{equation}
Substituting \eqref{eq:uxtdonemode} and \eqref{eq:deltaxtdonemode} into \eqref{eq:XTDdatascalar}, we obtain
\begin{equation}
\displaystyle
\label{eq:XTDdataonemode}
\int_{\Omega}\delta g_1^{(1)}\cdots g_n^{(1)} (u-g_1^{(1)}\cdots g_n^{(1)})\  d\mu_1\cdots d\mu_{n}=0
\end{equation}
By rearrangement,
\begin{equation}
\displaystyle
\label{eq:XTDdataonemode2}
\int_{\Omega}\delta g_1^{(1)}\cdots g_n^{(1)} u\  d\mu_1\cdots d\mu_{n}=\int_{\Omega_{\mu_1}}\delta g_1^{(1)}g_1^{(1)}\ d\mu_1 \cdots \int_{\Omega_{\mu_n}}g_n^{(1)} g_n^{(1)} \ d\mu_{n}
\end{equation}
The discretized form for the right hand side of \eqref{eq:XTDdataonemode2} can be written as
\begin{equation}
\displaystyle
\label{eq:XTDdatarhs}
\int_{\Omega_{\mu_1}}\delta g_1^{(1)}g_1^{(1)}\ d\mu_1 \cdots \int_{\Omega_{\mu_n}}g_n^{(1)} g_n^{(1)} \ d\mu_{n}=(\delta\boldsymbol{g}_1^{(1)})^T\boldsymbol{g}_1^{(1)}\cdots(\boldsymbol{g}_n^{(1)})^T\boldsymbol{g}_n^{(1)}
\end{equation}
The left hand side of \eqref{eq:XTDdataonemode2} requires a full integration but can be performed in a way of vector and matrix multiplications at a reasonable cost. We remark that ${g}_1^{(1)}$ is the only unknown in the equation \eqref{eq:XTDdataonemode2}.

For better understanding, we illustrate the complete discrete form of $\eqref{eq:XTDdataonemode2}$ in the two dimensional case, 
\begin{equation}
\displaystyle
(\delta\boldsymbol{g}_1^{(1)})^T\boldsymbol{\mathcal{U}}\boldsymbol{g}_2^{(1)}=(\delta\boldsymbol{g}_1^{(1)})^T\boldsymbol{g}_1^{(1)}(\boldsymbol{g}_2^{(1)})^T\boldsymbol{g}_2^{(1)}
\end{equation}
Thus,
\begin{equation}
\displaystyle
\boldsymbol{g}_1^{(1)}=\boldsymbol{\mathcal{U}}\boldsymbol{g}_2^{(1)}((\boldsymbol{g}_2^{(1)})^T\boldsymbol{g}_2^{(1)})^{-1}
\end{equation}

By solving \eqref{eq:XTDdataonemode2}, the estimate of $g_1^{(1)}$ is obtained. Similarly, we can compute the estimate for $g_2^{(1)}$ by considering  $g^{(1)}_i,\ \forall i\neq 2$ are constant and  $\Tilde{{u}}^{(1)}=0$, $\delta g^{(1)}_i=0,\ \forall i\neq 2$, and $\delta \Tilde{{u}}^{(1)}=0$.  For the remaining functions, the procedure is the same. This is known as the alternating fixed point algorithm. Once all the first estimates for $g_i^{(1)}$ are computed, we can iterate several times until $\|\Delta g_1^{(1)}\cdots g_n^{(1)}\|$ converges. 

Once the first mode converges, we can go to the next mode.  For $M\geq2,K=1$,
\begin{equation}
\displaystyle
u^{\text{XTD}}=\sum_{m=1}^{M-1} g_1^{(m)}\cdots g_n^{(m)}+g_1^{(M)}\cdots g_n^{(M)} + \Tilde{{u}}^{(1)}
\end{equation}
\begin{equation}
\displaystyle
\delta u^{\text{XTD}}=\delta g_1^{(M)}\cdots g_n^{(M)}+\cdots+g_1^{(M)}\cdots\delta g_n^{(M)} + \delta \Tilde{{u}}^{(1)}
\end{equation}

Analogically to the first mode, we can solve for $g_1^{(M)},\dots ,g_n^{(M)}$. As an example, the 2D formulation for $g_1^{(M)}$ reads
\begin{equation}
\displaystyle
\boldsymbol{g}_1^{(M)}=(\boldsymbol{\mathcal{U}}-\sum_{m=1}^{M-1}
\boldsymbol{g}_1^{(m)}\otimes\boldsymbol{g}_2^{(m)})\boldsymbol{g}_2^{(M)}((\boldsymbol{g}_2^{(M)})^T\boldsymbol{g}_2^{(M)})^{-1}
\end{equation}

The number of separated mode $M$ can increase until a predefined criterion is reached, i.e. 
\begin{equation}
\displaystyle
\label{eq:convergsepmode}
\frac{\|u-\sum_{m=1}^{M} g_1^{(m)}\cdots g_n^{(m)}\|_\infty}{\|u\|_\infty}\leq\epsilon_M
\end{equation}

{Note that the above algorithm corresponds to an alternating least squares method \cite{ten1987some,zhang2001rank} for finding rank-one approximations of tensors. The orthogonality of the final separated modes is usually not guaranteed, except in 2D cases. Under the assumption of  existence of orthogonal decomposition (not necessarily exists for high order tensors), the so-called generalized Rayleigh-Newton iteration \cite{zhang2001rank}, can be used, which requires the Jacobian matrix of \eqref{eq:XTDdatascalar} and yields higher convergence rates than the alternating least squares method. Other alternatives, such as the Jacobi-version of alternating least squares \cite{zhang2001rank}, can be used when parallel computing is desired for this procedure.}

\subsection{Update the extended mode $\Tilde{{u}}^{(1)}$}
In XTD, this criterion $\epsilon_M$ is not the final approximation accuracy. It is defined with  a relatively larger value (e.g. 0.1). Once the stopping criterion is reached, we can update the first extended mode  $\Tilde{{u}}^{(1)}$. The method is very straightforward and can be written as follows
\begin{equation}
\displaystyle
\label{eq:XTDdatanonsep}
\begin{aligned}
&\Tilde{{u}}^{(1)}=u - \sum_{m=1}^{M} g_1^{(m)}\cdots g_n^{(m)},\quad \text{if}\quad ({\mu_1},\dots,\mu_n)\in \mathcal{I}\\
&\Tilde{{u}}^{(1)}=0,\quad \text{otherwise}
\end{aligned}
\end{equation}
where $\mathcal{I}$ is the index of $l$ maximum values of the difference $|u-\sum_{m=1}^{M} g_1^{(m)}\cdots g_n^{(m)}|$. To do so, we might define the following operator
\begin{equation}
\displaystyle
\mathcal{I}=\{({\mu_1},\dots,\mu_n)\in \Omega\ |\ |\mathcal{I}|=l, |u-\sum_{m=1}^{M} g_1^{(m)}\cdots g_n^{(m)}|\geq \sup_{\Omega\backslash\mathcal{I}}|u-\sum_{m=1}^{M} g_1^{(m)}\cdots g_n^{(m)}|\}
\end{equation}
The value of $l$ is at choice, usually much smaller than the vector size of $\boldsymbol{g}_i$. This ensures the sparsity of the extended mode. 
\subsection{Further enrich modes until convergence }
Once it is done, the number of separated modes $M$ can continue increasing until another stopping point. This means computing new separated modes  with the updated extended mode $\Tilde{{u}}^{(1)}$.    A new threshold in \eqref{eq:convergsepmode} can be used for the new stopping point. New extended modes can be computed as well. The final convergence criterion is defined as
\begin{equation}
\displaystyle
\label{eq:convergemultimode}
\frac{\|u-\sum_{m=1}^{M} g_1^{(m)}\cdots g_n^{(m)}-\sum_{k=1}^{K}\Tilde{{u}}^{(k)}\|_\infty}{\|u\|_\infty}\leq\epsilon_{XTD}
\end{equation}
where $\epsilon_{XTD}$ stands for the final approximation accuracy.

For summary, the overall solution procedure is described as follows.
\begin{itemize}
    \item Step 1: $  M=1,K=1 $.
    \item Step 2: Initialization for $g_i^{(M)}$ and $\Tilde{{u}}^{(K)}=0$
    \item Step 3: Update alternatively  the $g_i^{(M)}$ for $i=1,n$ using the alternating fixed point algorithm by considering  $\delta g_j^{(M)}=0, \forall j\neq i$ and $\delta\Tilde{{u}}^{(K)}=0$. If $M\geq2$, $g_j^{(m)}$ are constant for $m<M$.
    \item Step 4: If $\|\Delta g_1^{(M)}\cdots g_n^{(M)}\|$ converges, go to Step 5. Otherwise, repeat Step 3.
    \item Step 5: Check the convergence with \eqref{eq:convergsepmode}. If so, update $\Tilde{{u}}^{(K)}$ with \eqref{eq:XTDdatanonsep} and go to Step 6. Otherwise,  go to next modes $M\leftarrow M+1$ and return to Step 2.
    \item Step 6: Check the convergence with \eqref{eq:convergemultimode}. If so, stop computations. Otherwise, $M\leftarrow M+1$, $K\leftarrow K+1$, $\epsilon_M\leftarrow \epsilon_M/c$  and return to Step 2.
\end{itemize}
We remark that the factor $c>1$ is used in Step 6 to decrease the convergence threshold for a new stopping point for the separated modes.

\subsection{Numerical results}
\begin{figure}[htbp]
\centering
\subfigure[Original data function ]{\includegraphics[scale=0.23]{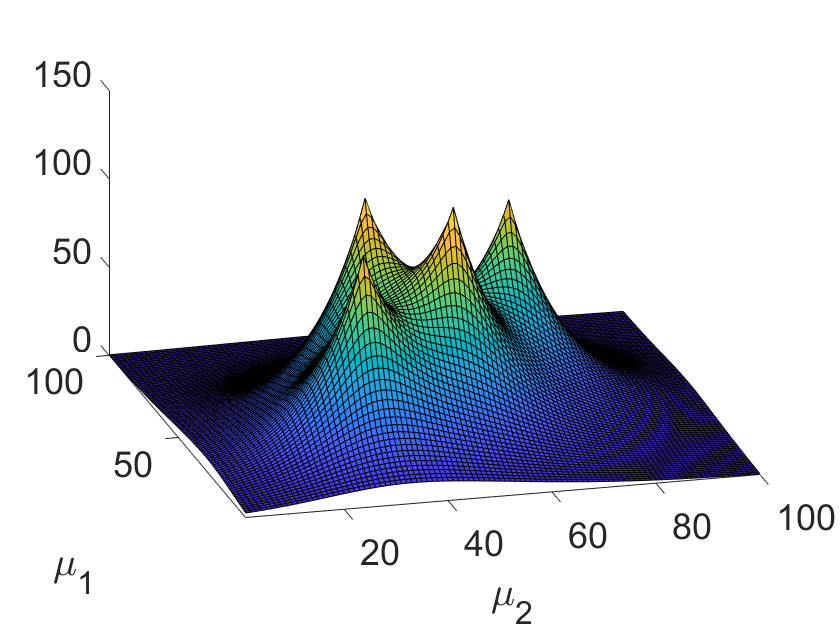}}\\
\subfigure[PGD 5 modes ]{\includegraphics[scale=0.23]{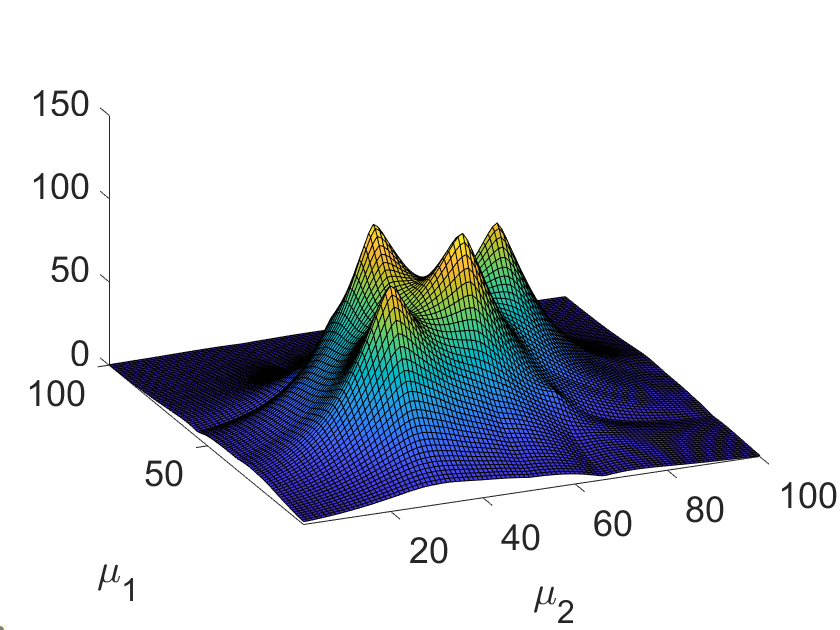}}
\subfigure[PGD 6 modes ]{\includegraphics[scale=0.23]{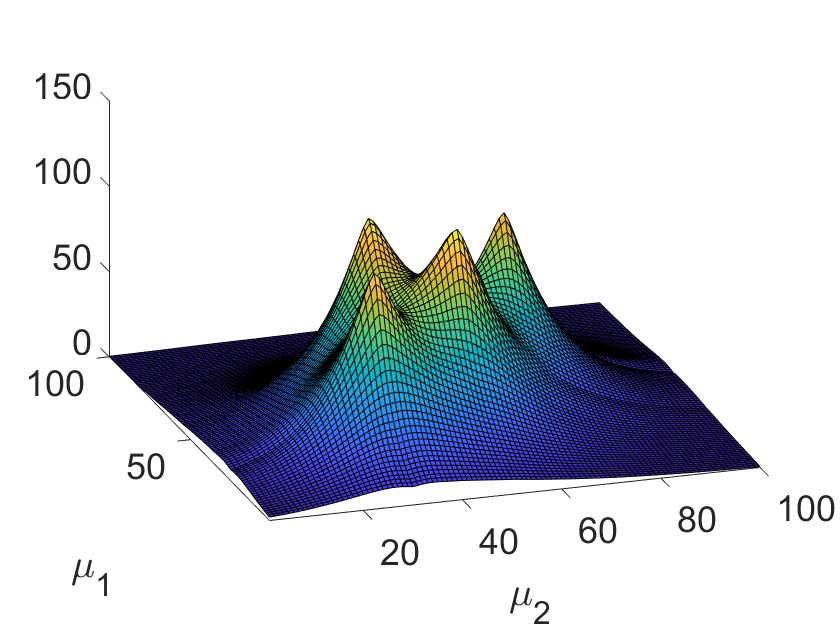}}
\subfigure[PGD 8 modes ]{\includegraphics[scale=0.23]{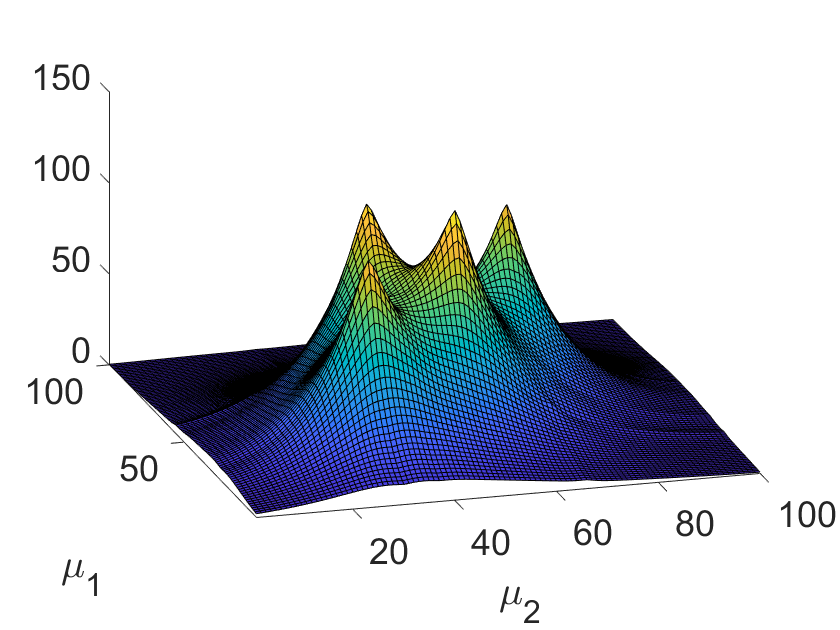}}\\
\subfigure[XTD 5 modes ]{\includegraphics[scale=0.23]{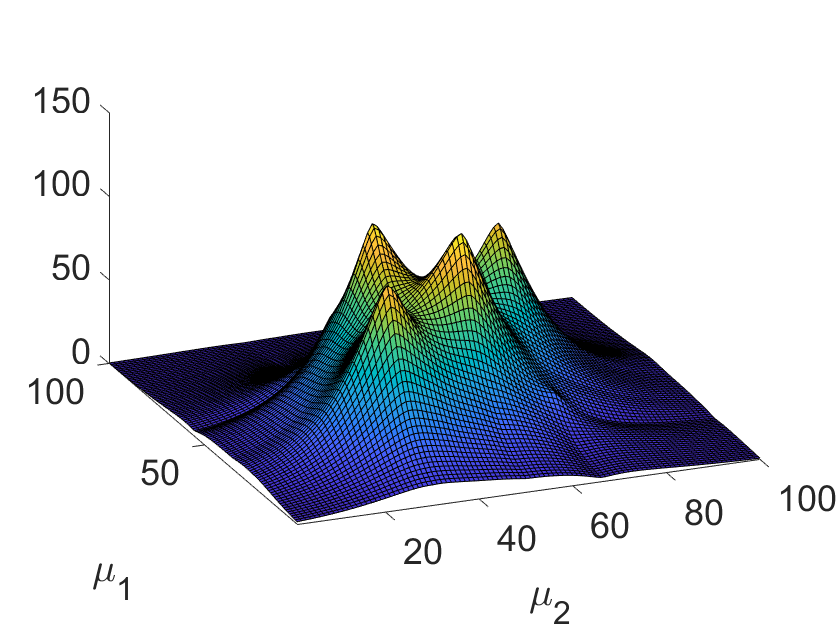}}
\subfigure[XTD 5 modes + 1 ext.  ]{\includegraphics[scale=0.23]{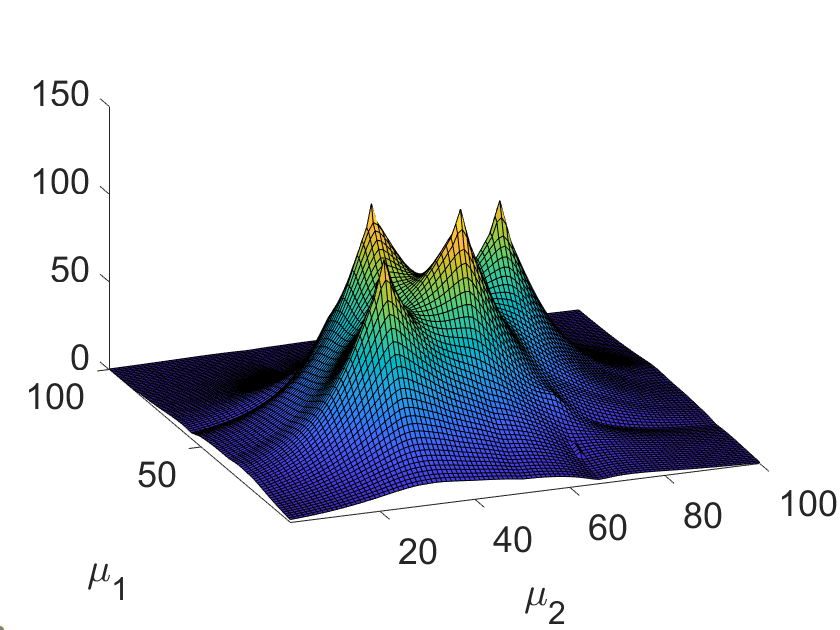}}
\subfigure[XTD 7 modes + 1 ext.  ]{\includegraphics[scale=0.23]{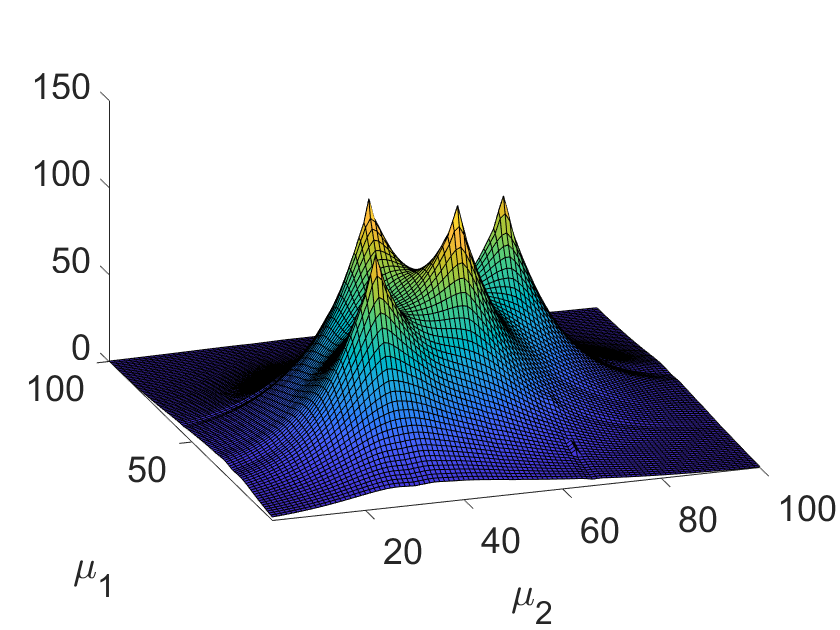}}

\caption{Data reproduced by XTD with comparison to PGD. Note: ext. standards for extended mode. }
\label{fig:xtddatafunction}
\end{figure}

We will illustrate a numerical example for a two dimensional function. The function for generating the data is written as 
\begin{equation}
\displaystyle
\begin{aligned}
u(\mu_1,\mu_2)=&100e^{\sqrt{(\mu_1-50)^2+(\mu_2-37)^2}/10}+100e^{\sqrt{(\mu_1-35)^2+(\mu_2-50)^2}/10}\\
&+100e^{\sqrt{(\mu_1-25)^2+(\mu_2-30)^2}/10}+100e^{\sqrt{(\mu_1-55)^2+(\mu_2-66)^2}/10}
\end{aligned}
\end{equation}
where the parameters $\mu_1\in [1,\ 100]$, $\mu_2\in [1,\ 100]$. The data  is then generated by evaluating the function $u$ for a uniform mesh of $100\times100$ in the 2D
domain $\Omega=\Omega_{\mu_1}\times\Omega_{\mu_2}$. \figurename~\ref{fig:xtddatafunction}(a) depicts the shape of this highly nonlinear function on the $100\times100$ mesh. We can observe the several sharp points around some predefined locations. This is the singularity that can challenge the reducibility of the function. 

\begin{figure}[htbp]
\centering
{\includegraphics[scale=0.4]{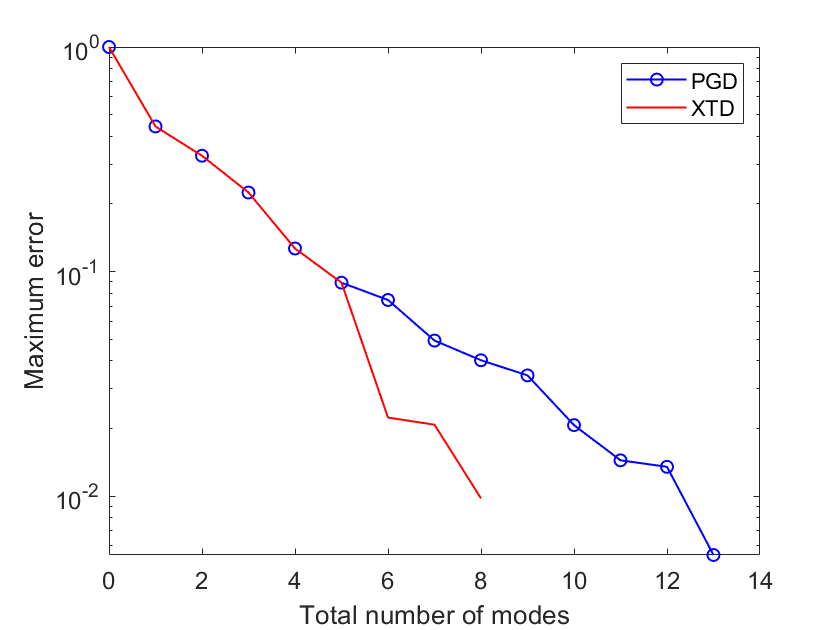}}
\caption{Evolution of approximation error against the number of modes}
\label{fig:xtdvspgd}
\end{figure}

If we use PGD (HOPGD, also equivalent to singular value decomposition in 2D cases) to reproduce this function, it requires 13 separated modes for an error less than $1\%$. However, if the XTD is used, only 7 separated modes and 1 sparse extended mode are needed for the same level of accuracy.  The evolution of error against the total number of modes is illustrated in \figurename~\ref{fig:xtdvspgd}. The first enrichment of the extended mode happens after 5 separated mode for the XTD method. We observe a significant decrease of error with this enrichment. This improvement is also clearly observed in the shape of the function reproduced by XTD (\figurename~\ref{fig:xtddatafunction}). It is shown that the separated modes can reproduce a smooth function, whereas the extended modes can efficiently enrich the local singularity. consequently, the  number of mode is reduced. 

\begin{figure}[htbp]
\centering
{\includegraphics[scale=0.4]{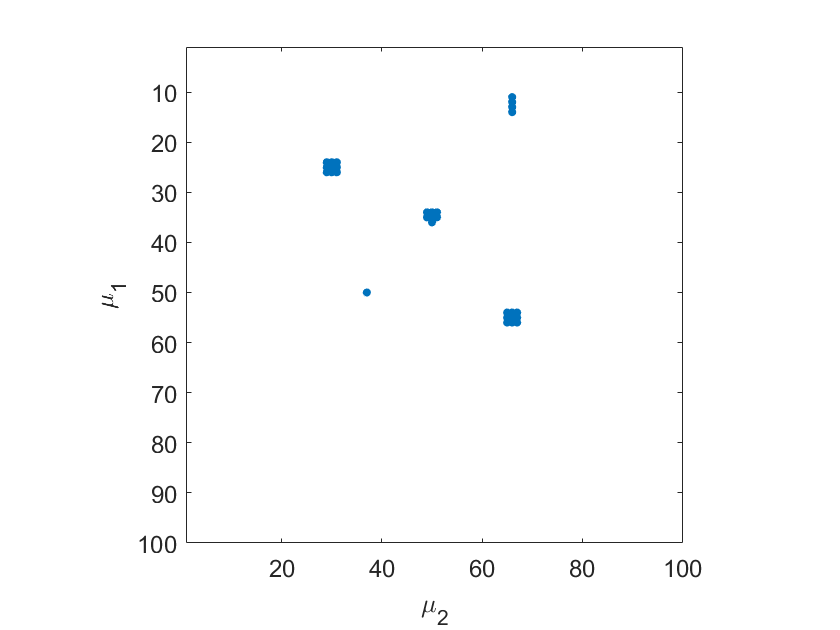}}
\caption{Location of non-zero terms for the extended mode. The sparsity (percentage of zero-valued elements) is 99.7$\%$. }
\label{fig:xtddatanonsep}
\end{figure}

\figurename~\ref{fig:xtddatanonsep} illustrates the location of the non-zero terms in the sparse extended mode. In this example, the initial  $\epsilon_M=10\%$, the number of local enrichment points $l=30$ which is much smaller than the total data size $100\times100$ and the size of one mode (1 pair of separated functions): $100+100$. Storing this sparse matrix is much cheaper than the modes.

This example shows that the singularity can  challenge the reducibility of given functions. This implies that the number of separated modes for a good approximation can be very high. Using the XTD method can significantly improve the approximation in terms of accuracy and separability. 

{\section{Implementation of the XTD  model reduction}\label{XTDmodel}}
For a better illustration, let us consider a problem with two design parameters. The displacement increment at $t_k$ reads
\begin{equation}
\displaystyle
\label{eq:XTDdu2d}
\Delta\ud(\X,{\mu_1},\mu_2)= \sum_{m=1}^{M}\boldsymbol{a}^{(m)}(\X)g_1^{(m)}(\mu_1)g_2^{(m)}(\mu_2) + \sum_{k=1}^{K}\Tilde{\boldsymbol{u}}^{(k)}(\X,{\mu_1},\mu_2)
\end{equation}
Given the  displacement $\ud_0$ at the previous step, the full displacement is then
\begin{equation}
\displaystyle
\ud(\X,{\mu_1},\mu_2)=\ud_0(\X,{\mu_1},\mu_2)+ \Delta\ud(\X,{\mu_1},\mu_2)
\end{equation}
and for the test function 
\begin{equation}
\displaystyle
{{\ud}^{*}}=\Delta\ud^*
\end{equation}
The weak form \eqref{eq:weakform} of the problem  becomes
\begin{equation}
\displaystyle
-\int_{\Omega_\X\times\Omega_{\mu_1}\times\Omega_{\mu_2}}{\sig(\ud_0+\Delta\ud) :\eps ({\Delta{\ud}^{*}})d\X d\mu_1d\mu_2}+\int_{\partial_{\textbf{F}} {{\Omega_\X}}\times\Omega_{\mu_1}\times\Omega_{\mu_2}}{\bar{\textbf{d}}\cdot{\Delta{\ud}^{*}}dsd\mu_1d\mu_2=0}
\end{equation}
Considering that $\bar{\textbf{d}}=\bar{\textbf{d}}_0+\Delta\bar{\textbf{d}}$ with $\bar{\textbf{d}}_0$ denoting the previous loading, we have 
\begin{equation}
\displaystyle
-\int_{\Omega_\X\times\Omega_{\mu_1}\times\Omega_{\mu_2}}{\sig(\ud_0) :\eps ({\Delta{\ud}^{*}})d\X d\mu_1d\mu_2}+\int_{\partial_{\textbf{F}} {{\Omega_\X}}\times\Omega_{\mu_1}\times\Omega_{\mu_2}}{\bar{\textbf{d}}_0\cdot{\Delta{\ud}^{*}}ds d\mu_1d\mu_2=0}
\end{equation}
Therefore, the following equation should be satisfied
\begin{equation}
\displaystyle
-\int_{\Omega_\X\times\Omega_{\mu_1}\times\Omega_{\mu_2}}{\Delta\sig :\eps ({\Delta{\ud}^{*}})d\X d\mu_1d\mu_2}+\int_{\partial_{\textbf{F}} {{\Omega_\X}}\times\Omega_{\mu_1}\times\Omega_{\mu_2}}{\Delta\bar{\textbf{d}}\cdot{\Delta{\ud}^{*}}ds d\mu_1d\mu_2=0}
\end{equation}
With $\Delta\sig=\textbf{D}:(\Delta\eps -\Delta\epsp)$ and $\Omega=\Omega_\X\times\Omega_{\mu_1}\times\Omega_{\mu_2}$ and $\partial_{\textbf{F}}\Omega=\partial_{\textbf{F}} {{\Omega_\X}}\times\Omega_{\mu_1}\times\Omega_{\mu_2}$
\begin{equation}
\displaystyle
\label{eq:weakincrementform}
-\int_{\Omega}{\textbf{D}:(\Delta\eps -\Delta\epsp):\eps ({\Delta{\ud}^{*}})d\X d\mu_1d\mu_2}+\int_{\partial_{\textbf{F}} {{\Omega}}}{\Delta\bar{\textbf{d}}\cdot{\Delta{\ud}^{*}}ds d\mu_1d\mu_2=0}
\end{equation}

\subsection{Compute initial separated modes with $\Tilde{\boldsymbol{u}}^{(k)}=0$ \label{initialsepmode}}
Now, we can start from computing some separated modes of XTD by solving the above problem. Following the incremental strategy explained in \ref{XTDdata}, we start from
$M=1, K=1$. In this case, 
\begin{equation}
\displaystyle
\Delta\ud(\X,{\mu_1},\mu_2)= \boldsymbol{a}^{(1)}(\X)g_1^{(1)}(\mu_1)g_2^{(1)}(\mu_2) + \Tilde{\boldsymbol{u}}^{(1)}
\end{equation}
\begin{equation}
\displaystyle
\Delta\ud^*= \boldsymbol{a}^{(1)*}g_1^{(1)}g_2^{(1)}+\boldsymbol{a}^{(1)}g_1^{(1)*}g_2^{(1)}+\boldsymbol{a}^{(1)}g_1^{(1)}g_2^{(1)*}+ \Tilde{\boldsymbol{u}}^{(1)*}
\end{equation}
Given the initial $\Tilde{\boldsymbol{u}}^{(1)}=0$, then
\begin{equation}
\displaystyle
\Delta\ud(\X,{\mu_1},\mu_2)= \boldsymbol{a}^{(1)}(\X)g_1^{(1)}(\mu_1)g_2^{(1)}(\mu_2) 
\end{equation}
\begin{equation}
\displaystyle
\Delta\ud^*= \boldsymbol{a}^{(1)*}g_1^{(1)}g_2^{(1)}+\boldsymbol{a}^{(1)}g_1^{(1)*}g_2^{(1)}+\boldsymbol{a}^{(1)}g_1^{(1)}g_2^{(1)*}
\end{equation}
For notation simplification, we omit the superscripts
\begin{equation}
\displaystyle
\Delta\ud(\X,{\mu_1},\mu_2)= \boldsymbol{a}(\X)g_1(\mu_1)g_2(\mu_2) 
\end{equation}
\begin{equation}
\displaystyle
\Delta\ud^*= \boldsymbol{a}^{*}g_1 g_2+\boldsymbol{a} g_1^{*}g_2+\boldsymbol{a}g_1g_2^{*}
\end{equation}
Now we solve first $\boldsymbol{a}$ by fixing $g_1, g_2$ as a constant. Hence, $g_1^{*}=0, g_2^{*}=0$.
\begin{equation}
\displaystyle
\Delta\ud^*= \boldsymbol{a}^{*}g_1g_2
\end{equation}
Assuming $\Delta\epsp=0$ and considering the FE discretization for $\boldsymbol{a}(\X)$
\begin{equation}
\displaystyle
\boldsymbol{a}(\X)=\textbf{N}(\X)\boldsymbol{\mathcal{A}} 
\end{equation}
The \eqref{eq:weakincrementform} becomes
\begin{equation}
\displaystyle
-\int_{\Omega}{g_2^{T} g_1^{T}\boldsymbol{\mathcal{A}}^{*T}\textbf{B}^T\textbf{D}\textbf{B}\boldsymbol{\mathcal{A}}g_1g_2\ d\X d\mu_1d\mu_2}+\int_{\partial_{\textbf{F}} {{\Omega}}}{g_2^{T} g_1^{T}\boldsymbol{\mathcal{A}}^{*T}\textbf{N}^T\Delta\bar{\textbf{d}}\ ds d\mu_1d\mu_2=0}
\end{equation}
Here we assume that the stiffness  $\textbf{D}$ is independent of design parameters, and the external load has been separated and can be approximated by one mode, i.e. $\Delta\bar{\textbf{d}}=\textbf{f}_\X(\X) \text{f}_{\mu_1}({\mu_1})\text{f}_{\mu_2}({\mu_2})$. By rearrangement, the above equation reads
\begin{equation}
\displaystyle
\begin{aligned}
\int_{\Omega_\X}{ \boldsymbol{\mathcal{A}}^{*T}\textbf{B}^T\textbf{D}\textbf{B}\boldsymbol{\mathcal{A}}\ d\X \int_{\Omega_{\mu_1}}g_1^{T}g_1\ d\mu_1 \int_{\Omega_{\mu_2}} g_2^{T}g_2\  d\mu_2}\\
=\int_{\partial_{\textbf{F}} {{\Omega_\X}}}{ \boldsymbol{\mathcal{A}}^{*T}\textbf{N}^T\textbf{f}_\X\ ds \int_{\Omega_{\mu_1}} g_1^{T}\text{f}_{\mu_1} \ d\mu_1 \int_{\Omega_{\mu_2}} g_2^{T}\text{f}_{\mu_2}\  d\mu_2}
\end{aligned}
\end{equation}
The final discretized form for solving $\boldsymbol{a}$ is obtained as
\begin{equation}
\displaystyle
\begin{aligned}
\textbf{K}\boldsymbol{\mathcal{A}} \boldsymbol{g}_1^{T}\boldsymbol{g}_1 \boldsymbol{g}_2^{T}\boldsymbol{g}_2
={\textbf{F}_\X }{  \boldsymbol{g}_1^{T}\textbf{F}_{\mu_1}  \boldsymbol{g}_2^{T}\textbf{F}_{\mu_2} }
\end{aligned}
\end{equation}
and
\begin{equation}
\displaystyle
\begin{aligned}
\boldsymbol{\mathcal{A}} 
=\textbf{K}^{-1}{\textbf{F}_\X }{  \boldsymbol{g}_1^{T}\textbf{F}_{\mu_1}  \boldsymbol{g}_2^{T}\textbf{F}_{\mu_2} }(\boldsymbol{g}_1^{T}\boldsymbol{g}_1 \boldsymbol{g}_2^{T}\boldsymbol{g}_2)^{-1}
\end{aligned}
\end{equation}
with 
\begin{equation}
\displaystyle
\begin{aligned}
&\textbf{F}_\X =\int_{\partial_{\textbf{F}} {{\Omega_\X}}}{ \textbf{N}^T\textbf{f}_\X\ ds }\\
&\textbf{F}_{\mu_1} = \text{f}_{\mu_1}(\Omega_{\mu_1}^h)  \\ 
&\textbf{F}_{\mu_2}  = \text{f}_{\mu_2}(\Omega_{\mu_2}^h)
\end{aligned}
\end{equation}
where $\Omega_{\mu_1}^h$ and $\Omega_{\mu_2}^h$ denote respectively the supporting mesh in each parameter domain.

After updating $\boldsymbol{a}$, we solve then $g_1$ by fixing $\boldsymbol{a}, g_2$ as a constant. Hence, $\boldsymbol{a}^{*}=0, g_2^{*}=0$.
\begin{equation}
\displaystyle
\Delta\ud^*= \boldsymbol{a}g_1^{*}g_2
\end{equation}
Following the same concept as previously, the final discretized form reads
\begin{equation}
\displaystyle
\begin{aligned}
\boldsymbol{\mathcal{A}}^T\textbf{K}\boldsymbol{\mathcal{A}} \boldsymbol{g}_1 \boldsymbol{g}_2^{T}\boldsymbol{g}_2
=\boldsymbol{\mathcal{A}}^T{\textbf{F}_\X }{  \textbf{F}_{\mu_1}  \boldsymbol{g}_2^{T}\textbf{F}_{\mu_2} }
\end{aligned}
\end{equation}
and
\begin{equation}
\displaystyle
\begin{aligned}
\boldsymbol{g}_1
=(\boldsymbol{\mathcal{A}}^T\textbf{K}\boldsymbol{\mathcal{A}})^{-1}\boldsymbol{\mathcal{A}}^T{\textbf{F}_\X }{  \textbf{F}_{\mu_1}  \boldsymbol{g}_2^{T}\textbf{F}_{\mu_2} }( \boldsymbol{g}_2^{T}\boldsymbol{g}_2)^{-1}
\end{aligned}
\end{equation}

Similarly we solve  $g_2$ by fixing $\boldsymbol{a}, g_1$ as a constant. The final discrete equation is
\begin{equation}
\displaystyle
\begin{aligned}
\boldsymbol{g}_2
=(\boldsymbol{\mathcal{A}}^T\textbf{K}\boldsymbol{\mathcal{A}})^{-1}\boldsymbol{\mathcal{A}}^T{\textbf{F}_\X }{  \boldsymbol{g}_1^{T}\textbf{F}_{\mu_1}  \textbf{F}_{\mu_2} }( \boldsymbol{g}_1^{T}\boldsymbol{g}_1)^{-1}
\end{aligned}
\end{equation}
The alternating solution procedure is iteratively performed until $\|\Delta \boldsymbol{\mathcal{A}}\otimes \boldsymbol{g}_1 \otimes \boldsymbol{g}_2\|$ converges. 

When the first mode is obtained, we can continue to solve for more modes. In the case of $M\geq2$, the $M-1$ modes have been known, the current unknown is again denoted by $\boldsymbol{a},g_1,g_2$. We have 
\begin{equation}
\displaystyle
\Delta\ud(\X,{\mu_1},\mu_2)= \sum_{m=1}^{M-1}\boldsymbol{a}^{(m)}(\X)g_1^{(m)}(\mu_1)g_2^{(m)}(\mu_2) + \boldsymbol{a}(\X)g_1(\mu_1)g_2(\mu_2)
\end{equation}
with the test function
\begin{equation}
\displaystyle
\Delta\ud^*= \boldsymbol{a}^{*}g_1 g_2+\boldsymbol{a} g_1^{*}g_2+\boldsymbol{a}g_1g_2^{*}
\end{equation}
The alternating algorithm gives 
\begin{equation}
\displaystyle
\begin{aligned}
\boldsymbol{\mathcal{A}} 
=\ &\textbf{K}^{-1}{\textbf{F}_\X }{  \boldsymbol{g}_1^{T}\textbf{F}_{\mu_1}  \boldsymbol{g}_2^{T}\textbf{F}_{\mu_2} }(\boldsymbol{g}_1^{T}\boldsymbol{g}_1 \boldsymbol{g}_2^{T}\boldsymbol{g}_2)^{-1}\\
&-\textbf{K}^{-1}\sum_{m=1}^{M-1}\textbf{K}\boldsymbol{\mathcal{A}}^{(m)} \boldsymbol{g}_1^{T}\boldsymbol{g}_1^{(m)}\boldsymbol{g}_2^{T}\boldsymbol{g}_2^{(m)}(\boldsymbol{g}_1^{T}\boldsymbol{g}_1 \boldsymbol{g}_2^{T}\boldsymbol{g}_2)^{-1}
\end{aligned}
\end{equation}
\begin{equation}
\displaystyle
\begin{aligned}
\boldsymbol{g}_1
=(&\boldsymbol{\mathcal{A}}^T\textbf{K}\boldsymbol{\mathcal{A}})^{-1}\boldsymbol{\mathcal{A}}^T{\textbf{F}_\X }{  \textbf{F}_{\mu_1}  \boldsymbol{g}_2^{T}\textbf{F}_{\mu_2} }( \boldsymbol{g}_2^{T}\boldsymbol{g}_2)^{-1}\\
&-(\boldsymbol{\mathcal{A}}^T\textbf{K}\boldsymbol{\mathcal{A}})^{-1}\sum_{m=1}^{M-1}\boldsymbol{\mathcal{A}}^T\textbf{K}\boldsymbol{\mathcal{A}}^{(m)}\boldsymbol{g}_1^{(m)}\boldsymbol{g}_2^{T}\boldsymbol{g}_2^{(m)}(\boldsymbol{g}_2^{T}\boldsymbol{g}_2)^{-1}
\end{aligned}
\end{equation}
\begin{equation}
\displaystyle
\begin{aligned}
\boldsymbol{g}_2
=(&\boldsymbol{\mathcal{A}}^T\textbf{K}\boldsymbol{\mathcal{A}})^{-1}\boldsymbol{\mathcal{A}}^T{\textbf{F}_\X }{  \boldsymbol{g}_1^{T}\textbf{F}_{\mu_1}  \textbf{F}_{\mu_2} }( \boldsymbol{g}_1^{T}\boldsymbol{g}_1)^{-1}\\
&-(\boldsymbol{\mathcal{A}}^T\textbf{K}\boldsymbol{\mathcal{A}})^{-1}\sum_{m=1}^{M-1}\boldsymbol{\mathcal{A}}^T\textbf{K}\boldsymbol{\mathcal{A}}^{(m)}\boldsymbol{g}_1^{T}\boldsymbol{g}_1^{(m)}\boldsymbol{g}_2^{(m)}(\boldsymbol{g}_1^{T}\boldsymbol{g}_1)^{-1}
\end{aligned}
\end{equation}
The number of modes keeps increasing until the convergence  
\begin{equation}
\displaystyle
\label{eq:XTDROMmodeconverg}
\frac{\|\boldsymbol{\mathcal{A}}\otimes \boldsymbol{g}_1 \otimes \boldsymbol{g}_2\|_\infty}{\|\Delta\ud\|_\infty}\leq\epsilon_M
\end{equation}
\subsection{Compute the extended mode}
Once the first set of separated modes is computed, we can solve the extended mode by updating the plastic strain $\epsp$ and $\Delta\epsp$.

Taking the current displacement increment as an estimate $\Delta \ud^{(0)}$
\begin{equation}
\displaystyle
\Delta\ud^{(0)}(\X,{\mu_1},\mu_2)= \Delta \ud
\end{equation}
we can check the local plastic region $\Omega_{\text{pl}}$ where $\Delta\epsp(\X,{\mu_1},\mu_2)\neq 0$.  This can be done using standard FE implementation for plasticity with a point-wise evaluation in the global domain $\Omega$.

The Newton-Raphson method is then used to update the plastic strain and the displacement in $\Omega_{\text{pl}}$ with the boundary condition imposed by the remaining region. The updated displacement increment $\Delta \ud$ should satisfy that
\begin{equation}
    \displaystyle
    {{\Fm}_{\text{ext}}(\ud_0+\Delta \ud)-{\Fm}_{\text{int}}(\ud_0+\Delta \ud)=0 \quad \text{in}\quad \Omega_{\text{pl}} \quad \text{with } \quad \Delta\ud|_{\partial \Omega_{\text{pl}} }=\Delta\ud|_{\partial \Omega_{\text{el}} }}
\end{equation}
where ${\Fm}_{\text{ext}}$ and ${\Fm}_{\text{int}}$ are computed with \eqref{eq:FEforce}. This solution can be solved easily in a parallel way. 

The sparse extended mode is obtained with
\begin{equation}
\displaystyle
\Tilde{\boldsymbol{u}}^{(K)}(\X,{\mu_1},\mu_2)= \Delta \ud - \Delta\ud^{(0)}(\X,{\mu_1},\mu_2)
\end{equation}
Since $\Delta \ud$ is locally updated, the extended mode is sparse.

An another outcome of this step is the updated plastic strain and the induced internal force, i.e.
\begin{equation}
\displaystyle
\textbf{F}_{\text{pl}}({\mu_1},\mu_2), \quad \Delta\textbf{F}_{\text{pl}}({\mu_1},\mu_2)
\end{equation}
as computed by \eqref{eq:Fpl}, which is also sparse.

\subsection{Enrich separated modes with $\Tilde{\boldsymbol{u}}^{(k)}\neq0$}
With the updated extended mode, we can keep increasing the number of separated modes $M$. Again, the new mode is denoted by $\boldsymbol{a},g_1,g_2$. 
\begin{equation}
\displaystyle
\Delta\ud(\X,{\mu_1},\mu_2)= \sum_{m=1}^{M-1}\boldsymbol{a}^{(m)}g_1^{(m)}g_2^{(m)} + \boldsymbol{a}g_1g_2 +\sum_{k=1}^{K}\Tilde{\boldsymbol{u}}^{(k)}
\end{equation}
with the test function
\begin{equation}
\displaystyle
\Delta\ud^*= \boldsymbol{a}^{*}g_1 g_2+\boldsymbol{a} g_1^{*}g_2+\boldsymbol{a}g_1g_2^{*}
\end{equation}

Now, we solve $\boldsymbol{a}$  by fixing the others. The \eqref{eq:weakincrementform} becomes
\begin{equation}
\displaystyle
\begin{aligned}
&\int_{\Omega}{g_2^{T} g_1^{T}\boldsymbol{\mathcal{A}}^{*T}\sum_{m=1}^{M-1}\textbf{B}^T \textbf{D}\textbf{B}\boldsymbol{\mathcal{A}}^{(m)}g_1^{(m)}g_2^{(m)}\ d\X d\mu_1d\mu_2}\\
&+\int_{\Omega}{g_2^{T} g_1^{T}\boldsymbol{\mathcal{A}}^{*T}\textbf{B}^T\textbf{D}\textbf{B}\boldsymbol{\mathcal{A}}g_1g_2\ d\X d\mu_1d\mu_2}\\
&+\int_{\Omega}{g_2^{T} g_1^{T}\boldsymbol{\mathcal{A}}^{*T}\sum_{k=1}^{K}\textbf{B}^T \textbf{D}\textbf{B}\Tilde{\boldsymbol{\mathcal{U}}}^{(k)}\ d\X d\mu_1d\mu_2}\\
&-\int_{\Omega_{\boldsymbol{\mu}}}{g_2^{T} g_1^{T}\boldsymbol{\mathcal{A}}^{*T}\Delta\textbf{F}_{\text{pl}}\ d\mu_1d\mu_2}\\
&=\int_{\partial_{\textbf{F}} {{\Omega}}}{g_2^{T} g_1^{T}\boldsymbol{\mathcal{A}}^{*T}\textbf{N}^T\Delta\bar{\textbf{d}}\ ds d\mu_1d\mu_2}
\end{aligned}
\end{equation}
where $\Omega_{\boldsymbol{\mu}}$ denotes $\Omega_{{\mu}_1}\times \Omega_{\mu_2}$. As we can see, $\Tilde{\boldsymbol{\mathcal{U}}}^{(k)}$ and $\Delta\textbf{F}_{\text{pl}}$ are not naturally separated,  the third and the fourth terms in above equation require the so-called full integration. Nevertheless, since they are sparse and the integration can be done by simple vector and matrix multiplications, the associated cost remains still low from authors' experience. 

Thus the solution is 
\begin{equation}
\displaystyle
\begin{aligned}
\boldsymbol{\mathcal{A}} 
=\ &\textbf{K}^{-1}{\textbf{F}_\X }{  \boldsymbol{g}_1^{T}\textbf{F}_{\mu_1}  \boldsymbol{g}_2^{T}\textbf{F}_{\mu_2} }(\boldsymbol{g}_1^{T}\boldsymbol{g}_1 \boldsymbol{g}_2^{T}\boldsymbol{g}_2)^{-1}\\
&-\textbf{K}^{-1}\sum_{m=1}^{M-1}\textbf{K}\boldsymbol{\mathcal{A}}^{(m)} \boldsymbol{g}_1^{T}\boldsymbol{g}_1^{(m)}\boldsymbol{g}_2^{T}\boldsymbol{g}_2^{(m)}(\boldsymbol{g}_1^{T}\boldsymbol{g}_1 \boldsymbol{g}_2^{T}\boldsymbol{g}_2)^{-1}\\
&-\textbf{K}^{-1}\textbf{Q}_{\boldsymbol{\Tilde{\boldsymbol{\mathcal{U}}}\mathcal{A}}}(\boldsymbol{g}_1^{T}\boldsymbol{g}_1 \boldsymbol{g}_2^{T}\boldsymbol{g}_2)^{-1}+\textbf{K}^{-1}\textbf{Q}_{\boldsymbol{{\text{pl}}\mathcal{A}}}(\boldsymbol{g}_1^{T}\boldsymbol{g}_1 \boldsymbol{g}_2^{T}\boldsymbol{g}_2)^{-1}
\end{aligned}
\end{equation}
with
\begin{equation}
\displaystyle
\begin{aligned}
&\textbf{Q}_{\boldsymbol{\Tilde{\boldsymbol{\mathcal{U}}}\mathcal{A}}}=\int_{\Omega}{g_2^{T} g_1^{T}\sum_{k=1}^{K}\textbf{B}^T \textbf{D}\textbf{B}\Tilde{\boldsymbol{\mathcal{U}}}^{(k)}\ d\X d\mu_1d\mu_2}\\
&\textbf{Q}_{\boldsymbol{{\text{pl}}\mathcal{A}}}=\int_{\Omega_{\boldsymbol{\mu}}}{g_2^{T} g_1^{T}\Delta\textbf{F}_{\text{pl}}\ d\mu_1d\mu_2}
\end{aligned}
\end{equation}

Analogically, the alternating algorithm gives 
\begin{equation}
\displaystyle
\begin{aligned}
\boldsymbol{g}_1
=(&\boldsymbol{\mathcal{A}}^T\textbf{K}\boldsymbol{\mathcal{A}})^{-1}\boldsymbol{\mathcal{A}}^T{\textbf{F}_\X }{  \textbf{F}_{\mu_1}  \boldsymbol{g}_2^{T}\textbf{F}_{\mu_2} }( \boldsymbol{g}_2^{T}\boldsymbol{g}_2)^{-1}\\
&-(\boldsymbol{\mathcal{A}}^T\textbf{K}\boldsymbol{\mathcal{A}})^{-1}\sum_{m=1}^{M-1}\boldsymbol{\mathcal{A}}^T\textbf{K}\boldsymbol{\mathcal{A}}^{(m)}\boldsymbol{g}_1^{(m)}\boldsymbol{g}_2^{T}\boldsymbol{g}_2^{(m)}(\boldsymbol{g}_2^{T}\boldsymbol{g}_2)^{-1}\\
&-(\boldsymbol{\mathcal{A}}^T\textbf{K}\boldsymbol{\mathcal{A}})^{-1}\textbf{Q}_{\boldsymbol{\Tilde{\boldsymbol{\mathcal{U}}}\boldsymbol{g}_1}}( \boldsymbol{g}_2^{T}\boldsymbol{g}_2)^{-1}+\textbf{K}^{-1}\textbf{Q}_{\boldsymbol{{\text{pl}}\boldsymbol{g}_1}}(\boldsymbol{g}_2^{T}\boldsymbol{g}_2)^{-1}
\end{aligned}
\end{equation}
with
\begin{equation}
\displaystyle
\begin{aligned}
&\textbf{Q}_{\boldsymbol{\Tilde{\boldsymbol{\mathcal{U}}}\boldsymbol{g}_1}}=\int_{\Omega}{g_2^{T} \boldsymbol{\mathcal{A}}^{T}\sum_{k=1}^{K}\textbf{B}^T \textbf{D}\textbf{B}\Tilde{\boldsymbol{\mathcal{U}}}^{(k)}\ d\X d\mu_1d\mu_2}\\
&\textbf{Q}_{\boldsymbol{{\text{pl}}\boldsymbol{g}_1}}=\int_{\Omega_{{\mu}_2}}{g_2^{T} \boldsymbol{\mathcal{A}}^{T}\Delta\textbf{F}_{\text{pl}}\ d\mu_2}
\end{aligned}
\end{equation}

and
\begin{equation}
\displaystyle
\begin{aligned}
\boldsymbol{g}_2
=(&\boldsymbol{\mathcal{A}}^T\textbf{K}\boldsymbol{\mathcal{A}})^{-1}\boldsymbol{\mathcal{A}}^T{\textbf{F}_\X }{  \boldsymbol{g}_1^{T}\textbf{F}_{\mu_1}  \textbf{F}_{\mu_2} }( \boldsymbol{g}_1^{T}\boldsymbol{g}_1)^{-1}\\
&-(\boldsymbol{\mathcal{A}}^T\textbf{K}\boldsymbol{\mathcal{A}})^{-1}\sum_{m=1}^{M-1}\boldsymbol{\mathcal{A}}^T\textbf{K}\boldsymbol{\mathcal{A}}^{(m)}\boldsymbol{g}_1^{T}\boldsymbol{g}_1^{(m)}\boldsymbol{g}_2^{(m)}(\boldsymbol{g}_1^{T}\boldsymbol{g}_1)^{-1}\\
&-(\boldsymbol{\mathcal{A}}^T\textbf{K}\boldsymbol{\mathcal{A}})^{-1}\textbf{Q}_{\boldsymbol{\Tilde{\boldsymbol{\mathcal{U}}}\boldsymbol{g}_2}}( \boldsymbol{g}_1^{T}\boldsymbol{g}_1)^{-1}+\textbf{K}^{-1}\textbf{Q}_{\boldsymbol{{\text{pl}}\boldsymbol{g}_2}}(\boldsymbol{g}_1^{T}\boldsymbol{g}_1)^{-1}
\end{aligned}
\end{equation}
with
\begin{equation}
\displaystyle
\begin{aligned}
&\textbf{Q}_{\boldsymbol{\Tilde{\boldsymbol{\mathcal{U}}}\boldsymbol{g}_2}}=\int_{\Omega}{g_1^{T} \boldsymbol{\mathcal{A}}^{T}\sum_{k=1}^{K}\textbf{B}^T \textbf{D}\textbf{B}\Tilde{\boldsymbol{\mathcal{U}}}^{(k)}\ d\X d\mu_1d\mu_2}\\
&\textbf{Q}_{\boldsymbol{{\text{pl}}\boldsymbol{g}_2}}=\int_{\Omega_{{\mu}_1}}{g_1^{T} \boldsymbol{\mathcal{A}}^{T}\Delta\textbf{F}_{\text{pl}}\ d\mu_1}
\end{aligned}
\end{equation}

A similar iterative procedure to \ref{initialsepmode} is then applied for this enrichment step.
The converged number of modes for this step is also determined by \ref{eq:XTDROMmodeconverg}.

\subsection{Overall convergence criterion}
The enrichment can be performed alternatively for the separated modes and extended modes until the following global convergence criterion is satisfied.
\begin{equation}
\displaystyle
\label{eq:XTDROMallconverg}
\frac{\|\Tilde{\boldsymbol{u}}^{(K)}\|_\infty}{\|\Delta\ud\|_\infty}\leq\epsilon_{\text{XTD}}
\end{equation}

%{\section{Formulation of the coupled XTD-SCA method}\label{SCA}}

\bibliographystyle{model1-num-names}
\bibliography{XTD.bib}

\begin{thebibliography}{38}
\expandafter\ifx\csname natexlab\endcsname\relax\def\natexlab#1{#1}\fi
\providecommand{\bibinfo}[2]{#2}
\ifx\xfnm\relax \def\xfnm[#1]{\unskip,\space#1}\fi
%Type = Article
\bibitem[{WANG and SHAN(2007)}]{wang2007review}
\bibinfo{author}{G.~G. WANG}, \bibinfo{author}{S.~SHAN},
\newblock \bibinfo{title}{Review of metamodeling techniques in support of
  engineering design optimization},
\newblock \bibinfo{journal}{Journal of mechanical design (1990)}
  \bibinfo{volume}{129} (\bibinfo{year}{2007}) \bibinfo{pages}{370--380}.
%Type = Article
\bibitem[{Liu et~al.(2016)Liu, Bessa, and Liu}]{liu2016self}
\bibinfo{author}{Z.~Liu}, \bibinfo{author}{M.~Bessa}, \bibinfo{author}{W.~K.
  Liu},
\newblock \bibinfo{title}{Self-consistent clustering analysis: an efficient
  multi-scale scheme for inelastic heterogeneous materials},
\newblock \bibinfo{journal}{Computer Methods in Applied Mechanics and
  Engineering} \bibinfo{volume}{306} (\bibinfo{year}{2016})
  \bibinfo{pages}{319--341}.
%Type = Article
\bibitem[{Yu et~al.(2021)Yu, Kafka, and Liu}]{yu2021multiresolution}
\bibinfo{author}{C.~Yu}, \bibinfo{author}{O.~L. Kafka}, \bibinfo{author}{W.~K.
  Liu},
\newblock \bibinfo{title}{Multiresolution clustering analysis for efficient
  modeling of hierarchical material systems},
\newblock \bibinfo{journal}{Computational Mechanics} \bibinfo{volume}{67}
  (\bibinfo{year}{2021}) \bibinfo{pages}{1293--1306}.
%Type = Article
\bibitem[{Willcox and Peraire(2002)}]{willcox2002balanced}
\bibinfo{author}{K.~Willcox}, \bibinfo{author}{J.~Peraire},
\newblock \bibinfo{title}{Balanced model reduction via the proper orthogonal
  decomposition},
\newblock \bibinfo{journal}{AIAA journal} \bibinfo{volume}{40}
  (\bibinfo{year}{2002}) \bibinfo{pages}{2323--2330}.
%Type = Article
\bibitem[{Goury et~al.(2016)Goury, Amsallem, Bordas, Liu, and
  Kerfriden}]{goury2016automatised}
\bibinfo{author}{O.~Goury}, \bibinfo{author}{D.~Amsallem},
  \bibinfo{author}{S.~P.~A. Bordas}, \bibinfo{author}{W.~K. Liu},
  \bibinfo{author}{P.~Kerfriden},
\newblock \bibinfo{title}{Automatised selection of load paths to construct
  reduced-order models in computational damage micromechanics: from
  dissipation-driven random selection to bayesian optimization},
\newblock \bibinfo{journal}{Computational Mechanics} \bibinfo{volume}{58}
  (\bibinfo{year}{2016}) \bibinfo{pages}{213--234}.
%Type = Article
\bibitem[{Kerfriden et~al.(2011)Kerfriden, Gosselet, Adhikari, and
  Bordas}]{kerfriden2011bridging}
\bibinfo{author}{P.~Kerfriden}, \bibinfo{author}{P.~Gosselet},
  \bibinfo{author}{S.~Adhikari}, \bibinfo{author}{S.~P.-A. Bordas},
\newblock \bibinfo{title}{Bridging proper orthogonal decomposition methods and
  augmented newton--krylov algorithms: an adaptive model order reduction for
  highly nonlinear mechanical problems},
\newblock \bibinfo{journal}{Computer methods in applied mechanics and
  engineering} \bibinfo{volume}{200} (\bibinfo{year}{2011})
  \bibinfo{pages}{850--866}.
%Type = Article
\bibitem[{Lu et~al.(2018)Lu, Blal, and Gravouil}]{lu2018space}
\bibinfo{author}{Y.~Lu}, \bibinfo{author}{N.~Blal},
  \bibinfo{author}{A.~Gravouil},
\newblock \bibinfo{title}{Space--time pod based computational vademecums for
  parametric studies: application to thermo-mechanical problems},
\newblock \bibinfo{journal}{Advanced Modeling and Simulation in Engineering
  Sciences} \bibinfo{volume}{5} (\bibinfo{year}{2018}) \bibinfo{pages}{1--27}.
%Type = Article
\bibitem[{Ryckelynck(2009)}]{ryckelynck2009hyper}
\bibinfo{author}{D.~Ryckelynck},
\newblock \bibinfo{title}{Hyper-reduction of mechanical models involving
  internal variables},
\newblock \bibinfo{journal}{International Journal for Numerical Methods in
  Engineering} \bibinfo{volume}{77} (\bibinfo{year}{2009})
  \bibinfo{pages}{75--89}.
%Type = Article
\bibitem[{Carlberg et~al.(2013)Carlberg, Farhat, Cortial, and
  Amsallem}]{carlberg2013gnat}
\bibinfo{author}{K.~Carlberg}, \bibinfo{author}{C.~Farhat},
  \bibinfo{author}{J.~Cortial}, \bibinfo{author}{D.~Amsallem},
\newblock \bibinfo{title}{The gnat method for nonlinear model reduction:
  effective implementation and application to computational fluid dynamics and
  turbulent flows},
\newblock \bibinfo{journal}{Journal of Computational Physics}
  \bibinfo{volume}{242} (\bibinfo{year}{2013}) \bibinfo{pages}{623--647}.
%Type = Article
\bibitem[{Zhang et~al.(2017)Zhang, Combescure, and
  Gravouil}]{zhang2017efficient}
\bibinfo{author}{Y.~Zhang}, \bibinfo{author}{A.~Combescure},
  \bibinfo{author}{A.~Gravouil},
\newblock \bibinfo{title}{Efficient hyper-reduced-order model (hrom) for
  thermal analysis in the moving frame},
\newblock \bibinfo{journal}{International Journal for Numerical Methods in
  Engineering} \bibinfo{volume}{111} (\bibinfo{year}{2017})
  \bibinfo{pages}{176--200}.
%Type = Article
\bibitem[{Lu et~al.(2020)Lu, Jones, Gan, and Liu}]{lu2020adaptive}
\bibinfo{author}{Y.~Lu}, \bibinfo{author}{K.~K. Jones},
  \bibinfo{author}{Z.~Gan}, \bibinfo{author}{W.~K. Liu},
\newblock \bibinfo{title}{Adaptive hyper reduction for additive manufacturing
  thermal fluid analysis},
\newblock \bibinfo{journal}{Computer Methods in Applied Mechanics and
  Engineering} \bibinfo{volume}{372} (\bibinfo{year}{2020})
  \bibinfo{pages}{113312}.
%Type = Article
\bibitem[{Ladev{\`e}ze et~al.(2010)Ladev{\`e}ze, Passieux, and
  N{\'e}ron}]{ladeveze2010latin}
\bibinfo{author}{P.~Ladev{\`e}ze}, \bibinfo{author}{J.-C. Passieux},
  \bibinfo{author}{D.~N{\'e}ron},
\newblock \bibinfo{title}{The latin multiscale computational method and the
  proper generalized decomposition},
\newblock \bibinfo{journal}{Computer Methods in Applied Mechanics and
  Engineering} \bibinfo{volume}{199} (\bibinfo{year}{2010})
  \bibinfo{pages}{1287--1296}.
%Type = Article
\bibitem[{Boucinha et~al.(2014)Boucinha, Ammar, Gravouil, and
  Nouy}]{boucinha2014ideal}
\bibinfo{author}{L.~Boucinha}, \bibinfo{author}{A.~Ammar},
  \bibinfo{author}{A.~Gravouil}, \bibinfo{author}{A.~Nouy},
\newblock \bibinfo{title}{Ideal minimal residual-based proper generalized
  decomposition for non-symmetric multi-field models--application to transient
  elastodynamics in space-time domain},
\newblock \bibinfo{journal}{Computer Methods in Applied Mechanics and
  Engineering} \bibinfo{volume}{273} (\bibinfo{year}{2014})
  \bibinfo{pages}{56--76}.
%Type = Article
\bibitem[{Giacoma et~al.(2016)Giacoma, Dureisseix, and
  Gravouil}]{giacoma2016efficient}
\bibinfo{author}{A.~Giacoma}, \bibinfo{author}{D.~Dureisseix},
  \bibinfo{author}{A.~Gravouil},
\newblock \bibinfo{title}{An efficient quasi-optimal space-time pgd application
  to frictional contact mechanics},
\newblock \bibinfo{journal}{Advanced Modeling and Simulation in Engineering
  Sciences} \bibinfo{volume}{3} (\bibinfo{year}{2016}) \bibinfo{pages}{1--17}.
%Type = Article
\bibitem[{Ammar et~al.(2006)Ammar, Mokdad, Chinesta, and
  Keunings}]{ammar2006new}
\bibinfo{author}{A.~Ammar}, \bibinfo{author}{B.~Mokdad},
  \bibinfo{author}{F.~Chinesta}, \bibinfo{author}{R.~Keunings},
\newblock \bibinfo{title}{A new family of solvers for some classes of
  multidimensional partial differential equations encountered in kinetic theory
  modeling of complex fluids},
\newblock \bibinfo{journal}{Journal of non-Newtonian fluid Mechanics}
  \bibinfo{volume}{139} (\bibinfo{year}{2006}) \bibinfo{pages}{153--176}.
%Type = Article
\bibitem[{Chinesta et~al.(2011)Chinesta, Ladeveze, and
  Cueto}]{chinesta2011short}
\bibinfo{author}{F.~Chinesta}, \bibinfo{author}{P.~Ladeveze},
  \bibinfo{author}{E.~Cueto},
\newblock \bibinfo{title}{A short review on model order reduction based on
  proper generalized decomposition},
\newblock \bibinfo{journal}{Archives of Computational Methods in Engineering}
  \bibinfo{volume}{18} (\bibinfo{year}{2011}) \bibinfo{pages}{395--404}.
%Type = Article
\bibitem[{Chinesta et~al.(2013)Chinesta, Leygue, Bordeu, Aguado, Cueto,
  Gonz{\'a}lez, Alfaro, Ammar, and Huerta}]{chinesta2013pgd}
\bibinfo{author}{F.~Chinesta}, \bibinfo{author}{A.~Leygue},
  \bibinfo{author}{F.~Bordeu}, \bibinfo{author}{J.~V. Aguado},
  \bibinfo{author}{E.~Cueto}, \bibinfo{author}{D.~Gonz{\'a}lez},
  \bibinfo{author}{I.~Alfaro}, \bibinfo{author}{A.~Ammar},
  \bibinfo{author}{A.~Huerta},
\newblock \bibinfo{title}{Pgd-based computational vademecum for efficient
  design, optimization and control},
\newblock \bibinfo{journal}{Archives of Computational Methods in Engineering}
  \bibinfo{volume}{20} (\bibinfo{year}{2013}) \bibinfo{pages}{31--59}.
%Type = Article
\bibitem[{Relun et~al.(2013)Relun, N{\'e}ron, and Boucard}]{relun2013model}
\bibinfo{author}{N.~Relun}, \bibinfo{author}{D.~N{\'e}ron},
  \bibinfo{author}{P.-A. Boucard},
\newblock \bibinfo{title}{A model reduction technique based on the pgd for
  elastic-viscoplastic computational analysis},
\newblock \bibinfo{journal}{Computational Mechanics} \bibinfo{volume}{51}
  (\bibinfo{year}{2013}) \bibinfo{pages}{83--92}.
%Type = Article
\bibitem[{Heyberger et~al.(2012)Heyberger, Boucard, and
  N{\'e}ron}]{heyberger2012multiparametric}
\bibinfo{author}{C.~Heyberger}, \bibinfo{author}{P.-A. Boucard},
  \bibinfo{author}{D.~N{\'e}ron},
\newblock \bibinfo{title}{Multiparametric analysis within the proper
  generalized decomposition framework},
\newblock \bibinfo{journal}{Computational Mechanics} \bibinfo{volume}{49}
  (\bibinfo{year}{2012}) \bibinfo{pages}{277--289}.
%Type = Incollection
\bibitem[{Ladev{\`e}ze et~al.(2018)Ladev{\`e}ze, Paillet, and
  N{\'e}ron}]{ladeveze2018extended}
\bibinfo{author}{P.~Ladev{\`e}ze}, \bibinfo{author}{C.~Paillet},
  \bibinfo{author}{D.~N{\'e}ron},
\newblock \bibinfo{title}{Extended-pgd model reduction for nonlinear solid
  mechanics problems involving many parameters},
\newblock in: \bibinfo{booktitle}{Advances in Computational Plasticity},
  \bibinfo{publisher}{Springer}, \bibinfo{year}{2018}, pp.
  \bibinfo{pages}{201--220}.
%Type = Article
\bibitem[{Lu et~al.(2018)Lu, Blal, and Gravouil}]{lu2018multi}
\bibinfo{author}{Y.~Lu}, \bibinfo{author}{N.~Blal},
  \bibinfo{author}{A.~Gravouil},
\newblock \bibinfo{title}{Multi-parametric space-time computational vademecum
  for parametric studies},
\newblock \bibinfo{journal}{Finite Elements in Analysis and Design}
  \bibinfo{volume}{139} (\bibinfo{year}{2018}) \bibinfo{pages}{62--72}.
%Type = Article
\bibitem[{Modesto et~al.(2015)Modesto, Zlotnik, and Huerta}]{modesto2015proper}
\bibinfo{author}{D.~Modesto}, \bibinfo{author}{S.~Zlotnik},
  \bibinfo{author}{A.~Huerta},
\newblock \bibinfo{title}{Proper generalized decomposition for parameterized
  helmholtz problems in heterogeneous and unbounded domains: Application to
  harbor agitation},
\newblock \bibinfo{journal}{Computer Methods in Applied Mechanics and
  Engineering} \bibinfo{volume}{295} (\bibinfo{year}{2015})
  \bibinfo{pages}{127--149}.
%Type = Article
\bibitem[{Lu et~al.(2018)Lu, Blal, and Gravouil}]{lu2018adaptive}
\bibinfo{author}{Y.~Lu}, \bibinfo{author}{N.~Blal},
  \bibinfo{author}{A.~Gravouil},
\newblock \bibinfo{title}{Adaptive sparse grid based hopgd: Toward a
  nonintrusive strategy for constructing space-time welding computational
  vademecum},
\newblock \bibinfo{journal}{International Journal for Numerical Methods in
  Engineering} \bibinfo{volume}{114} (\bibinfo{year}{2018})
  \bibinfo{pages}{1438--1461}.
%Type = Article
\bibitem[{Lu et~al.(2019)Lu, Blal, and Gravouil}]{lu2019datadriven}
\bibinfo{author}{Y.~Lu}, \bibinfo{author}{N.~Blal},
  \bibinfo{author}{A.~Gravouil},
\newblock \bibinfo{title}{Datadriven hopgd based computational vademecum for
  welding parameter identification},
\newblock \bibinfo{journal}{Computational mechanics} \bibinfo{volume}{64}
  (\bibinfo{year}{2019}) \bibinfo{pages}{47--62}.
%Type = Article
\bibitem[{Hornik et~al.(1989)Hornik, Stinchcombe, and
  White}]{hornik1989multilayer}
\bibinfo{author}{K.~Hornik}, \bibinfo{author}{M.~Stinchcombe},
  \bibinfo{author}{H.~White},
\newblock \bibinfo{title}{Multilayer feedforward networks are universal
  approximators},
\newblock \bibinfo{journal}{Neural networks} \bibinfo{volume}{2}
  (\bibinfo{year}{1989}) \bibinfo{pages}{359--366}.
%Type = Article
\bibitem[{Lu et~al.(2019)Lu, Giovanis, Yvonnet, Papadopoulos, Detrez, and
  Bai}]{lu2019data}
\bibinfo{author}{X.~Lu}, \bibinfo{author}{D.~G. Giovanis},
  \bibinfo{author}{J.~Yvonnet}, \bibinfo{author}{V.~Papadopoulos},
  \bibinfo{author}{F.~Detrez}, \bibinfo{author}{J.~Bai},
\newblock \bibinfo{title}{A data-driven computational homogenization method
  based on neural networks for the nonlinear anisotropic electrical response of
  graphene/polymer nanocomposites},
\newblock \bibinfo{journal}{Computational Mechanics} \bibinfo{volume}{64}
  (\bibinfo{year}{2019}) \bibinfo{pages}{307--321}.
%Type = Article
\bibitem[{Raissi et~al.(2019)Raissi, Perdikaris, and
  Karniadakis}]{raissi2019physics}
\bibinfo{author}{M.~Raissi}, \bibinfo{author}{P.~Perdikaris},
  \bibinfo{author}{G.~E. Karniadakis},
\newblock \bibinfo{title}{Physics-informed neural networks: A deep learning
  framework for solving forward and inverse problems involving nonlinear
  partial differential equations},
\newblock \bibinfo{journal}{Journal of Computational Physics}
  \bibinfo{volume}{378} (\bibinfo{year}{2019}) \bibinfo{pages}{686--707}.
%Type = Article
\bibitem[{Liu et~al.(2019)Liu, Wu, and Koishi}]{liu2019deep}
\bibinfo{author}{Z.~Liu}, \bibinfo{author}{C.~Wu}, \bibinfo{author}{M.~Koishi},
\newblock \bibinfo{title}{A deep material network for multiscale topology
  learning and accelerated nonlinear modeling of heterogeneous materials},
\newblock \bibinfo{journal}{Computer Methods in Applied Mechanics and
  Engineering} \bibinfo{volume}{345} (\bibinfo{year}{2019})
  \bibinfo{pages}{1138--1168}.
%Type = Article
\bibitem[{Zhang et~al.(2021)Zhang, Cheng, Li, Gao, Yu, Domel, Yang, Tang, and
  Liu}]{zhang2021hierarchical}
\bibinfo{author}{L.~Zhang}, \bibinfo{author}{L.~Cheng},
  \bibinfo{author}{H.~Li}, \bibinfo{author}{J.~Gao}, \bibinfo{author}{C.~Yu},
  \bibinfo{author}{R.~Domel}, \bibinfo{author}{Y.~Yang},
  \bibinfo{author}{S.~Tang}, \bibinfo{author}{W.~K. Liu},
\newblock \bibinfo{title}{Hierarchical deep-learning neural networks: finite
  elements and beyond},
\newblock \bibinfo{journal}{Computational Mechanics} \bibinfo{volume}{67}
  (\bibinfo{year}{2021}) \bibinfo{pages}{207--230}.
%Type = Article
\bibitem[{Saha et~al.(2021)Saha, Gan, Cheng, Gao, Kafka, Xie, Li, Tajdari, Kim,
  and Liu}]{saha2021hierarchical}
\bibinfo{author}{S.~Saha}, \bibinfo{author}{Z.~Gan},
  \bibinfo{author}{L.~Cheng}, \bibinfo{author}{J.~Gao}, \bibinfo{author}{O.~L.
  Kafka}, \bibinfo{author}{X.~Xie}, \bibinfo{author}{H.~Li},
  \bibinfo{author}{M.~Tajdari}, \bibinfo{author}{H.~A. Kim},
  \bibinfo{author}{W.~K. Liu},
\newblock \bibinfo{title}{Hierarchical deep learning neural network (hidenn):
  An artificial intelligence (ai) framework for computational science and
  engineering},
\newblock \bibinfo{journal}{Computer Methods in Applied Mechanics and
  Engineering} \bibinfo{volume}{373} (\bibinfo{year}{2021})
  \bibinfo{pages}{113452}.
%Type = Article
\bibitem[{Zhang et~al.(2021)Zhang, Lu, Tang, and Liu}]{zhang2021hidenn}
\bibinfo{author}{L.~Zhang}, \bibinfo{author}{Y.~Lu}, \bibinfo{author}{S.~Tang},
  \bibinfo{author}{W.~K. Liu},
\newblock \bibinfo{title}{Hidenn-pgd: reduced-order hierarchical deep learning
  neural networks},
\newblock \bibinfo{journal}{arXiv preprint arXiv:2105.06363}
  (\bibinfo{year}{2021}).
%Type = Article
\bibitem[{Kolda and Bader(2009)}]{kolda2009tensor}
\bibinfo{author}{T.~G. Kolda}, \bibinfo{author}{B.~W. Bader},
\newblock \bibinfo{title}{Tensor decompositions and applications},
\newblock \bibinfo{journal}{SIAM review} \bibinfo{volume}{51}
  (\bibinfo{year}{2009}) \bibinfo{pages}{455--500}.
%Type = Article
\bibitem[{Tucker(1966)}]{tucker1966some}
\bibinfo{author}{L.~R. Tucker},
\newblock \bibinfo{title}{Some mathematical notes on three-mode factor
  analysis},
\newblock \bibinfo{journal}{Psychometrika} \bibinfo{volume}{31}
  (\bibinfo{year}{1966}) \bibinfo{pages}{279--311}.
%Type = Article
\bibitem[{Kiers(2000)}]{kiers2000towards}
\bibinfo{author}{H.~A. Kiers},
\newblock \bibinfo{title}{Towards a standardized notation and terminology in
  multiway analysis},
\newblock \bibinfo{journal}{Journal of Chemometrics: A Journal of the
  Chemometrics Society} \bibinfo{volume}{14} (\bibinfo{year}{2000})
  \bibinfo{pages}{105--122}.
%Type = Article
\bibitem[{Denlinger(2015)}]{denlinger2015thermo}
\bibinfo{author}{E.~R. Denlinger},
\newblock \bibinfo{title}{Thermo-mechanical model development and experimental
  validation for metallic parts in additive manufacturing},
\newblock \bibinfo{journal}{Thesis}  (\bibinfo{year}{2015}).
%Type = Article
\bibitem[{Moulinec and Suquet(1998)}]{moulinec1998numerical}
\bibinfo{author}{H.~Moulinec}, \bibinfo{author}{P.~Suquet},
\newblock \bibinfo{title}{A numerical method for computing the overall response
  of nonlinear composites with complex microstructure},
\newblock \bibinfo{journal}{Computer methods in applied mechanics and
  engineering} \bibinfo{volume}{157} (\bibinfo{year}{1998})
  \bibinfo{pages}{69--94}.
%Type = Article
\bibitem[{Ten~Berge et~al.(1987)Ten~Berge, De~Leeuw, and
  Kroonenberg}]{ten1987some}
\bibinfo{author}{J.~M. Ten~Berge}, \bibinfo{author}{J.~De~Leeuw},
  \bibinfo{author}{P.~M. Kroonenberg},
\newblock \bibinfo{title}{Some additional results on principal components
  analysis of three-mode data by means of alternating least squares
  algorithms},
\newblock \bibinfo{journal}{Psychometrika} \bibinfo{volume}{52}
  (\bibinfo{year}{1987}) \bibinfo{pages}{183--191}.
%Type = Article
\bibitem[{Zhang and Golub(2001)}]{zhang2001rank}
\bibinfo{author}{T.~Zhang}, \bibinfo{author}{G.~H. Golub},
\newblock \bibinfo{title}{Rank-one approximation to high order tensors},
\newblock \bibinfo{journal}{SIAM Journal on Matrix Analysis and Applications}
  \bibinfo{volume}{23} (\bibinfo{year}{2001}) \bibinfo{pages}{534--550}.

\end{thebibliography}






% \end{thebibliography}

\end{document}